%% file: manuscript_Arxiv.tex
\documentclass[draft]{agujournal2019}
\usepackage{url} %this package should fix any errors with URLs in refs.
\usepackage{lineno}
\usepackage[inline]{trackchanges} %for better track changes. finalnew option will compile document with changes incorporated.
\usepackage{soul}
\include{macros}

\draftfalse

\journalname{ResearchGate}

\begin{document}

%%%%%%%%%%%%%%%%%%%%%%%%%%%%%%%%%%%%%%%%%%%%%%%
%  TITLE
%
% (A title should be specific, informative, and brief. Use
% abbreviations only if they are defined in the abstract. Titles that
% start with general keywords then specific terms are optimized in
% searches)
%
%%%%%%%%%%%%%%%%%%%%%%%%%%%%%%%%%%%%%%%%%%%%%%%

% Example: \title{This is a test title}

\title{Multiscale Modeling Framework using Element-based Galerkin Methods for Moist Atmospheric Limited-Area Simulations}

%%%%%%%%%%%%%%%%%%%%%%%%%%%%%%%%%%%%%%%%%%%%%%%
%
%  AUTHORS AND AFFILIATIONS
%
%%%%%%%%%%%%%%%%%%%%%%%%%%%%%%%%%%%%%%%%%%%%%%%

% List authors by first name or initial followed by last name and
% separated by commas. Use \affil{} to number affiliations, and
% \thanks{} for author notes.
% Additional author notes should be indicated with \thanks{} (for
% example, for current addresses).

% Example: \authors{A. B. Author\affil{1}\thanks{Current address, Antartica}, B. C. Author\affil{2,3}, and D. E.
% Author\affil{3,4}\thanks{Also funded by Monsanto.}}

\authors{Soonpil Kang\affil{1}, James F. Kelly\affil{2}, Anthony P. Austin\affil{1}, Francis X. Giraldo\affil{1}}

\affiliation{1}{Department of Applied Mathematics, Naval Postgraduate School, U.S.}
\affiliation{2}{Space Science Division, U.S. Naval Research Laboratory, Washington, D.C., U.S.}

% Corresponding author mailing address and e-mail address:

\correspondingauthor{Soonpil Kang}{soonpil.kang.ks@nps.edu}

%%%%%%%%%%%%%%%%%%%%%%%%%%%%%%%%%%%%%%%%%%%%%%%
% KEY POINTS
%%%%%%%%%%%%%%%%%%%%%%%%%%%%%%%%%%%%%%%%%%%%%%%

\begin{keypoints}
\item A multiscale modeling framework (MMF) for moist atmospheric limited-area simulations is developed.
\item The large-scale and small-scale processes in the framework are modeled using a nonhydrostatic model and an element-based Galerkin method.
\item The MMF improves the representation of moist dynamics for cloud processes.
\end{keypoints}

%% \begin{abstract} starts the second page

\begin{abstract}

This paper presents a multiscale modeling framework (MMF) to model moist atmospheric limited-area weather. The MMF resolves large-scale convection using a coarse grid while simultaneously resolving local features through numerous fine local grids and coupling them seamlessly. Both large- and small-scale processes are modeled using the compressible Navier-Stokes equations within the Nonhydrostatic Unified Model of the Atmosphere (NUMA), and they are discretized using a continuous element-based Galerkin method (spectral elements) with high-order basis functions. Consequently, the large-scale and small-scale models share the same dynamical core but have the flexibility to be adjusted individually. The proposed MMF method is tested in 2D and 3D idealized limited-area weather problems involving storm clouds produced by squall line and supercell simulations. The MMF numerical results showed enhanced representation of cloud processes compared to the coarse model.

\end{abstract}

\section*{Plain Language Summary}

Convential Numerical Weather Models (NWPs) represent convective processes using a subgrid-scale parameterization. The large-scale convective processes typically use grid with relatively lower resolution, while resolving small-scale cloud processes requires much higher resolution to accurately predict the weather. To efficiently address this multiscale problem, the multiscale modeling framework (MMF) simulates the large-scale process on a lower resolution grid, while incorporating multiple higher resolution grids to simulate local cloud processes. We present a multiscale modeling framework in which the large-scale and small-scale models are constructed using similar mathematical models and high-order discretization methods. 

%%%%%%%%%%%%%%%%%%%%%%%%%%%%%%%%%%%%%%%%%%%%%%%
%
%  BODY TEXT
%
%%%%%%%%%%%%%%%%%%%%%%%%%%%%%%%%%%%%%%%%%%%%%%%

% Introduction
\input{sections/Introduction}

% Numerical Method
\input{sections/Governing_Equations}
\input{sections/MMF}
\input{sections/Numerical_Method.tex}
\input{sections/EBG.tex}
\input{sections/SGS.tex}

\input{sections/Sponge.tex}

\input{sections/Nonconforming_Grids.tex}
\input{sections/Implementation.tex}

% Complexity
\input{sections/Complexity.tex}

% Numerical Results
\input{sections/Numerical_Results}
\input{sections/Squall_Line}

\input{sections/Supercell}

% Conclusions
\input{sections/Conclusions}

%%%%%%%%%%%%%%%%%%%%%%%%%%%%%%%%%%%%%%%%%%%%%%%
%
% DATA SECTION and ACKNOWLEDGMENTS
%
%%%%%%%%%%%%%%%%%%%%%%%%%%%%%%%%%%%%%%%%%%%%%%%

\section*{Data Availability Statement}
Figures of instantaneous fields in the squall line test were made using Matlab. Figures of instantaneous fields in the supercell test were made using Paraview. Other figures were produced using Matplotlib. The figures and data for squall line and supercell simulations associated with this manuscript are avaliable on Zenodo repository (\url{https://doi.org/10.5281/zenodo.11166368}).

\acknowledgments
This work was supported by the Office of Naval Research under Grant No.\
N0001419WX00721. F.\ X.\ Giraldo was also supported by the National Science Foundation under grant AGS-1835881. This work was performed when Soonpil Kang held a National Academy of Sciences’ National Research Council Fellowship at the Naval Postgraduate School. The authors  wish to thank Wojciech W.\ Grabowski (National Center for Atmospheric Research) and Walter Hannah (Lawrence Livermore National Laboratory) for helpful discussions.

%%%%%%%%%%%%%%%%%%%%%%%%%%%%%%%%%%%%%%%%%%%%%%%
% REFERENCES and BIBLIOGRAPHY
%%%%%%%%%%%%%%%%%%%%%%%%%%%%%%%%%%%%%%%%%%%%%%%

\bibliography{references}

\end{document}

%% file: macros.tex
\usepackage{amsmath}
\usepackage{amssymb}
\usepackage{algorithm, algorithmic}
\usepackage{bm}
\usepackage{graphicx}
\usepackage{physics}
\usepackage{subcaption}
\allowdisplaybreaks % Page break of equations

\newcommand{\vc}[1]{\mbox{\boldmath$#1$\unboldmath}}
\newcommand{\dt}{\Delta t}
\newcommand{\DT}{\Delta T}
\newcommand{\dx}{\Delta x}
\newcommand{\dy}{\Delta y}
\newcommand{\dz}{\Delta z}
\newcommand{\fb}{\vc{f}}
\newcommand{\Fb}{\vc{F}}
\newcommand{\kb}{\vc{k}}

\newcommand{\Mhat}{\widehat{M}}
\newcommand{\Mtilde}{\widetilde{M}}
\newcommand{\qb}{\vc{q}}

\newcommand{\psihat}{\widehat{\psi}}
\newcommand{\psitilde}{\widetilde{\psi}}
\newcommand{\Qb}{\vc{Q}}

\newcommand{\ub}{\vc{u}}

\newcommand{\xb}{\vc{x}}
\newcommand{\zhat}{\widehat{\zeta}}
\newcommand{\ztilde}{\widetilde{\zeta}}

\newcommand{\genericdel}[4]{%
  \ifcase#3\relax
  \ifx#1.\else#1\fi#4\ifx#2.\else#2\fi\or
  \bigl#1#4\bigr#2\or
  \Bigl#1#4\Bigr#2\or
  \biggl#1#4\biggr#2\or
  \Biggl#1#4\Biggr#2\else
  \left#1#4\right#2\fi
}
\newcommand{\evalx}[2][-1]{\genericdel.|{#1}{#2}}

\newcommand{\intO}[1]{\int_{\Omega}#1d\Omega}
\newcommand{\intOe}[1]{\int_{\Omega_e}#1d\Omega_e}

\renewcommand\arraystretch{1.3}

%% file: sections/Introduction.tex
\section{Introduction}

Atmospheric flows, within the context of weather and climate, encompass a wide spectrum of spatial and temporal scales. These scales range from thousands of kilometers for general circulation to hundreds of meters for precipitation microphysics. The multiscale nature of these phenomena poses a significant computational challenge for accurately simulating the atmosphere. To address this challenge, deep moist convection can be represented by conventional subgrid-scale parameterizations or by simulating them explicitly. General circulation models (GCMs) for the moist atmosphere typically resolve the large-scale process (LSP), while parameterizing unresolved information from the small-scale processes (SSPs). Conventional parameterizations contain modeling approximations that are difficult to quantify and consequently require significant efforts for tuning free parameters. For an operational weather model, the calibration process needs to be performed for each numerical setup and different suit of parameterizations. In contrast, the multiscale modeling framework (MMF) aims to explicitly resolve both the LSP and SSP on different grids and couples them through their interaction. The LSP represents the general convective process, while the SSP represents the cloud-scale process. Within this framework, the LSP model forces the SSP model and receives feedback from it. Since the SSP model is based on the dynamics equations, the modeling bias is minimized, unlike traditional parameterizations.

% SP & MMF

MMF methods have been proposed for weather and climate simulations to replace traditional semi-empirical parameterizations used by cloud-resolving models (CRMs). Grabowski and Smolarkiewicz \cite{grabowski1999crcp} proposed a computational framework termed cloud resolving convection parameterization (CRCP) and later referred to it as superparameterization \cite{grabowski2004improved}. This framework is designed to capture explicitly both large-scale and cloud-scale dynamics and couple them through relaxation and has been shown to improve the representation of small- and mesoscale processes. Khairoutdinov and Randall \cite{khairoutdinov2001cloud, khairoutdinov2005simulations} developed a superparameterized version of the NCAR Community Atmosphere Model (SP-CAM). It was shown that SP-CAM improves the variability of precipitation over the land and ocean, closely matching the observed precipitation frequency. Recent investigations on the performance of SP-CAM variants with a particular focus on rain statistics are available in the literature \cite{kooperman2016robust}. This improvement of MMF over conventional parameterizations for representing precipitation was also demonstrated in \cite{pritchard2009empirical}. The Department of Energy Superparameterized Energy Exascale Earth System Model (SP-E3SM) \cite{hannah2020initial, lin2022mesoscale} combines a spectral-element based dynamical core with the System for Atmospheric Modeling (SAM) \cite{khairoutdinov2003cloud} as a CRM. SP-E3SM produces more accurate mesoscale convective systems and realistic tropical waves compared to its standard version. This work has recently been extended to incorporate land surface processes, such as soil moisture, into their CRM \cite{lin2023modeling}, which is crucial for global climate modeling. Majda and Grooms \cite{majda2014new} presented a formulation of stochastic superparameterization wherein an eddy closure is derived from the stochastic modeling of the SSPs. The effectiveness of this formulation is examined on idealized mathematical test cases. Han et al. \cite{han2020moist} applied a deep convolutional residual neural network technique to emulate SP-CAM model for moist physics parameterization in GCMs.

These MMF works commonly employ different dynamics models (i.e., equation sets) for the LSP and SSP. In SSP modeling, the anelastic approximation is widely adopted to avoid numerical issues associated with fast-propagating acoustic waves \cite{grabowski2016towards}; however, Arakawa and Konor \cite{arakawa2009unification} emphasized the necessity of a unified modeling approach for GCM and CRM to ensure consistency and provide flexibility in choosing resolution based on the objective of the applications. Many MMF works have incorporated 2D CRMs for SSP modeling, motivated by reasonable success of early MMF works. While 2D CRMs are considered capable of representing cloud processes in each column of the large-scale model, as documented in \cite{grabowski2001coupling}, there are efforts to replicate 3D effects of SSP or directly implement 3D SSP models. Jung and Arakawa \cite{jung2010development} proposed a quasi-3D CRM to couple with a GCM, where the SSP domains are narrow channels with a small thickness and embedded in perpendicular directions, achieving comparable results to a full 3D CRM. Jansson et al. \cite{jansson2019regional} presented a coupling scheme for OpenIFS and 3D DALES models in selective regions, which were successfully integrated despite significant differences in terms of domain height, units for quantities, and prognostic variables.

This paper presents a novel MMF method designed for moist atmospheric limited-area simulations. The LSP and SSP models within the current MMF share the same dynamical core, employing consistent governing equations and discretization methods. Both models are constructed using the nonhydrostatic compressible Navier-Stokes equations with moisture and discretized by the element-based Galerkin method (e.g., see \cite{giraldo_introduction_2020}). We use the Nonhydrostatic Unified Model of the Atmosphere (NUMA) \cite{kelly2012continuous} as a base model for both LSP and SSP. To implement this MMF simulation, we have restructured the NUMA code into a fully object-oriented form, enabling the instantiation of numerous simulations (we call this new model xNUMA). We test the MMF method by applying it to two idealized limited-area weather problems: 1) a 2D squall line \cite{gabervsek2012dry} and 2) a 3D supercell thunderstorm \cite{tissaoui2023non}.

The element-based Galerkin method, specifically the spectral element method, has been adopted for the dynamical core of limited area models (LAMs) and GCMs for numerical weather prediction, as listed in \cite{marras2016review}. Among the MMF models, SP-E3SM \cite{hannah2020initial} utilizes the spectral element method for its GCM. The element-based Galerkin method enables the adjustment of spatial resolution via element size ($h$-refinement) and the order of basis functions ($p$-refinement). It has demonstrated high-order accuracy and efficiency across various applications of dry and moist atmosphere \cite{yi2020vertical, gabervsek2012dry, tissaoui2023non}. This paper represents the first application of the element-based Galerkin method for both LSP and SSP models for limited-area atmospheric MMF.

% \rc{FXG: Nice! This is a great statement showing the novelty of this work.}

% Outline

An outline of the paper is as follows. The governing equations for the nonhydrostatic moist atmosphere are presented in Section \ref{sec:governing}. Section \ref{sec:mmf} presents the coupling algorithm for the LSP and SSP model. Section \ref{sec:method} details the numerical methods including the spatial and temporal discretizations, sponge layer, nonconforming vertical grids, and parallel implementation. Section \ref{sec:complexity} evaluates the complexity of the algorithms for the standard and MMF simulations in terms of arithmetic intensity. The numerical test cases for the 2D squall line and 3D supercell are investigated in Section \ref{sec:results}, and conclusions are drawn in Section \ref{sec:conclusion}.

%% file: sections/Governing_Equations.tex
\section{Governing equations}\label{sec:governing}

Let us consider a fixed spatial domain $\Omega$ and a time interval $(0,t_f]$. The governing equations for the nonhydrostatic model of the moist atmosphere can be written as follows:
\begin{subequations}\label{eq:governing}
	\begin{align}
		&\pdv{\rho'}{t} + \div\qty\big((\rho_0+\rho')\ub) = 0, \label{eq:governing:rho}
		\\
		&\pdv{\ub}{t} + \ub\cdot\grad\ub + \dfrac{1}{\rho_0+\rho'} \grad p' + g\kb \qty(\dfrac{\rho'}{\rho_0+\rho'} - \varepsilon q_v'+q_c+q_r) = \nu\grad^2\ub, \label{eq:governing:u}
		\\
		&\pdv{\theta_v'}{t} + \ub\cdot\grad(\theta_{v0}+\theta_v') = S_{\theta_v'} + \nu\grad^2\theta_v', \label{eq:governing:theta}
		\\
		&\pdv{q_v'}{t} + \ub\cdot\grad (q_{v0} + q_v') = S_v + \nu\grad^2 q_v', \label{eq:governing:qv}\\
		&\pdv{q_c}{t} + \ub\cdot\grad q_c = S_c + \nu\grad^2 q_c, \label{eq:governing:qc}\\
		&\pdv{q_r}{t} + \ub\cdot\grad q_r = S_r + \nu\grad^2 q_r,	\label{eq:governing:qr}
	\end{align}
\end{subequations}
where $\rho_0$ is the reference density, $\rho'$ is the density perturbation, $\ub$ is the velocity vector, $\theta_{v0}$ is the reference virtual potential temperature, $\theta_v'$ represents the virtual potential temperature perturbation, $p$ is the pressure, $g$ is the gravitational acceleration, $\kb$ is the upward-pointing unit vector, $\varepsilon=R_v/R_d-1\approx 0.608$ with the specific gas constants of dry air and water vapor denoted by $R_d$ and $R_v$, and $\nu$ is the kinematic viscosity.
In the last three equations, $q_v$, $q_c$, and $q_r$ represent the mixing ratios of water vapor, cloud water, and rain, respectively. Equations \eqref{eq:governing:rho}--\eqref{eq:governing:theta} represent local mass conservation, momentum balance, and thermodynamics, respectively. The transport of water species and moist microphysical processes are represented by Eqs.\ \eqref{eq:governing:qv}-\eqref{eq:governing:qr}, the latent heat release term $S_{\theta_v'}$ in Eq.\ \eqref{eq:governing:theta}, and the buoyancy terms in Eq.\ \eqref{eq:governing:u}. We assume that the reference fields are in hydrostatic balance and dependent only on the vertical coordinate $z$, i.e., $\dv{p_0}{z}=-\rho_0 g$. The source terms, $S_{\theta_v}$, $S_v$, $S_q$, and $S_r$ describe the effects of the microphysical processes, including condensation, autoconversion, accretion, evaporation, and sedimentation. These terms are defined using the Kessler microphysics model, as summarized in \cite{klemp1978simulation}.

The prognostic variables for the system \eqref{eq:governing} are $\qb=(\rho', \ub^T, \theta_v', q_v', q_c, q_r)^\mathcal{T}$, where the superscript $\mathcal{T}$ denotes the transpose operator.  Eq.\ \eqref{eq:governing} can be rewritten more concisely in the compact form as
\begin{equation}\label{eq:governing:compact}
	\pdv{\qb}{t} = \mathcal{S}(\qb),
\end{equation}
where $\mathcal{S}(\qb)$ contains all terms in the equations apart from the time derivatives. The pressure of moist air is calculated using the thermodynamic equation of state,
\begin{equation}
  p = \rho R_d T_v,
\end{equation}
where the virtual temperature $T_v$ is defined as $T_v = T \qty(1+\varepsilon q_v)$, with the air temperature $T$ \cite{vallis2005atmospheric}.

Let us now describe the MMF strategy that we use to solve the governing equations.

%% file: sections/MMF.tex
\section{Multiscale modeling framework}\label{sec:mmf}

\begin{figure}
  \centering
  \includegraphics[width=0.6\textwidth]{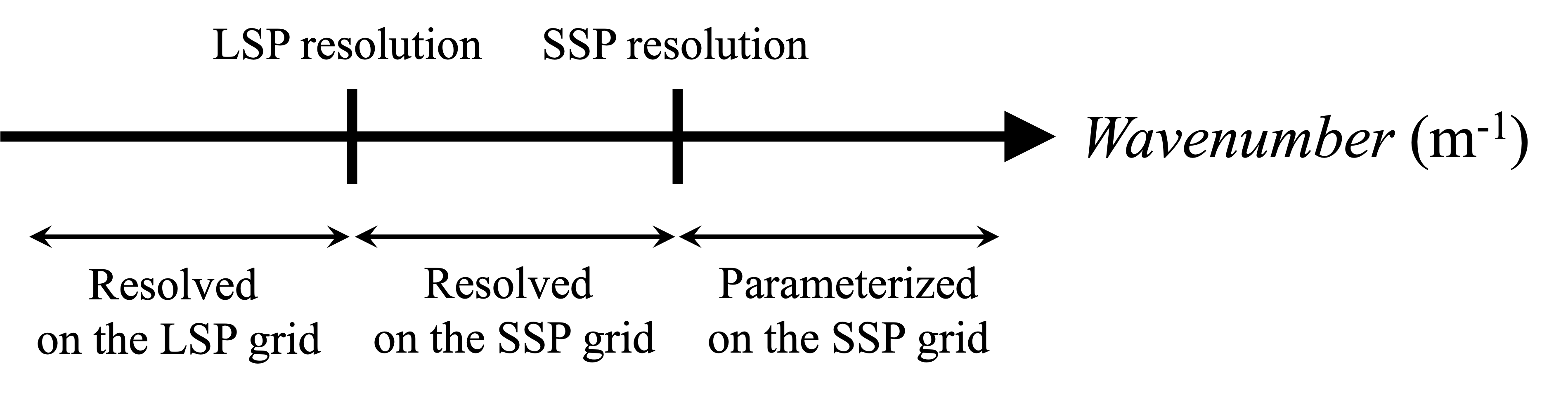}
  \caption{Decomposition of the scales modeled via MMF along the wavenumber axis.}\label{fig:scale_decomposition}
\end{figure}

We consider a limited-area domain for modeling the large-scale process (LSP) and local domains for modeling the small-scale processes (SSP) locally within each column of the LSP. Figure \ref{fig:scale_decomposition} illustrates the three regions of scales along the conceptual wavenumber axis that are handled differently in our MMF strategy. The LSP (such as the global circulation) is resolved on the LSP grid while the SSPs (such as cloud-scale processes) are explicitly resolved on the SSP grid. Subgrid-scale processes including precipitation microphysics are parameterized on the SSP grid. Since the parameterization is performed on a finer grid with the SSP resolution than in the standard simulation, it is expected that this will enhance accuracy in representing fine-scale features.

%------------------------------------------------------%
%									 LSP and SSP coupling
%------------------------------------------------------%
\subsection{LSP and SSP coupling}

\begin{figure}
  \centering
  \includegraphics[width=0.8\textwidth]{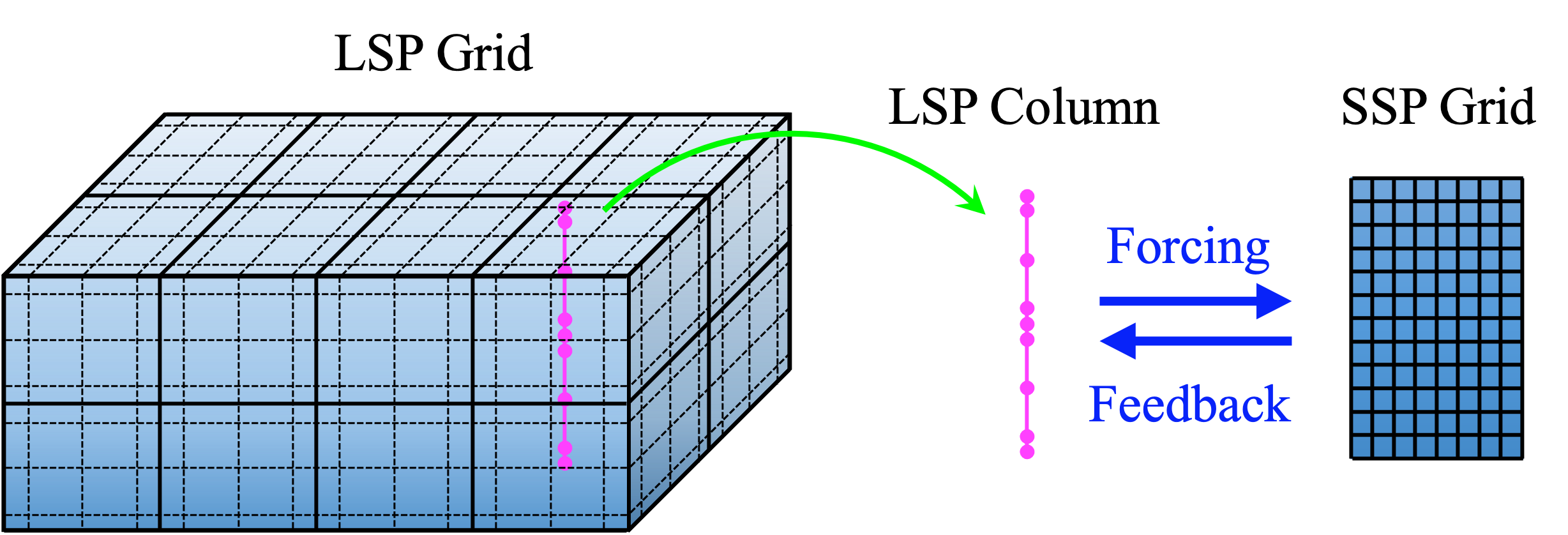}
  \caption{Configuration of the LSP and SSP grids for MMF.}\label{fig:MMF_config}
\end{figure}

In order to couple the LSP and SSP models, we employ the coupling method used in the original superparameterization scheme \cite{grabowski2001coupling, grabowski2004improved}. Let us denote the LSP and SSP domains as $\Omega_{LSP}$ and $\Omega_{SSP}$, with the prognostic variables for LSP and SSP represented as $\Qb$ and $\qb$, respectively. The LSP domain is discretized with low resolution as a global circulation model (GCM), while the SSP domains are discretized with high resolution as a cloud-resolving model (CRM), as shown in Figure \ref{fig:MMF_config}. In this work, we consider a two-dimensional SSP model generated at each grid column of a 2D or 3D LSP model. For three-dimensional problems, the SSP domains align with the background wind direction, usually the $x$ direction in weather applications. The coupling condition that is imposed is
\begin{equation}
  \label{eq:coupling_condition}
  \Qb\qty(X,Y,Z,t) = \evalx{\expval{\qb(x,y,z,t)}}_{(X,Y)},
\end{equation}
which dictates that the LSP variable equals the horizontal average of the SSP variable for each level in the LSP column. We assume periodicity in the horizontal direction \cite{grabowski2001coupling}. Therefore, this SSP model represents the small-scale processes that occur around each LSP column. The two-dimensional horizontal averaging operator $\evalx{\expval{\cdot}}_{(X,Y)}$ is defined as
\begin{equation}
  \evalx{\expval{\qb(x,0,z,t)}}_{(X,Y)} = \dfrac{1}{L_x}\int^{L_x/2}_{-L_x/2} \qb(\xi,0,z,t)|_{(X,Y)} d\xi,
\end{equation}
where $L_x$ is the length of the SSP domain in the $x$ direction.

The LSP and SSP models are coupled in two ways. LSP variables drive the flows in the SSP models through a forcing term, and each LSP column receives feedback from the corresponding SSP models. The coupling of the LSP and SSP models is achieved by adding forcing or feedback terms to the governing equations \eqref{eq:governing:compact}, resulting in the following coupled model:
\begin{alignat}{2}
  \pdv{\Qb}{t} &= \mathcal{S}(\Qb) + \Fb(\Qb,\qb) \label{eq:govering-LSP} \\
  \pdv{\qb}{t} &= \mathcal{S}(\qb) + \fb(\qb,\Qb), \label{eq:govering-SSP}
\end{alignat}
where $\mathcal{S}(\Qb)$ and $\mathcal{S}(\qb)$ are the operators in the governing equations in \eqref{eq:governing} in terms of the LSP and SSP variables. We use the same governing equations for modeling both LSP and SSP. Consequently, the MMF proposed in this work yields a nonhydrostatic model for both LSP and SSP and a mathematically consistent method in the treatment of the large-scale and small-scale features. 

%------------------------------------------------------%
%							 Temporal discretization
%------------------------------------------------------%
\subsection{Temporal discretization}

The update of LSP and SSP variables progresses from time level $n$ to $n+1$ over the LSP time step $\DT$ using an implicit-explicit (IMEX) time integrator \cite{giraldo2013implicit} as follows:
\begin{alignat}{2}
  \Qb^{n+1} &= \Qb^{n} + \DT \qty\Big[ \mathcal{S}_{\delta}\qty(\Qb^n,\Qb^{n+1}) + \Fb^n(\Qb^n,\qb^n) ], \\
  \qb^{n+1} &= \qb^{n} + \sum_{m=1}^{M} \dt\qty\Big[\mathcal{S}_{\delta}(\qb^{n+(m-1)/M},\qb^{n+m/M}) + \fb^n(\qb^n,\Qb^{n+1}) ], 
\end{alignat}
where $\dt$ is the SSP time step, $M=\DT/\dt$ is the number of SSP time steps per LSP time step, and $\mathcal{S}_{\delta}$ represents the IMEX operator defined as $\mathcal{S}_{\delta} \qty(\Qb^n, \Qb^{n+1} ) = \qty\big[ \mathcal{S} \qty(\Qb^n ) - \delta \mathcal{L} \qty(\Qb^n ) ] + \delta \mathcal{L} \qty(\Qb^{n+1} )$ with the linear operator $\mathcal{L}$ given in \cite{giraldo2013implicit} and the parameter $\delta\in\qty{0,1}$ switching the formulation between forward Euler and linear IMEX. The forcing term $\Fb^n$ and the feedback term $\fb^n$ are defined as follows:
\begin{align}
  \Fb^n(\Qb^n,\qb^n) &= \dfrac{\expval{\qb^{n}}-\Qb^{n}}{\DT}, \label{eq:forcing}\\
  \fb^n(\qb^n,\Qb^{n+1}) &= \dfrac{\Qb^{n+1}-\expval{\qb^{n}}}{\DT},\label{eq:feedback}
\end{align}
where $\qb^n$ denotes the SSP state variables from the previous time step $n$, while $\Qb^n$ and $\Qb^{n+1}$ are the LSP state variables from the previous time step $n$ and the current time step $n+1$, respectively. Eqs.\ \eqref{eq:forcing}--\eqref{eq:feedback} introduce tendencies that relax each LSP state variable toward their SSP counterparts and vice versa, ensuring that both do not drift appreciably during time integration. Both \eqref{eq:forcing} and \eqref{eq:feedback} will vanish if the coupling condition in Eq.\ \eqref{eq:coupling_condition} is exactly satisfied at each time step.

In the IMEX method, the terms responsible for fast waves, such as acoustic and gravity waves, are discretized implicitly, while the other terms for slower dynamics are discretized explicitly. This IMEX strategy enables the use of larger time steps by eliminating the tight CFL restriction due to the fast waves. In this work, we employ the second-order additive Runge-Kutta (ARK2) scheme proposed in \cite{giraldo2013implicit}, which belongs to a family of linear multistage schemes. This time integrator is also used in \cite{giraldo2023performance} by referring to ARK(2,3,2)b. We solve the linear system that arises from the implicit part of the IMEX formulation using the iterative Krylov subspace GMRES \cite[pg. 266]{trefethen1997} in a matrix-free fashion. While the MMF method offers flexibility in selecting different time integrators for LSP and SSP models, we choose to use the same time integrators for both the models in MMF simulations in this work. In future work, we will explore the use of different time integrators and orders of accuracy for the SSP models. 

%------------------------------------------------------%
%							 Implementation details
%------------------------------------------------------%
\subsection{Implementation details}

We couple the LSP and SSP models in terms of horizontal velocity, virtual potential temperature, and the mixing ratios of water vapor, cloud water, and rain; density and vertical velocity are not coupled. The idea here is to mimic as closely as possible the procedure used with a standard parameterization. We consider grid refinement primarily in the horizontal direction via the MMF but also in the vertical direction using the nonconforming vertical discretization described in Sec.\ \ref{sec:nonconforming}. Consequently, the forcing and feedback terms in Eqs.\ \eqref{eq:forcing}--\eqref{eq:feedback} are added to the right-hand side of the corresponding balance equations in \eqref{eq:governing}.

% \rc{(SK: We may need to explain about not coupling the vertical velocity. This is straightforward for the anelastic SSPs with periodicity, but not clear for the compressible SSPs.)}

To ensure an adequate representation of spatial scales around the LSP resolution, as shown in Figure \ref{fig:scale_decomposition}, it is necessary to select an SSP domain size large enough to capture the largest wavenumber (smallest wavelength) that can be resolved at the LSP resolution.

Algorithm \ref{algorithm:time_loop} presents pseudocode for the time integration loop in the MMF. The outermost loop iterates over time steps for the LSP model. During each time step, we begin by updating the LSP state variables $\Qb$ by integrating the equation in \eqref{eq:govering-LSP}. Afterward, we update the SSP state variables $\qb$ within each SSP domain through multiple sub-time steps. A smaller time-step size may be necessary for the SSP due to its finer grid resolution and CFL constraint. In principle, we have flexibility to select different time-integration methods to independently handle the stability and efficiency of the LSP and SSP models. The numerical test case in Sec.\ \ref{section:squall} demonstrates that these staggered updates impose the coupling conditions accurately. This staggered approach also reduces computational overhead. In distributed-memory parallel computations, inter-processor communications are required only during the LSP update, while the subsequent SSP updates occur on-process. 

\begin{algorithm}
	\caption{Time integration loop. ($nstep$ is the number of time steps in the LSP model, $M$ is the number of SSP time steps per LSP time step, and the coefficient $\delta\in\qty{0,1}$ determines whether the time-integration method is explicit or IMEX.)}\label{algorithm:time_loop}
  \begin{algorithmic}
    \FOR{$n=1,nstep$}
    \STATE{Update LSP state variables: $\Qb^{n+1} = \Qb^{n} + \DT \qty[ \mathcal{S}_{\delta}\qty(\Qb^n,\Qb^{n+1}) + \dfrac{\expval{\qb^n}-\Qb^n}{\DT}]$. }
      \FOR{$m=1,M$}
      \STATE{Update SSP state variables: \\
        $\qb^{n+m/M} = \qb^{n+(m-1)/M} + \dt \qty[ \mathcal{S}_{\delta}\qty(\qb^{n+(m-1)/M},\qb^{n+m/M}) + \dfrac{\Qb^{n+1}-\expval{\qb^n}}{\DT} ]$. }
      \STATE{Update SSP state variables by the contribution of the source term.}
      \ENDFOR
    \ENDFOR
  \end{algorithmic}
\end{algorithm}

%% file: sections/Numerical_Method.tex
\section{Numerical methods}\label{sec:method}

This section describes the numerical methods employed in our limited-area simulations using the standard and MMF approaches. In the MMF, we use the continuous version of the element-based Galerkin method (spectral elements) \cite{giraldo_introduction_2020} to spatially descretize both the LSP and SSP domains. We add the microphysics model to parameterize the SSP models. An implicit sponge layer \cite{klemp2008upper} serves as an upper absorbing boundary condition (ABC). The coupling of the LSP and SSP grids accommodates nonconforming discretization in the vertical direction. Sec.\ \ref{sec:implementation} describes the parallel implementation of our MMF simulations.

%% file: sections/EBG.tex
\subsection{Element-based Galerkin method}

The element-based Galerkin method decomposes the spatial domain $\Omega\subset\mathbb{R}^d$ ($d=$2 or 3) into $N_e$ disjoint elements $\Omega_e$ such as $\Omega = \bigcup_{e=1}^{N_e}\Omega_e$. In our model, the $\Omega_e$ are either quadrilaterals (in 2D) or hexahedra (in 3D). Within each element $\Omega_e$, the prognostic vector $\qb$ is approximated as a finite-dimensional projection $\qb_N$ using basis functions $\psi_j(\xb)$ and element-wise nodal coefficients $\qb_j^{(e)}(t)$ such that
\begin{equation}\label{eq:approximation_q}
  \qb_N(\xb,t) = \sum_{j=1}^{M_N} \psi_j(\xb) \qb_j^{(e)}(t),
\end{equation}
where $N$ denotes the order of the basis functions, $M_N=(N+1)^d$ is the number of nodes per element, and the superscript $(e)$ denotes element-wise quantities. The basis functions are constructed as tensor products of $N$-th order Lagrange polynomials $h_\alpha$ associated with the Legendre-Gauss-Lobatto (LGL) points, given by
\begin{equation}
  \psi_i(\xi,\eta,\zeta) = h_\alpha(\xi)\otimes h_\beta(\eta) \otimes h_\gamma(\zeta),
\end{equation}
where $\alpha,\beta,\gamma\in \qty{1,\cdots,N+1}$, and $(\xi,\eta,\zeta)$ are the element coordinates in the reference domain $E=[-1,1]^d$ that are mapped from the physical coordinate $\xb$ via metric terms. 

For a continuous Galerkin (CG) method, let us consider a finite-dimensional Sobolev space $\mathcal{V_N}$ defined as
\begin{equation}
  \mathcal{V_N} = \qty\big{ \psi \; | \; \psi\in H^1(\Omega) \text{ and } \psi\in\mathcal{P}_N \left( \Omega_e \right)},
\end{equation}
where $\mathcal{P}_N \left( \Omega_e \right)$ is the set of all polynomials of degree less than or equal to $N$ on the element $\Omega_e$. The discrete weak form of the governing equation in \eqref{eq:governing:compact} is obtained by multiplying it by a test function and integrating over the domain, which is stated as follows: Find $\qb_N\in\mathcal{V_N}$ such that for all $\psi\in\mathcal{V_N}$,
\begin{equation}\label{eq:governing:weak}
  \sum_{e=1}^{N_e}\intOe{\psi_i\pdv{\qb_N}{t}} = \sum_{e=1}^{N_e}\intOe{\psi_i \mathcal{S}(\qb_N)}.
\end{equation}
In 3D, the integral over the domain $\Omega$ is computed using a quadrature rule based on the LGL points as follows:
\begin{equation}\label{eq:integral}
  \intO{(\,\cdot\,)} = \sum_{e=1}^{N_e}\intOe{(\,\cdot\,)} = \sum_{e=1}^{N_e}\sum_{i=1}^{N_{\xi}+1}\sum_{j=1}^{N_{\eta}+1}\sum_{k=1}^{N_{\zeta}+1}w^{(e)}_{ijk} J^{(e)}_{ijk}(\,\cdot\,),
\end{equation}
where $N_{\xi}$, $N_{\eta}$, and $N_{\zeta}$ are the orders of polynomial basis functions in each direction, and $w^{(e)}_{ijk}$ and $J^{(e)}_{ijk}$ denote the weight and Jacobian determinant associated with the LGL points, respectively.

The discrete weak form in \eqref{eq:governing:weak} yields a matrix-vector form written as
\begin{equation}
  M_{IJ}\pdv{\qb_J}{t} = R_I(\qb_N),
\end{equation}
where $I,J\in\qty{1,\cdots,N_p}$ and $N_p$ is the total number of global grid points, and the standard summation convention is assumed unless otherwise specified. The global mass matrix $M_{IJ}$ and the right-hand side (RHS) vector $R_I$ are constructed by applying the global assembly or direct stiffness summation (DSS) operator $\bigwedge_{e=1}^{N_e}$ to the element-wise mass matrices and RHS vectors as follows:
\begin{align}
  M_{IJ} &= \bigwedge_{e=1}^{N_e} \intOe{\psi_i (\xb) \psi_j (\xb)}, \\
  R_{I}(\qb_N) &= \bigwedge_{e=1}^{N_e} \intOe{\psi_i (\xb) \mathcal{S}(\qb_N)}, \label{eq:RHSvec}
\end{align}
where the DSS operator performs summation through mapping from local indices to global indices as $(i,e)\rightarrow I$ and $(j,e)\rightarrow J$. We choose the integration points to be co-located with the interpolation points (inexact integration), which yields accurate numerical integration for $N\ge4$ \cite{giraldo1998lagrange}. This choice makes the mass matrix $M_{IJ}$ diagonal due to the cardinal property of the basis functions, simplifying inversion. Consequently, we obtain
\begin{equation}\label{eq:governing:matrix}
  \pdv{\qb_I}{t} = M^{-1}_{IJ}R_J(\qb_N).
\end{equation}
Eq.\ \eqref{eq:governing:matrix} is used to calculate the solutions for both the LSP and SSP problems in the MMF simulations, where $\mathcal{S}$ in Eq.\ \eqref{eq:RHSvec} is augmented with the corresponding coupling term, $\Fb$ or $\fb$ given in Eqs.\ \eqref{eq:forcing}--\eqref{eq:feedback}.

%% file: sections/SGS.tex
% \subsection{Subgrid scale model}

% In the current MMF approach, the subgrid scales that fall below the resolution of the SSP model in Figure \ref{fig:MMF_config} are realized through parameterizations. In this work, we employ the Smagorinsky-Lilly model to represent the effects of subgrid scale fluid motions and the Kessler scheme to represent microphysical processes.

% In the Smagorinsky-Lilly model \cite[pg. 587]{pope2000}, the turbulent eddy viscosity is defined as
% \begin{equation}
%   \nu_t =(C_S\Delta)^2\sqrt{2\Sb:\Sb},
% \end{equation}
% where $\Sb=\frac{1}{2}\qty\big(\grad\ub+\grad\ub^T)$ is the strain rate tensor, $C_S=0.14$ is the Smagorinsky coefficient, and $\Delta$ is the filter width that is defined as $\Delta=\max(\Delta x,\Delta y,\Delta z)/N$ where $N$ is the order of basis function. The viscosity parameter for the $\theta_v$ equation is defined following \cite{marras2015stabilized, reddy2022comparison} as
% \begin{align}
%   \nu_{\theta} & =\frac{\text{Pr}}{\gamma-1}\nu_t,
% \end{align}
% where $\text{Pr}$ is an artificial Prandtl number that is selected as $\text{Pr}=0.7$ for dry air, and $\gamma=c_p/c_v$ where $c_p=1004.5\text{ Jkg}^{-1}\text{K}^{-1}$ and $c_v=717.5\text{ Jkg}^{-1}\text{K}^{-1}$ are specific heat of dry air at constant pressure and at constant volume, respectively.

\subsection{Microphysics model}

For modeling the Kessler microphysics, the microphysical source terms on the right-hand side of Eqs.\ \eqref{eq:governing:theta} and \eqref{eq:governing:qv}--\eqref{eq:governing:qr} are present only in the SSP model. Meanwhile, the LSP model transports the mixing ratios for the three species of water. The microphysical source terms in Eq.\ \eqref{eq:governing} have been implemented in the NUMA model \cite{giraldo2013implicit} using both the column-based approach \cite{gabervsek2012dry} and the non-column-based approach \cite{tissaoui2023non}. We adopt the column-based approach for modeling moist microphysics in this work, since we use vertical columns to couple the LSP and SSP models.

%% file: sections/Sponge.tex
\subsection{Implicit sponge layer}

Numerical weather problems commonly require a non-reflecting boundary condition at the top of the domain to mimic the infinite vertical extent of the atmosphere. We employ the implicit Rayleigh damping technique \cite{klemp2008upper} for absorbing upper gravity waves, which adds a damping term to the right-hand side of the momentum balance equations in \eqref{eq:governing:u} as follows:
\begin{equation}
  \pdv{\ub}{t} = \mathcal{S}_{\bm{u}}(\qb) - R_w(z) (\kb\cdot\ub)\kb,
\end{equation}
where $\mathcal{S}_{\bm{u}}(\qb)$ represents the operators in the momentum balance equations, $\kb$ is the upward-pointing
unit vector, and the vertically varying damping profile $R_w(z)$ is defined as
\begin{equation}
  R_w(z) = R_{\max}\sin^2\qty[\frac{\pi}{2}\qty(\dfrac{z-z_b}{z_t-z_b})],
\end{equation}
where $z$ is the vertical coordinate, $z_b$ is the height at the bottom of the sponge layer, $z_t$ is the model top, and the coefficient $R_{\max}$ controls the magnitude of the damping to minimize the reflection while maximizing absorption. Like the buoyancy term, the damping term is treated implicitly in our IMEX scheme. Since no damping is applied to the continuity equation \eqref{eq:governing:rho}, this scheme preserves local mass conservation. We use a thickness ($\Delta z_s = z_t - z_b$) of 6 km for the standard simulation and also LSP/SSP models in MMF simulations of the squall line and supercell in Sec.\ \ref{sec:results}.

%% file: sections/Nonconforming_Grids.tex
\subsection{Nonconforming vertical grids}\label{sec:nonconforming}

Our MMF method allows the coupling of nonconforming vertical discretizations of LSP and SSP models. This capability is incorporated into the MMF models to achieve grid refinement effects not only in the horizontal direction but also in the vertical. This flexibility is particularly useful when the LSP and SSP grids have different number of elements or different orders of basis functions in the vertical. The coupling of the LSP and SSP variables, evaluated in different spaces, is achieved by projecting them onto their counterpart spaces using the element-wise least-squares projection. This projection method is adopted from the adaptive mesh refinement strategy for the element-based Galerkin method \cite{giraldo_introduction_2020, kopriva2009implementing, kopera:2014, kopera:2015, kopriva:1996b}. 

Let us consider a 1D domain, $[0,z_{top}]$, corresponding to the columns of LSP and SSP grids. Each column of the LSP and SSP model is discretized using 1D spectral elements. We assume that an LSP element is divided by multiple SSP elements, where the number of SSP elements per LSP element is denoted by $N_{S/L}\ge 1$. Let $\zhat$ and $\ztilde$ denote the natural coordinates on the LSP and SSP elements along the vertical axis in the reference element, respectively, and $\psihat$ and $\psitilde$ denote the basis functions of the LSP and SSP elements, respectively. The following relationship between the natural coordinates on the LSP and SSP elements holds:
\begin{equation}
	\zhat = s \ztilde^{(k)} + o^{(k)} \quad(k=1,\cdots,N_{S/L}),
\end{equation}
where the coefficients are determined by $s=1/N_{S/L}$ and $o^{(k)}=(2k-1)s-1$. 

The projections between the LSP and SSP variables are defined as
\begin{subequations}
  \begin{align}
    Q_i &= P^{S\rightarrow L}_{ij} q_j, \\
    q_i &= P^{L\rightarrow S}_{ij} Q_j.
  \end{align}  
\end{subequations}
The element-wise projection matrices are defined as
\begin{subequations}
  \begin{align}
    P^{S\rightarrow L}_{ij} &= \Mhat^{-1}_{ik}\qty( \int_{-1}^1 \psihat_k \psitilde_j d\zhat ), \label{eq:projection:S2L}\\
    P^{L\rightarrow S}_{ij} &= \Mtilde^{-1}_{ik} \qty(\int_{-1}^1 \psitilde_k\psihat_j d\ztilde ). \label{eq:projection:L2S}
  \end{align}
\end{subequations}
where $\Mhat_{ij}=\int_{-1}^{1}\psihat_i\psihat_j d\zhat$ and $\Mtilde_{ij}=\int_{-1}^{1}\psitilde_i\psitilde_j d\ztilde$ are 1D mass matrices. The matrices in the parenthesis in Eqs.\ \eqref{eq:projection:S2L}--\eqref{eq:projection:L2S} are calculated as follows:
\begin{subequations}
  \begin{align}
    \int_{-1}^1 \psihat_i \psitilde_j d\zhat
    &= s \sum_{k=1}^{N_{S/L}} \psihat_i \qty(s \ztilde^{(k)}_j + o^{(k)}) w_j, \\
    \int_{-1}^1 \psitilde_i\psihat_j d\ztilde
    &= \psihat_j \qty( s \ztilde^{(k)}_i + o^{(k)} ) w_i,
  \end{align}
\end{subequations}
where $w_i$ is the weight of the quadrature rule. Consequently, the projection matrices $P^{S\rightarrow L}_{ij}$ and $P^{L\rightarrow S}_{ij}$ are not square, and their dimensions are $(N+1)\times N_{S/L}(N+1)$ and $N_{S/L}(N+1)\times (N+1)$, respectively, where $N$ is the order of the basis functions.

\begin{figure}
  \centering
  \includegraphics[width=0.35\textwidth]{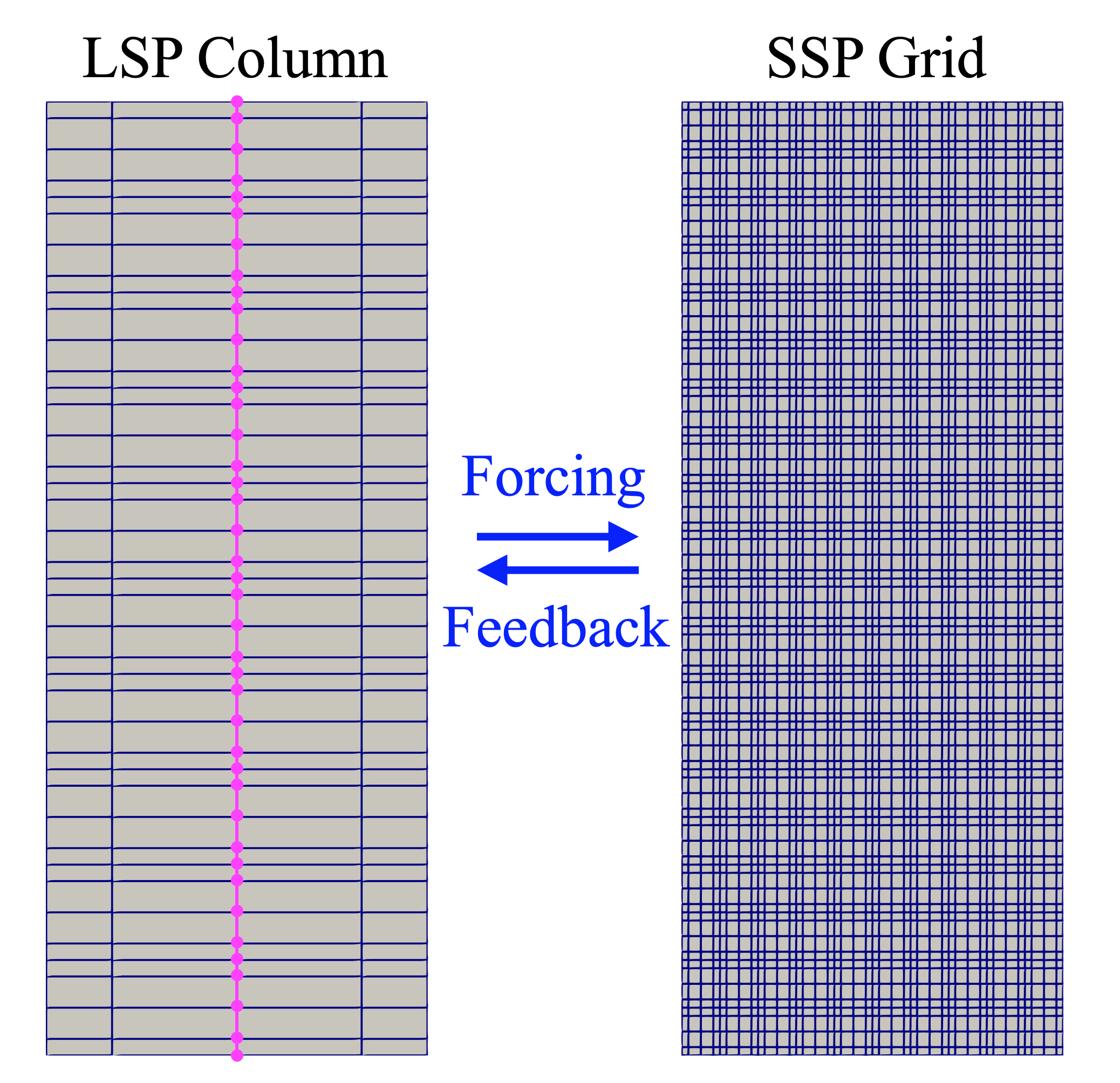}
  \caption{Coupling of nonconforming LSP and SSP grids for MMF.}\label{fig:nonconforming}
\end{figure}

Figure \ref{fig:nonconforming} illustrates an example of nonconforming LSP and SSP grids. In this figure, the LSP and SSP grids have 5 and 10 elements in the vertical direction, respectively, both with fourth order basis functions. Consequently, $N_{S/L}=2$ at each LSP column.

%% file: sections/Implementation.tex
\subsection{Parallel implementation}\label{sec:implementation}

This section describes our parallel implementation of the numerical methods for MMF simulations. The MMF simulation instantiates a single LSP simulation and multiple SSP simulations simultaneously. Each LSP and SSP simulation utilizes the dynamical core of the nonhydrostatic unified model of the atmosphere, referred to as NUMA \cite{kelly2012continuous,giraldo2013implicit}. To facilitate this task, we wrote an object-oriented
version of NUMA with multiple, independent simulator objects which we call xNUMA.  Each simulator object contains objects for input, spatial discretization, time integration, and solvers. Consequently, all the data and functions that constitute a simulation are members of the simulator objects.

\begin{figure}
  \centering
  \includegraphics[width=0.8\textwidth]{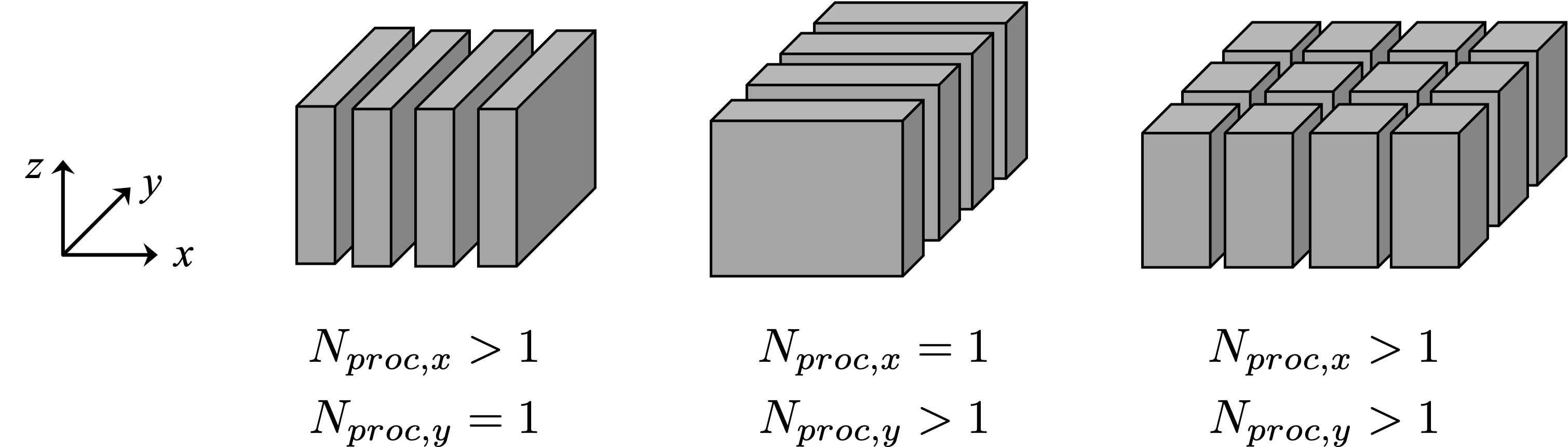}
  \caption{Layouts of domain decomposition for limited-area models. ($N_{proc,x}$ and $N_{proc,y}$ are the number of processors in the $x$ and $y$ directions.)}\label{fig:parallel:decomposition}
\end{figure}

The current MMF model is implemented via distributed memory parallelism using MPI. We subdivide the LSP domain for limited-area simulations into multiple subdomains of equal size using non-overlapping element-based partitioning. In this configuration, each element is exclusively assigned to a single partition, and each partition is assigned to an MPI process (or core). As a result, workload is well balanced among MPI processes. As we adopt the column-based model for moist microphysics, although we could also use the non-column based moist microphysics presented in Tissaoui et al.\ \cite{tissaoui2023non}, these partitions encompass all grid points and elements in the vertical direction. Figure \ref{fig:parallel:decomposition} illustrates three types of layouts for domain decomposition. An SSP grid is generated at each column, resulting in a total $N_{p,x} \times N_{p,y}$ of SSP simulations for an MMF simulation, where $N_{p,x}$ and $N_{p,y}$ are the number of grid points in the $x$ and $y$ directions. Each SSP simulator runs serially on the MPI process that contains the corresponding LSP column, and consequently, the communication between the LSP and SSP models only occur within each processor. In the MMF model, inter-process communications are required only for the global DSS operation over the LSP partitions. Therefore, the communication load for the MMF model is significantly lower than that of the standard model for the same grid size, as discussed in Section \ref{sec:communication}.

The parallelization of the LSP model is similar to what we implemented for the standard model. Details about the computational stencil for the element-based Galerkin method and the global DSS process can be found in \cite{kelly2012continuous}.

%% file: sections/Complexity.tex
\section{Complexity analysis}\label{sec:complexity}

A motivation for MMF from a computing perspective is to achieve more efficient parallelism. Numerical atmospheric models tend to be limited in performance due to being memory-bandwidth-bound rather than compute-bound in distributed memory parallelism. This has been a bottleneck in fully exploiting massively parallel systems, such as graphics processing units (GPUs). Therefore, numerical algorithms that require more computations but less memory communication will be favorable for such computer architectures. For this objective, an MMF strategy is attractive since it offloads computations from the parallel LSP to the serial SSPs, while keeping inter-process communications only at the level of the LSP model. As modern high-performance computing (HPC) increasingly shifts toward GPU-based machines, an MMF strategy becomes a promising pathway. In this section, we analyze the complexity of the MMF algorithm in terms of the number of floating-point operations, memory communications, and arithmetic intensity. 

In this section, we compare the complexity of MMF simulations to simulations without MMF, which are referred to as standard simulations.

% The MMF version of E3SM model has been ported to GPUs using a directive-based approach \cite{norman2022unprecedented} to achieve strong and weak scaling. 

%------------------------------------------------------%
%   									FLOP
%------------------------------------------------------%
\subsection{Floating-point operations}

% \rc{Use another term than FLOP. It can be confused with FLOPs.}
% \rc{Provide some background for the numbers 816 and 4635.}

We first evaluate the total number of floating-point arithmetic operations, denoted as $F$, throughout the duration of the simulations. The number of operations for the standard and MMF simulation are evaluated as
\begin{subequations}\label{eq:F1}
\begin{align}
  F^S &= N^{S}_t N^{S}_e N_p^3 (816 N_p + 4635), \\
  F^M &= N^{M}_t N^{M}_e \qty( 1 + R_t R_x R_z N_p ) N_p^3 (816 N_p + 4635) ,
\end{align}
\end{subequations}
where $N_t$ is the number of time steps, $N_e$ is the number of elements, and $N_p$ is the number of points in the $x$, $y$, $z$ directions within an element. It is assumed that an element has the same number of points in all directions. The superscripts $S$ and $M$ represent the standard and the LSP of the MMF simulations. In a single simulation, the asymptotic complexity $\mathcal{O}(N_p^4)$ arises from the gradient and divergence kernels, which constitute the dominant cost in constructing the right-rand side vector. For an MMF, $R_t$, $R_x$, and $R_z$ represent the temporal and spatial refinement ratios between the LSP and SSP ($R_t, R_x, R_z\ge 1$). Note that there is no refinement in the $y$ direction, as we only employ 2D SSP models in this work.  

In order to express $F$ in terms of grid size, let us define the number of elements $N_e$ and the number of time steps $N_t$ for standard and MMF simulations as follows:
\begin{equation}\label{eq:N_e}
  N^S_e = \dfrac{L_x L_y L_z}{N^3 \dx \dy \dz}, \qquad
  N^M_e = \dfrac{L_x L_y L_z}{N^3 R_x R_z \dx \dy \dz},
\end{equation}
where $L_x$, $L_y$, and $L_z$ are the lengths of the domain in the $x$, $y$, $z$ directions, respectively. The grid sizes $\dx$, $\dy$, and $\dz$ for the standard simulation are set equal to the grid size of the SSP in the MMF simulation. The number of time steps for the time duration $T$ are counted as
\begin{equation}\label{eq:N_t}
  N^S_t = \dfrac{T}{\dt}, \qquad
  N^M_t = \dfrac{T}{R_t\dt},
\end{equation}
where the same time step size $\dt$ is assumed for the standard and the SSP of MMF simulations.

Substituting Eqs.\ \eqref{eq:N_e} and \eqref{eq:N_t} into Eq.\ \eqref{eq:F1}, we obtain the following expressions for the number of floating-poing ($F$) operations:
\begin{subequations}\label{eq:F:final}
  \begin{align}
    F^S &= \qty(\dfrac{T}{\dt}) \qty(\dfrac{L_x L_y L_z}{N^3\dx \dy \dz}) N_p^3 \qty(816 N_p + 4635), \\
    F^M &= \qty(\dfrac{T}{\dt}) \qty( \dfrac{ L_x L_y L_z}{N^3\dx \dy \dz} ) \qty( \dfrac{1}{R_t R_x R_z} + N_p) N_p^3 \qty(816 N_p + 4635).
  \end{align}
\end{subequations}
It can be seen that the MMF simulation requires more floating-point operations than the standard simulation by a factor of $(1/(R_t R_x R_z) + N_p)$.

%------------------------------------------------------%
%   								Communication
%------------------------------------------------------%
\subsection{Communication}\label{sec:communication}

\begin{figure}
  \centering
  \includegraphics[width=0.4\textwidth]{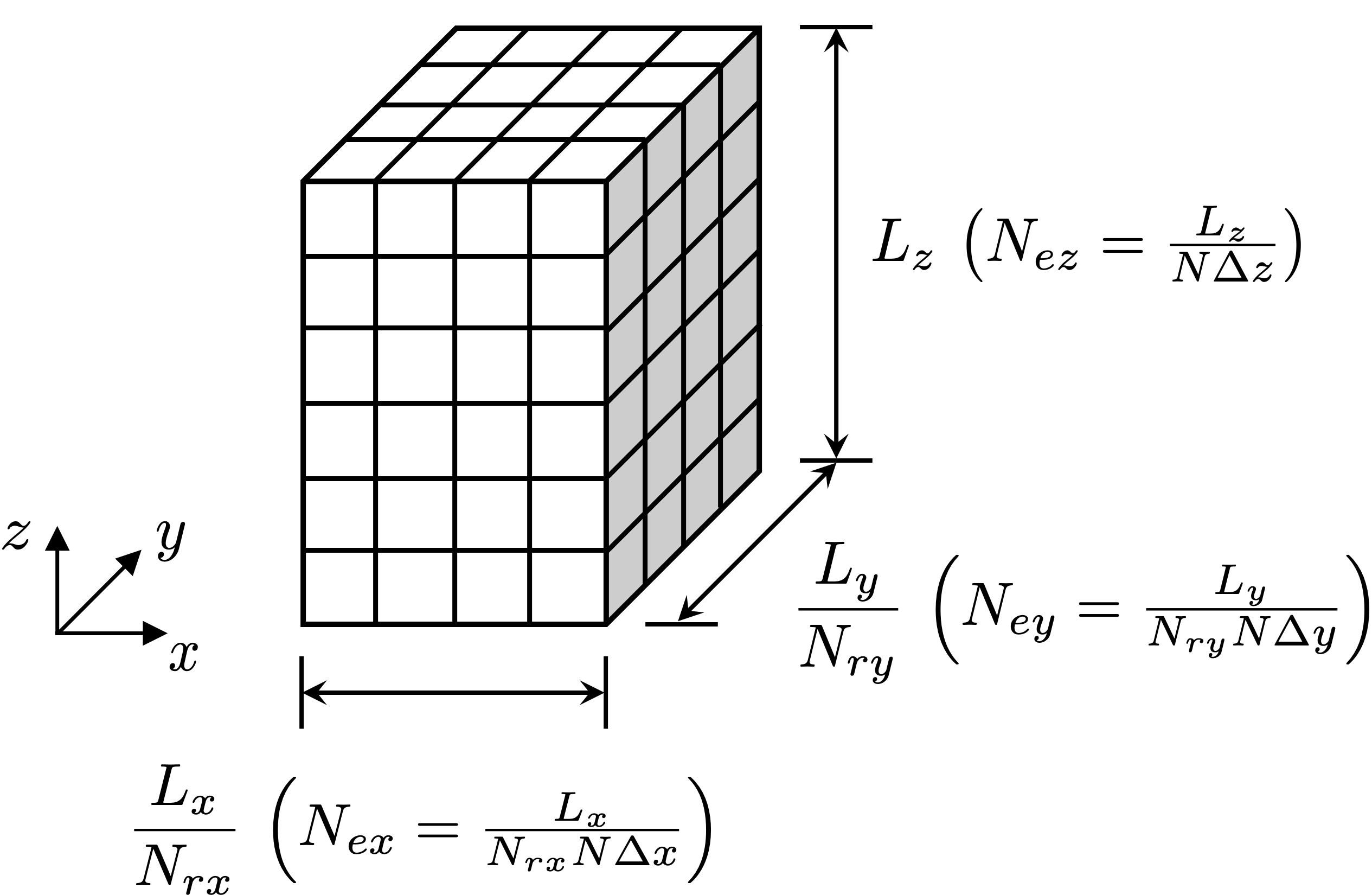}
  \caption{Dimensions of a grid partition assigned to an MPI rank and the number of elements along each direction. ($N_{ex}$, $N_{ey}$, and $N_{ez}$ are the number of elements along the $x,y,z$ directions.)}\label{fig:partition}
\end{figure}

In this analysis, we only consider communications among MPI processes in the distributed-memory system for comparisons. The transfer of shared data within cores are neglected. The data transfer occurs at the points on the lateral boundaries of the grid partitions in Figure \ref{fig:parallel:decomposition}. The total amount of floating-point numerical data to transfer or bytes ($B$) is evaluated as follows:
\begin{subequations}\label{eq:B}
  \begin{align}
    B^S &= N_t^S N_r \qty(784 N_{p,\Gamma}^S),\\
    B^M &= N_t^M N_r \qty(784 N_{p,\Gamma}^M),
  \end{align}
\end{subequations}
where $N_r$ is the number of ranks (the same for the standard and the MMF) and we assume the use of double precision so that each floating-point number occupies 8 bytes. The number of points at the lateral boundaries of the grid partition $N_{p,\Gamma}$ within a rank is
\begin{subequations}\label{eq:N_pGamma}
  \begin{align}
    N_{p,\Gamma}^S &= 2\qty( \dfrac{L_x}{N_{rx}N\dx} + \dfrac{L_y}{N_{ry}N\dy}) \qty(\dfrac{L_z}{N\dz}) N_p^2,\\
    N_{p,\Gamma}^M &= 2\qty( \dfrac{L_x}{R_x N_{rx} N\dx} + \dfrac{L_y}{N_{ry} N\dy}) \qty(\dfrac{L_z}{R_z N\dz}) N_p^2,
  \end{align}
\end{subequations}
where $N_{rx}$ and $N_{ry}$ are the number of ranks in the $x$ and $y$ directions ($N_r=N_{rx}\times N_{ry}$), respectively. Figure \ref{fig:partition} shows the number of elements along each direction at the communicating boundaries. Since inter-process communications occur only in the LSP model for MMF, there are fewer communicating points and, consequently, less data movement in the MMF. Substituting Eqs.\ \eqref{eq:N_t} and \eqref{eq:N_pGamma} in Eq.\ \eqref{eq:B}, we obtain
\begin{subequations}\label{eq:B:final}
  \begin{align}
    B^S &= 1568 N_r \qty(\dfrac{T}{\dt}) \qty( \dfrac{L_x}{N_{rx}N\dx} + \dfrac{L_y}{N_{ry}N\dy}) \qty( \dfrac{L_z}{N\dz} )N_p^2, \\
    B^M &= \dfrac{1568 N_r}{R_t R_z}\qty(\dfrac{T}{ \dt}) \qty( \dfrac{L_x}{R_x N_{rx}N\dx} + \dfrac{L_y}{N_{ry}N\dy}) \qty( \dfrac{L_z}{N\dz} ) N_p^2.
  \end{align}
\end{subequations}

%------------------------------------------------------%
%   								Arithmetic Intensity
%------------------------------------------------------%
\subsection{Arithmetic intensity}

Arithmetic intensity $I$ is the ratio of total floating-point ($F$) operations to total data movement/bytes ($B$), calculated as $I=F/B$. For simplicity, we assume the following: $N_{rx}=N_{ry}$, $L_x=L_y=L$, $\dx=\dy$, $R_x = R_t = R$ and $R_z=1$. Consequently, the arithmetic intensity for the standard and MMF simulations is calculated using Eqs.\ \eqref{eq:F:final} and \eqref{eq:B:final} as
\begin{subequations}\label{eq:I}
  \begin{align}
    I^S &= \dfrac{F^S}{B^S} =  \dfrac{1} { 2\qty( \dfrac{N_{rx}N\dx}{L} ) } N_p (0.520 N_p + 2.956), \\
    I^M &= \dfrac{F^M}{B^M} = \dfrac{  \qty( \dfrac{1}{R} + R N_p) }
    { \qty( \dfrac{1}{R} + 1) \qty( \dfrac{N_{rx}N\dx}{L} ) } N_p (0.520 N_p + 2.956).
  \end{align}
\end{subequations}
Equation \eqref{eq:I} illustrates that the arithmetic intensity for MMF is much higher than that of the standard simulation, i.e., $I^M > I^S$. This suggests that the MMF algorithm has the potential to overcome the memory-bound limitations that large-scale simulations typically encounter.

%% file: sections/Numerical_Results.tex
\section{Numerical Results}\label{sec:results}

We test the MMF method on two benchmark problems in numerical weather modeling. In this section, we refer to the simulations without the use of the MMF approach as standard simulations. We assess the performance of the MMF simulation by comparing it with these standard simulations. For the standard simulations, we consider two levels of resolution: a coarse one matching the resolution of the GCM, which is similar to the resolution of the LSP model in the MMF, and a fine one at a much higher resolution for CRM, akin to the resolution of the SSP model in the MMF. We define the grid size as $\Delta x = h/N$, where $h$ is the width of spectral element and $N$ is the order of basis function.

The simulations were conducted using the 2.0 GHz AMD EPYC 7662 processors on the NPS (Naval Postgraduate School) Supercomputer for parallel runs of the MMF simulation, each element column is assigned to a single processor. As a result, we use 18 cores for the MMF squall line and 216 cores for the MMF supercell case, matching the number of element columns in the LSP domains. 

%% file: sections/Squall_Line.tex
\subsection{2D Squall line}\label{section:squall}

We test the MMF method using the squall line test case \cite{gabervsek2012dry,tissaoui2023non}, an idealized benchmark weather problem. This test case has been used in previous studies on superparameterization and MMF methods \cite{xing2009new, majda2010new, majda2014new}. The computational domain is a 2D slab with a horizontal length of 150 km and a height of 24 km. The $x$ axis is in the streamwise direction, and the $z$ axis is in the vertical direction. We impose impermeable boundary conditions at the bottom and periodic boundary conditions at the lateral boundaries. Additionally, we apply an implicit sponge layer with a thickness of 6 km and damping coefficient $R_{\max}=0.25$ s$^{-1}$ at the top to absorb the gravity and acoustic waves \cite{klemp2008upper}. The domain is spatially discretized using fourth-order basis functions. In the horizontal direction, the grid size is approximately 200 m for the standard fine grid and approximately 4.2 km for the standard coarse grid. The grid size in the vertical direction is approximately 200 m and 400 m for the standard fine and coarse grids, respectively. The governing equations are integrated in time from $t=0$ to $t=8$ hours using the semi-implicit second-order additive Runge-Kutta (ARK2) method \cite{giraldo2013implicit} with a time-step size of $\Delta t_f=0.2$ seconds for the fine simulation and $\Delta t_c=2$ seconds for the coarse simulation. An artificial viscosity with $\nu=200$ m$^2$/s is added to the model, which is consistent with \cite{gabervsek2012dry,tissaoui2023non}. The Boyd-Vandeven filter \cite{boyd1996} is applied with strength of 1\%.

For the MMF simulation, the LSP domain retains the same dimensionality as the standard simulations, while the SSP grids are generated at every LSP grid column, each with a length of 8 km and a height of 24 km. The LSP domain is discretized with the same spatial and temporal resolutions as the standard coarse grid, while the SSP domain is discretized the same as the standard fine grid. Consequently, the LSP and SSP grids have different resolutions in the vertical direction, and there are 10 sub-steps for SSP per single LSP time step. Equal amount of artificial viscosity is added to both LSP and SSP models, and the same filter is added to the LSP model. We set up the LSP simulation identically to the standard coarse simulation to examine the properties of the MMF simulation enriched by the SSP models.

The initial state consists of the reference field for a saturated boundary layer and a thermal perturbation. The reference field is obtained from the atmosphere sounding data provided in the appendix of \cite{tissaoui2023non}, while the potential temperature perturbation at the initial time $\theta_0$ is defined as follows:
\begin{equation}
  \label{eq:squal:theta_pert}
  \theta_0 = \left\{ 
    \begin{matrix}
      \theta_c \cos^2\qty(\frac{\pi r}{2}) & r<r_c,\\
      0 & r\ge r_c,
    \end{matrix}
    \right.
\end{equation}
where 
\begin{equation}
  r = \sqrt{ \qty(\frac{x-x_c}{r_x})^2 + \qty(\frac{z-z_c}{r_z})^2 },
\end{equation}
and the parameters are set as $\theta_c=3$ K, $r_c=1$, $x_c=75$ km, $z_c=2$ km, $r_x=10$ km, and $r_z=1.5$ km, respectively. The SSP models are initialized with the profiles of the total variables along the corresponding LSP column. Consequently, the initial states of the SSP models are horizontally uniform. Additionally, we perturb these horizontally uniform potential temperature fields within the SSP models using a uniform random variable with amplitude of 0.3 K. This random perturbation is necessary to generate small-scale convection, as discussed in \cite{grabowski2001coupling}, and is defined as follows: 
\begin{equation}\label{eq:rand_pert}
  \theta_{\text{rand}}=(0.3 \text{ K})\frac{\theta_0}{\theta_c}U,
\end{equation}
where $U$ is a uniform random variable on $[-1,1]$, resulting in $-0.3 \text{ K}\le\theta_{\text{rand}}\le 0.3 \text{ K}$. We chose the amplitude used in \cite{grabowski2001coupling}.

% Instantaneous fields

Figure \ref{fig:squall:snapshot} displays snapshots of the instantaneous virtual potential temperature perturbations along with the cloud and rain mixing ratio contours at $t=1500$, 3000, 6000, and 9000 seconds, illustrating the evolution of the storm. The thermal bubble ascends, forming several cloud towers that merge into a large cloud cell. These clouds are convected eastward by the prevailing background wind. The evolution of the storm, particularly the cloud water, is highly dependent on the horizontal resolution. It is observed that the MMF produces cloud patterns closer to those of the standard fine (SF) simulation than the standard coarse (SC) simulation. This is noticeable in the shape of the cloud cell in an updraft region at $t=3000$ seconds and the cloud anvil at $t=9000$ seconds. 
\newcommand\figurewidth{0.30}
\begin{figure}
  \centering
  \begin{subfigure}[b]{\figurewidth\textwidth}
    \includegraphics[width=\textwidth]{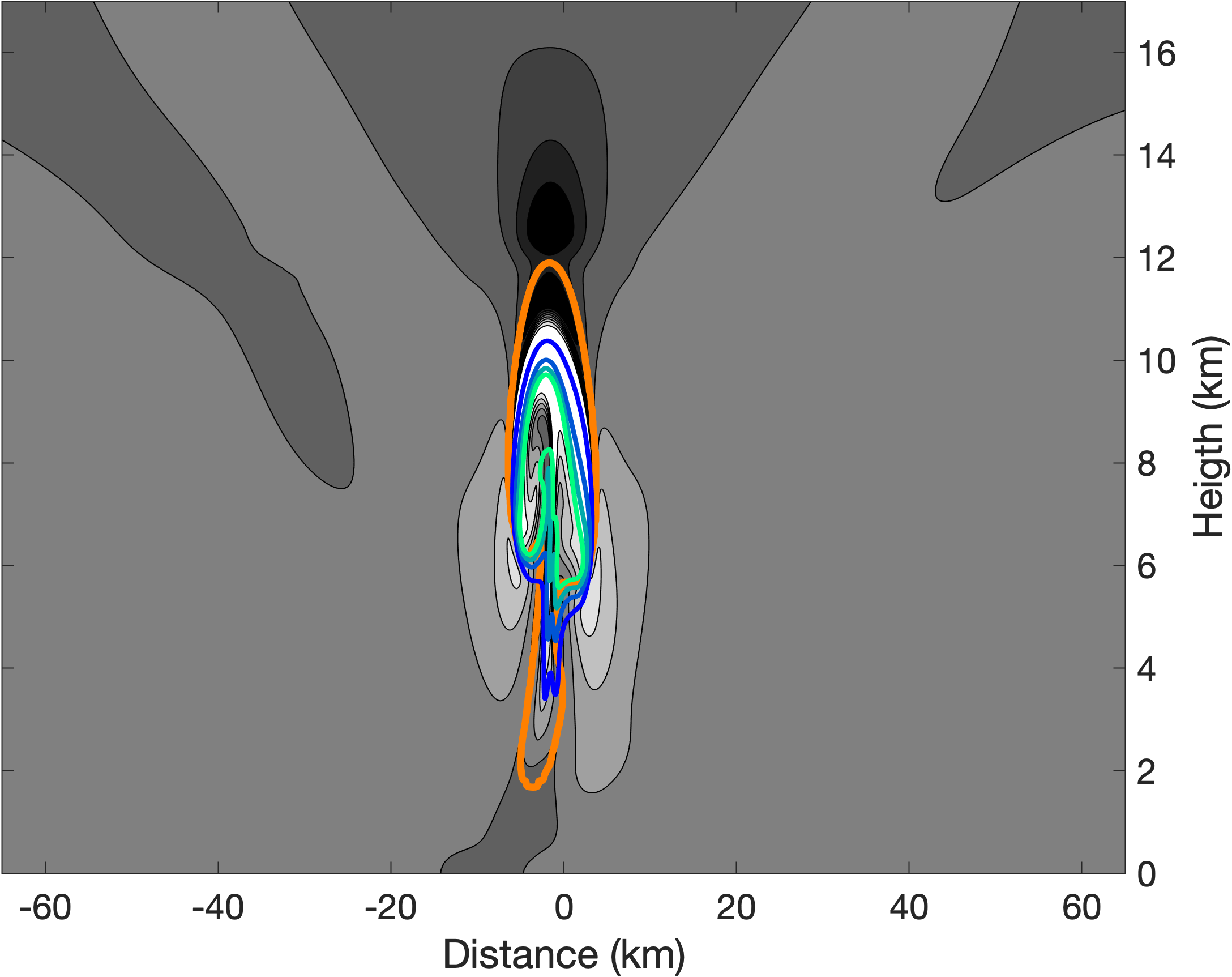}
    \caption{SF ($t=1500$)}
  \end{subfigure}
  \begin{subfigure}[b]{\figurewidth\textwidth}
    \includegraphics[width=\textwidth]{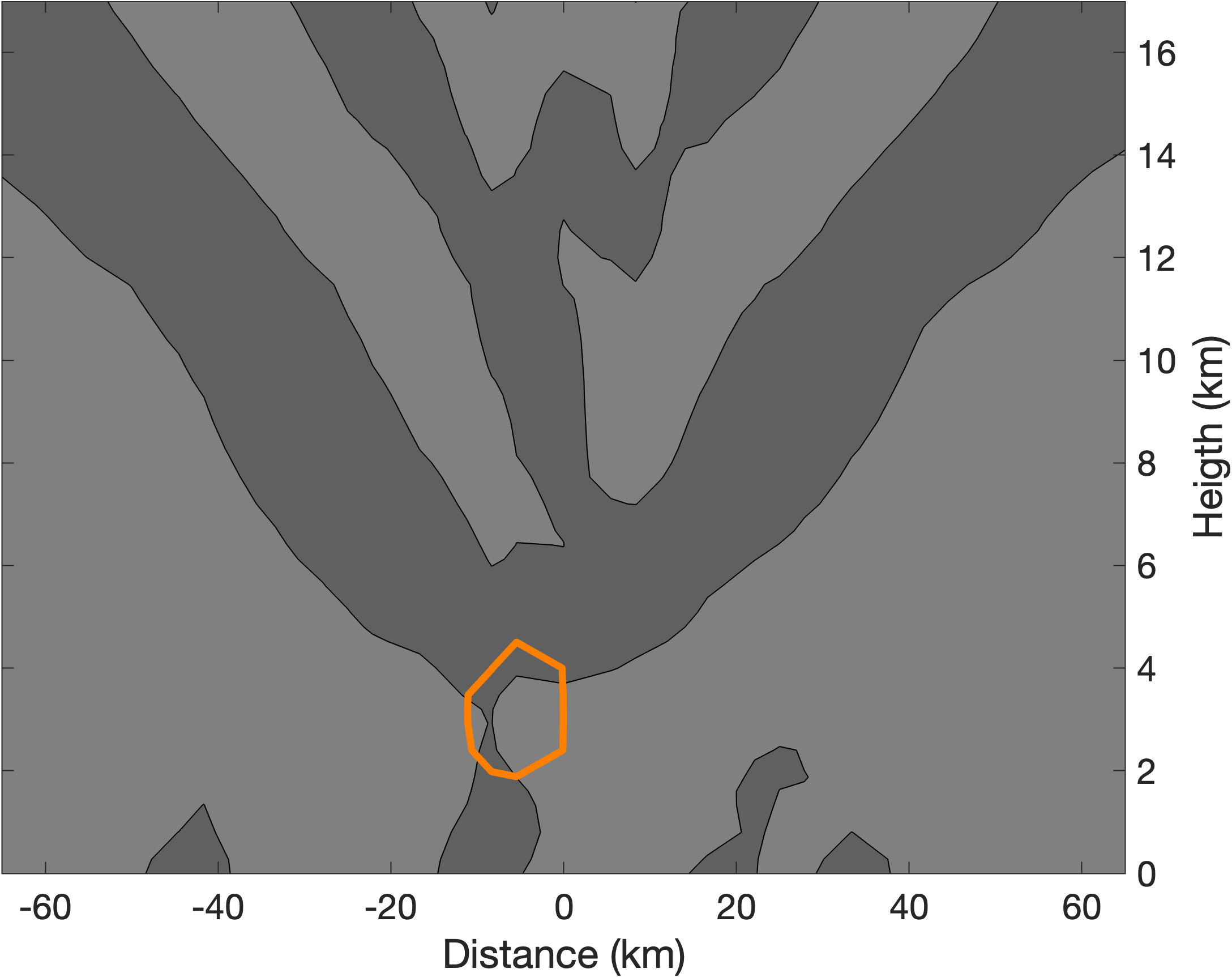}
    \caption{SC ($t=1500$)}	
  \end{subfigure}
  \begin{subfigure}[b]{\figurewidth\textwidth}
    \includegraphics[width=\textwidth]{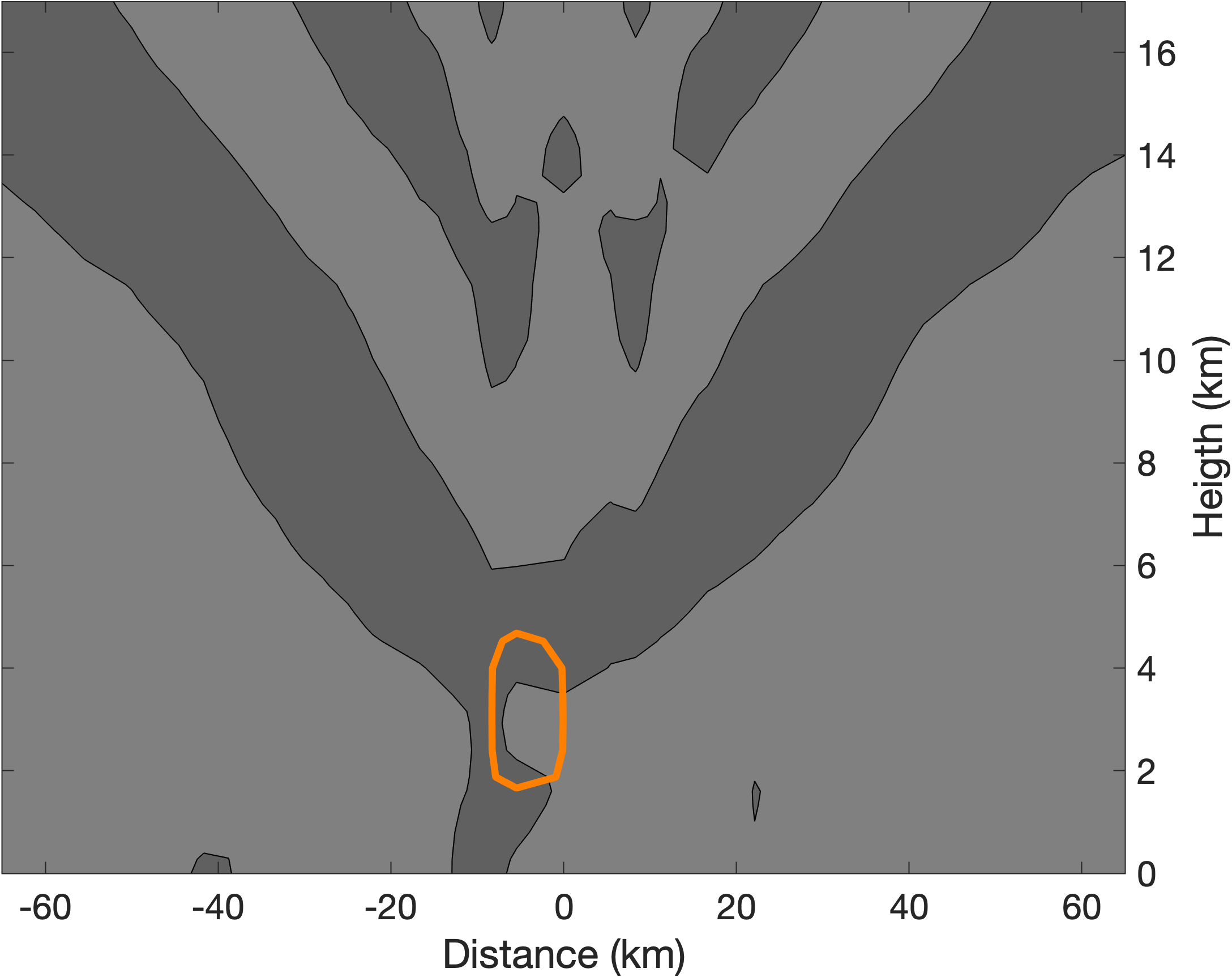}
    \caption{MMF ($t=1500$)}	
  \end{subfigure}  
  \begin{subfigure}[b]{\figurewidth\textwidth}
    \includegraphics[width=\textwidth]{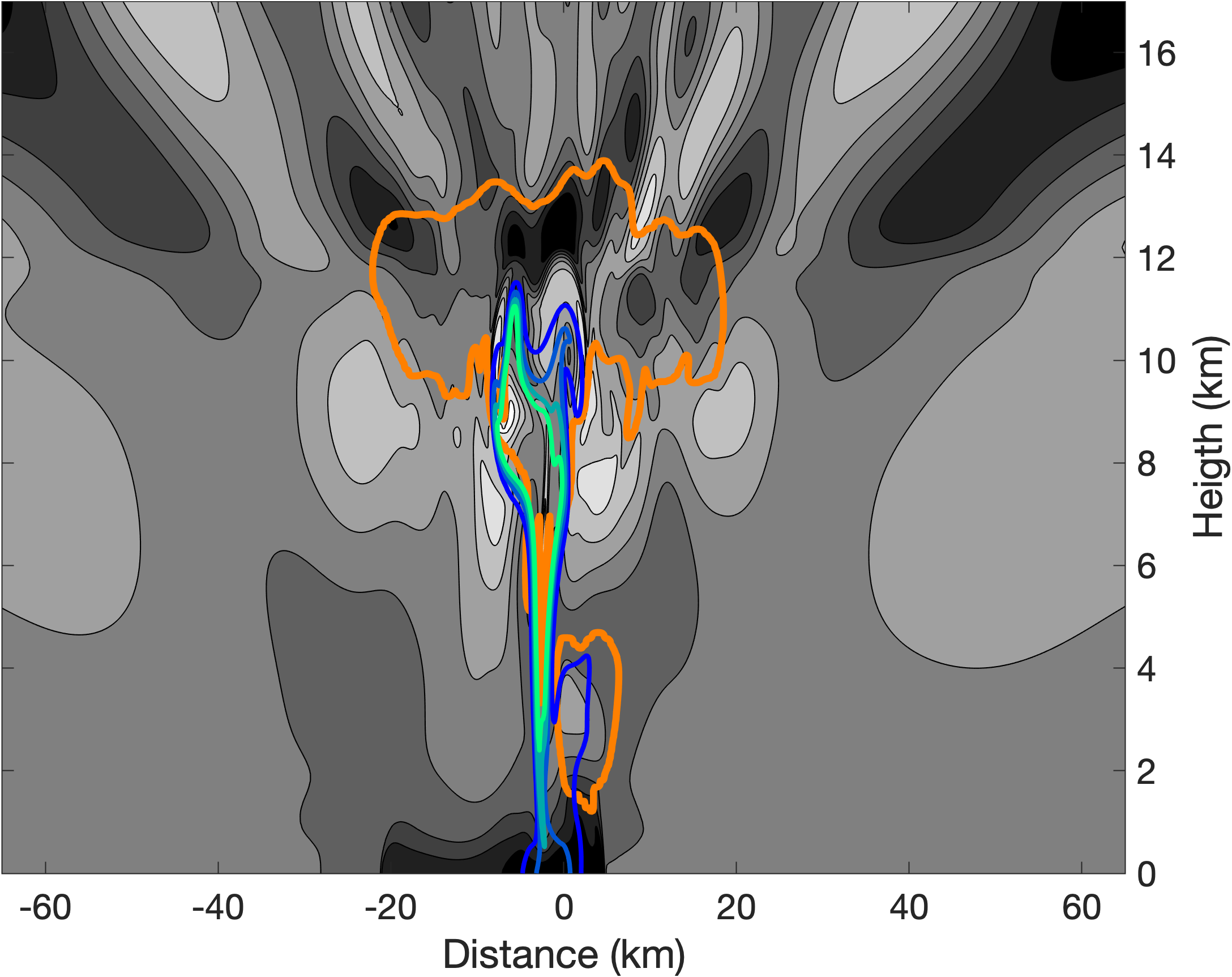}
    \caption{SF ($t=3000$)}	
  \end{subfigure}
  \begin{subfigure}[b]{\figurewidth\textwidth}
    \includegraphics[width=\textwidth]{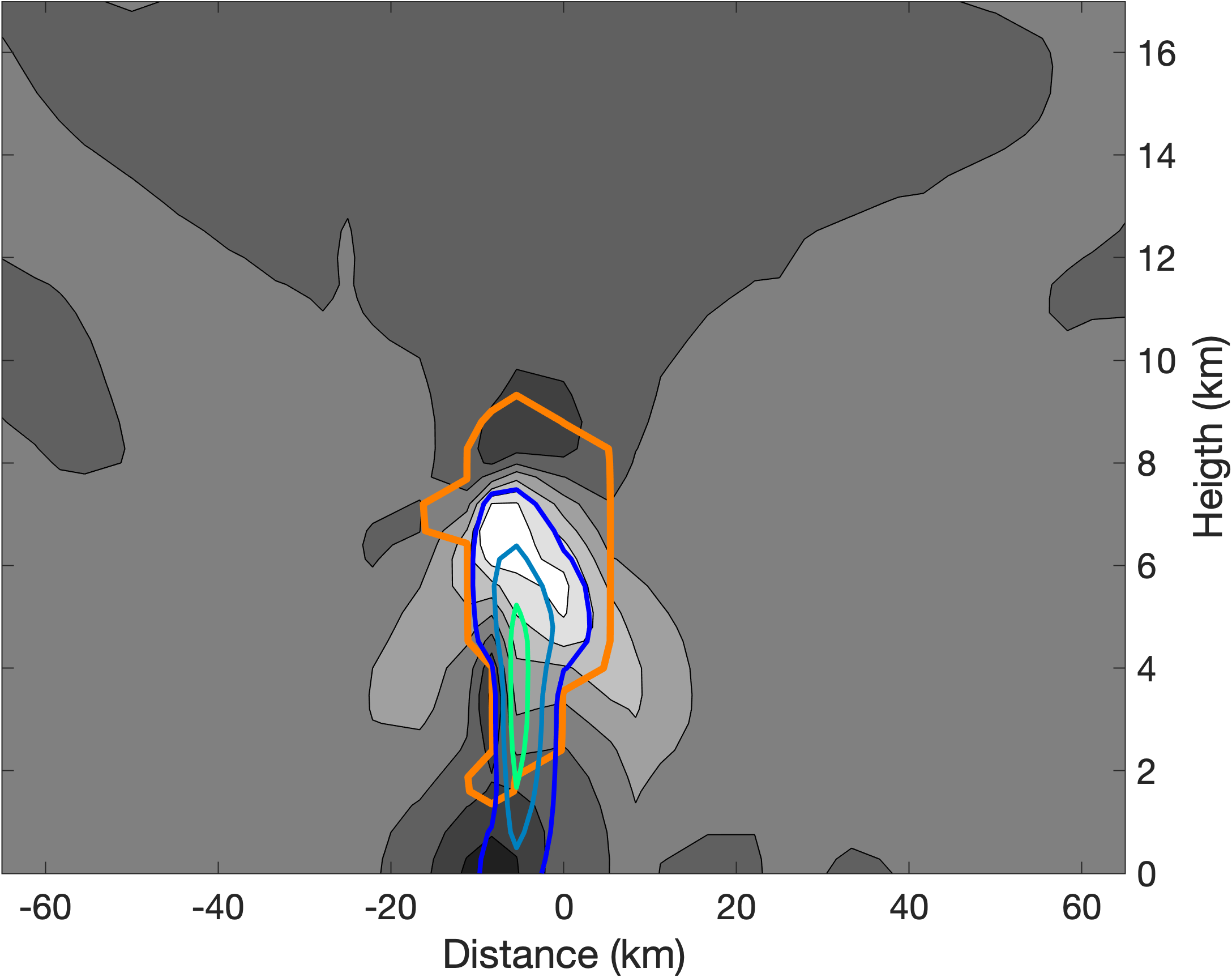}
    \caption{SC ($t=3000$)}	
  \end{subfigure}
  \begin{subfigure}[b]{\figurewidth\textwidth}
    \includegraphics[width=\textwidth]{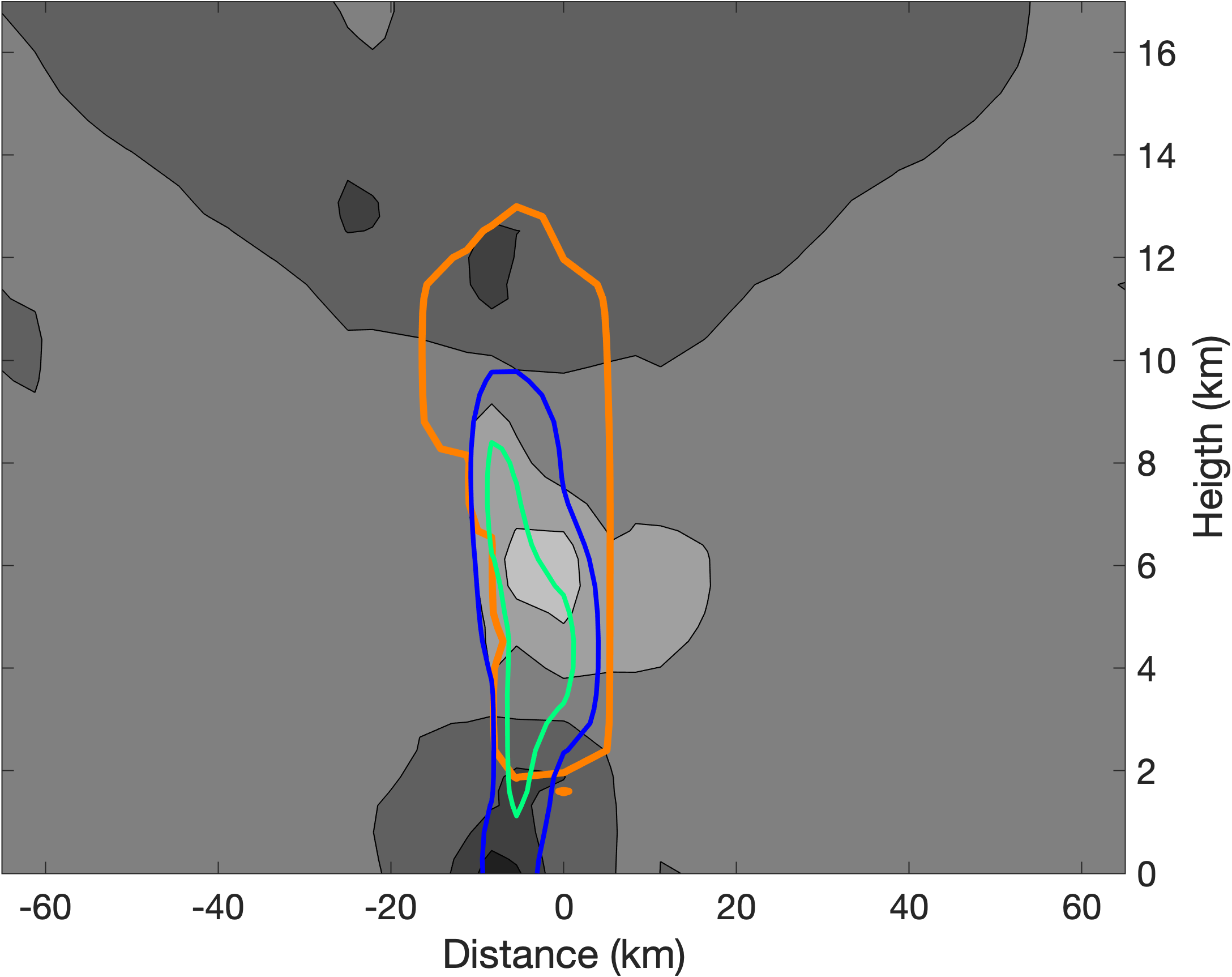}
    \caption{MMF ($t=3000$)}	
  \end{subfigure}    
  \begin{subfigure}[b]{\figurewidth\textwidth}
    \includegraphics[width=\textwidth]{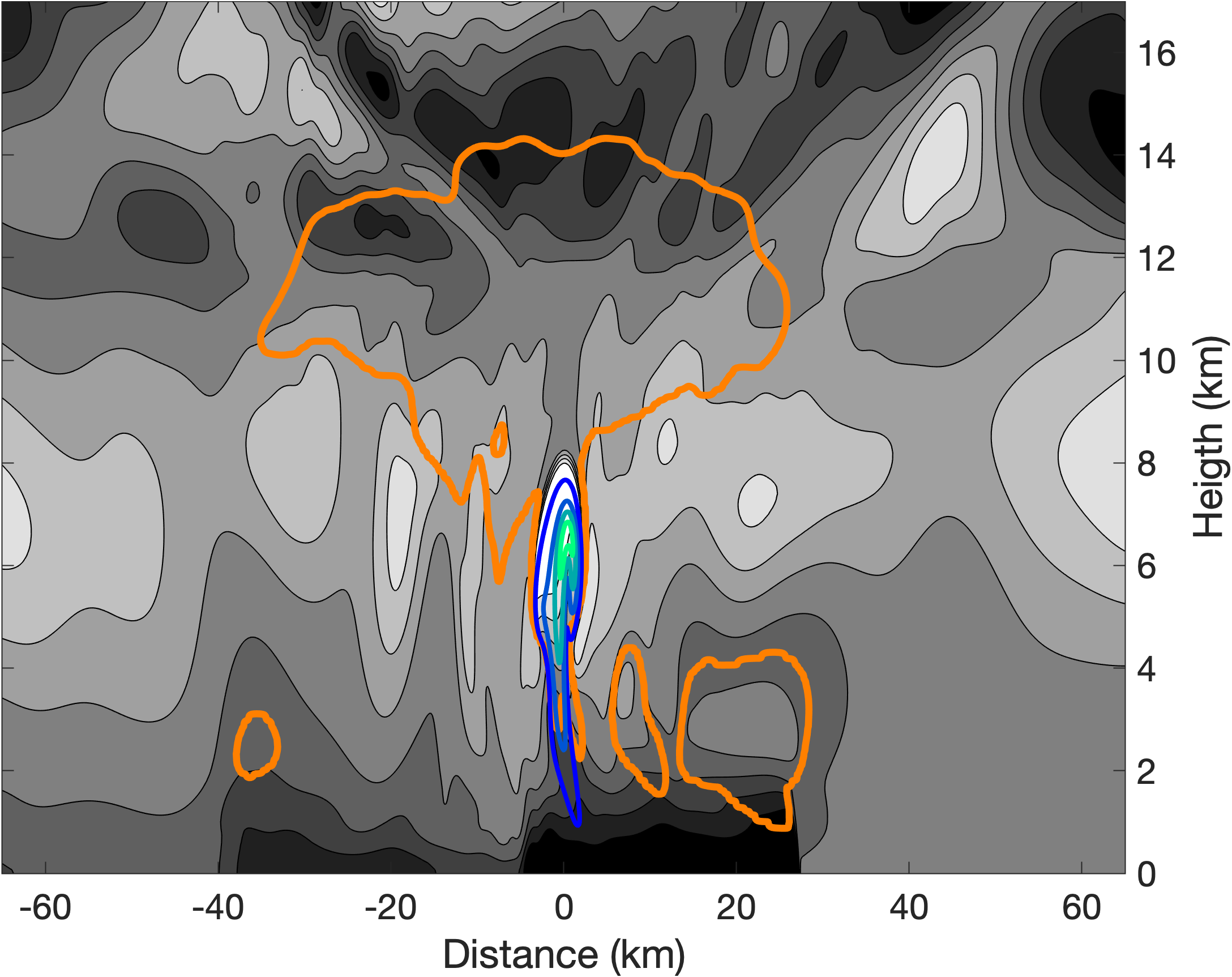}
    \caption{SF ($t=6000$)}	
  \end{subfigure}  
  \begin{subfigure}[b]{\figurewidth\textwidth}
    \includegraphics[width=\textwidth]{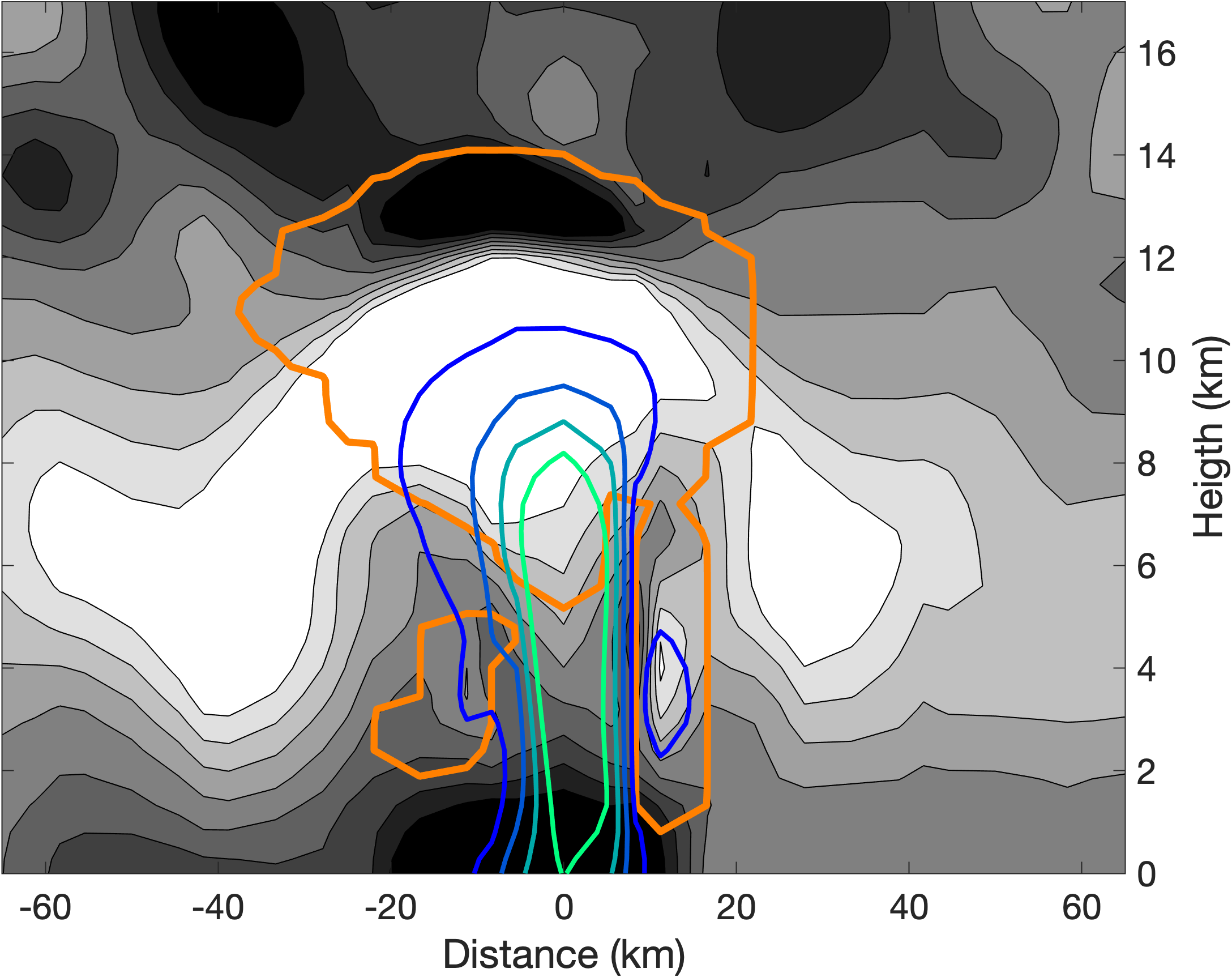}
    \caption{SC ($t=6000$)}	
  \end{subfigure}
  \begin{subfigure}[b]{\figurewidth\textwidth}
    \includegraphics[width=\textwidth]{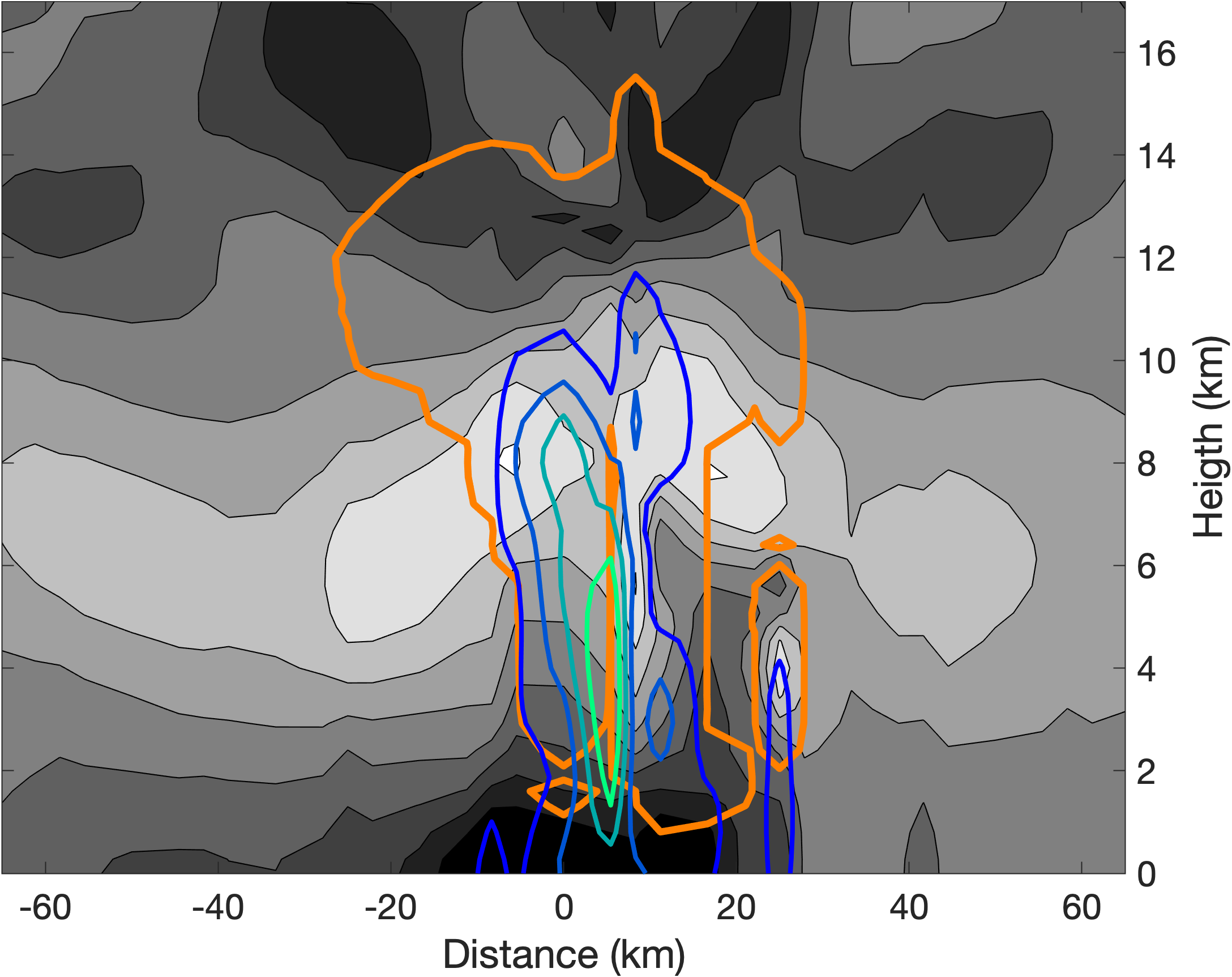}
    \caption{MMF ($t=6000$)}	
  \end{subfigure}    
  \begin{subfigure}[b]{\figurewidth\textwidth}
    \includegraphics[width=\textwidth]{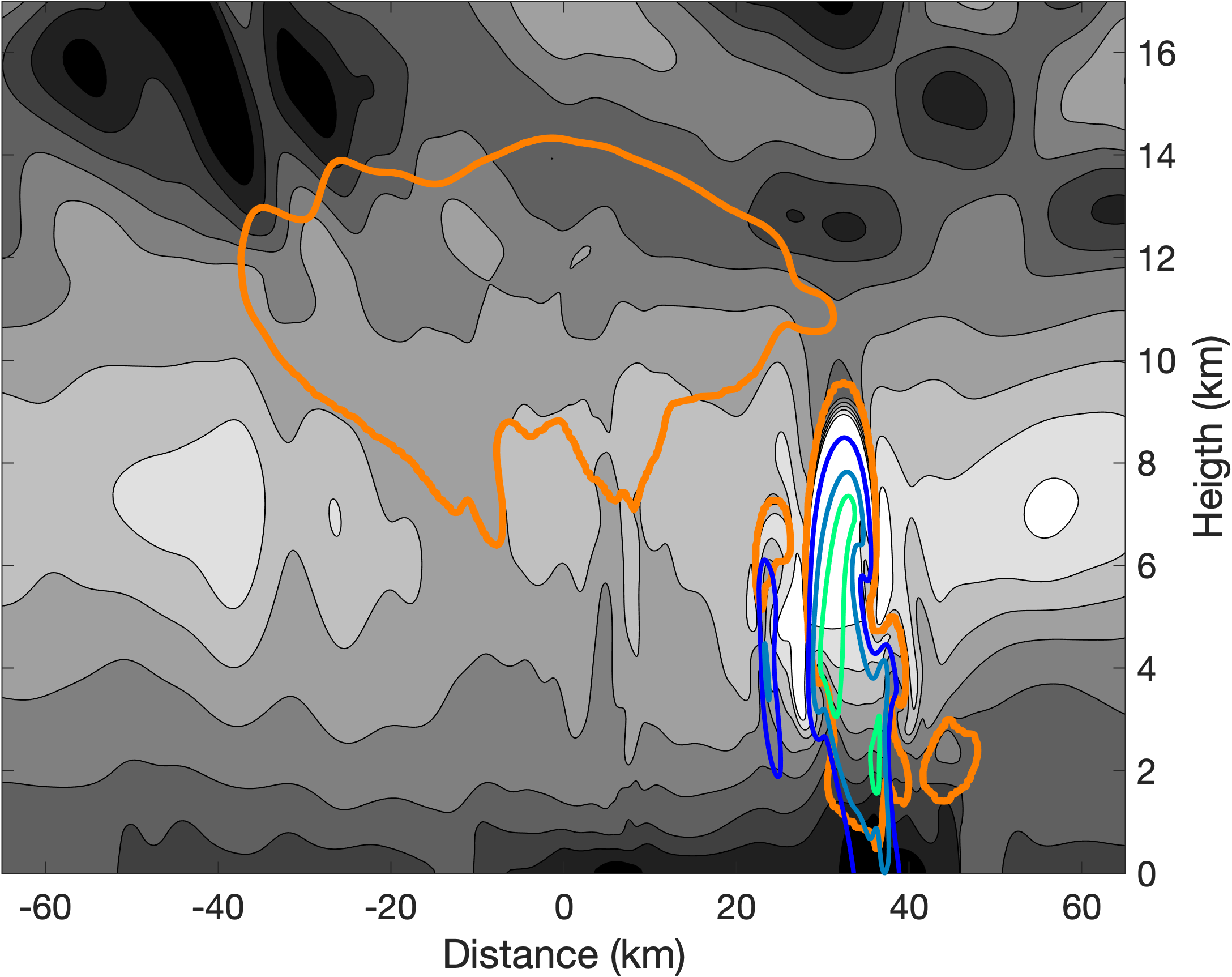}
    \caption{SF ($t=9000$)}	
  \end{subfigure}
  \begin{subfigure}[b]{\figurewidth\textwidth}
    \includegraphics[width=\textwidth]{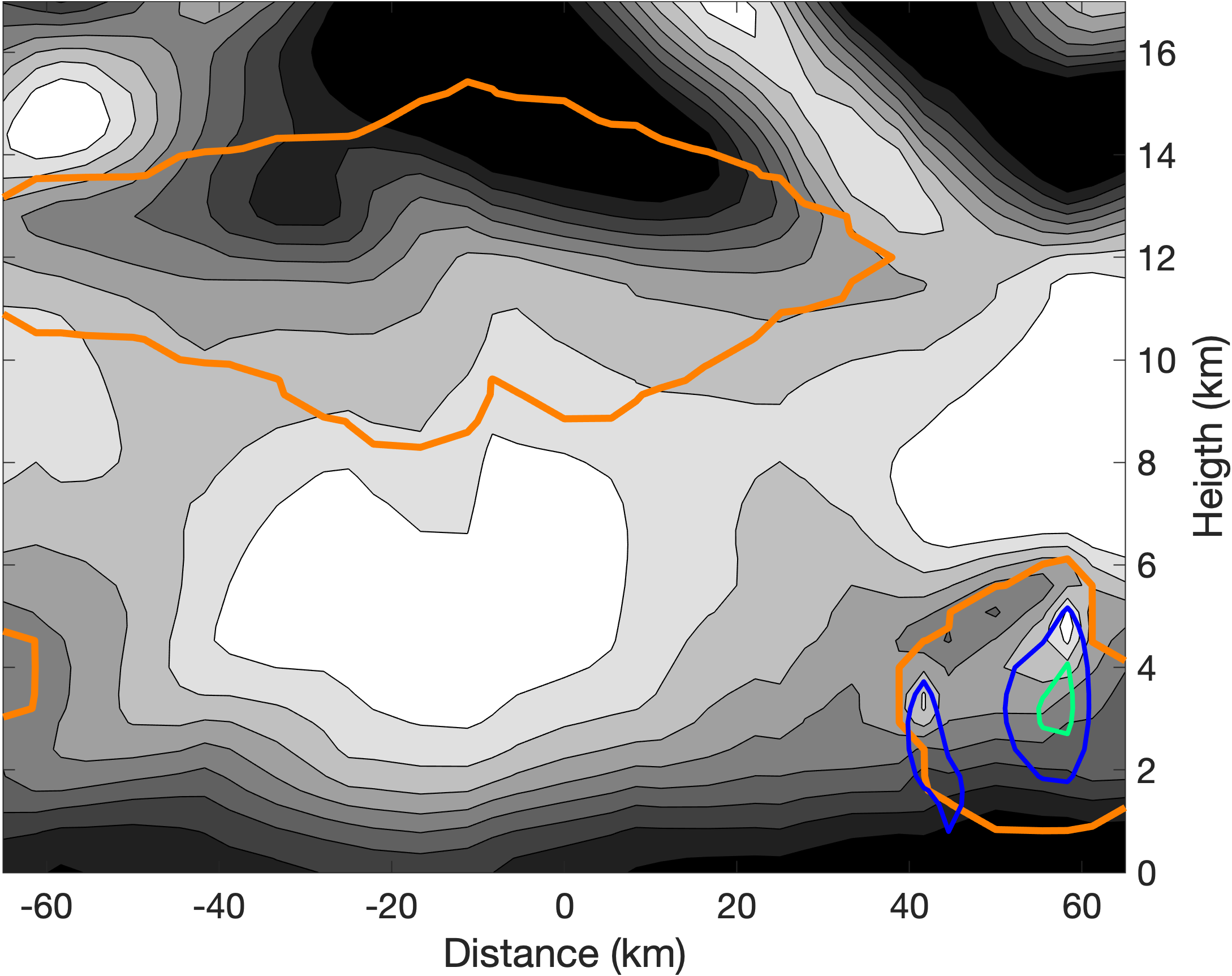}
    \caption{SC ($t=9000$)}	
  \end{subfigure}
  \begin{subfigure}[b]{\figurewidth\textwidth}
    \includegraphics[width=\textwidth]{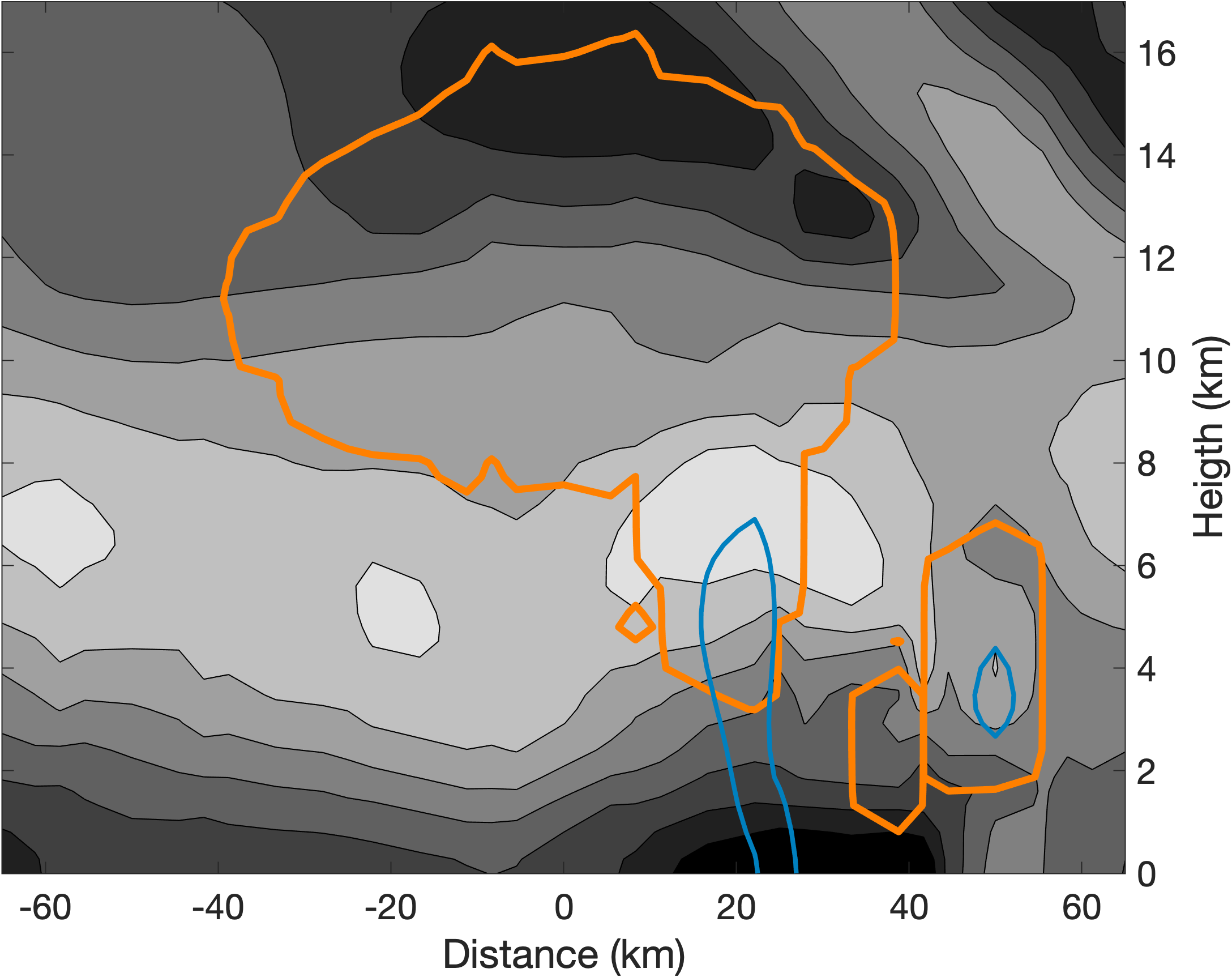}
    \caption{MMF ($t=9000$)}	
  \end{subfigure}    
  \begin{subfigure}[b]{0.32\textwidth}
    \includegraphics[width=\textwidth]{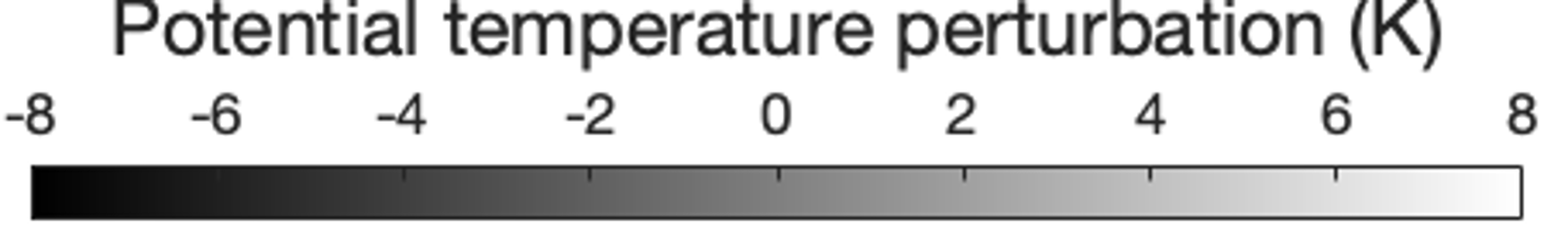}
  \end{subfigure}    
  \begin{subfigure}[b]{0.32\textwidth}
    \includegraphics[width=\textwidth]{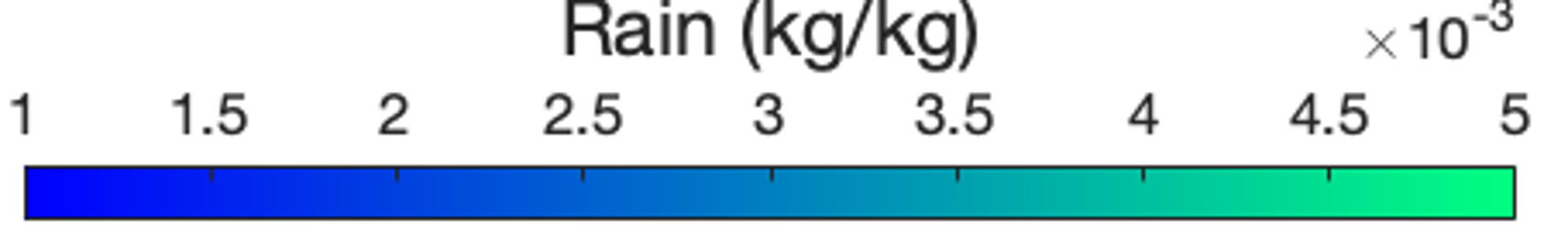}
  \end{subfigure}
  \caption{Instantaneous virtual potential temperature perturbations at $t$=1500, 3000, 6000, and 9000 seconds computed from the standard fine (SF), standard coarse (SC), and MMF simulations for the squall line. The orange-colored line represents the contour of the cloud mixing ratio at $q_c=10^{-5}$ and the lines colored in blue-green scale present the contours of the rain mixing ratio.}\label{fig:squall:snapshot}
\end{figure}

% Surface precipitation

Figure \ref{fig:squall:rain} illustrates contours of the surface precipitation on the space-time plane to show the propagating squall line from three different simulations, the standard simulations using fine and coarse grids and the MMF simulation. We observe that the pattern of surface precipitation from the MMF simulation is similar to that of the standard coarse and fine simulations.

\newcommand\figureSLRainWidth{0.32}
\begin{figure}
  \centering
  \begin{subfigure}[b]{\figureSLRainWidth\textwidth}
    \includegraphics[width=\textwidth]{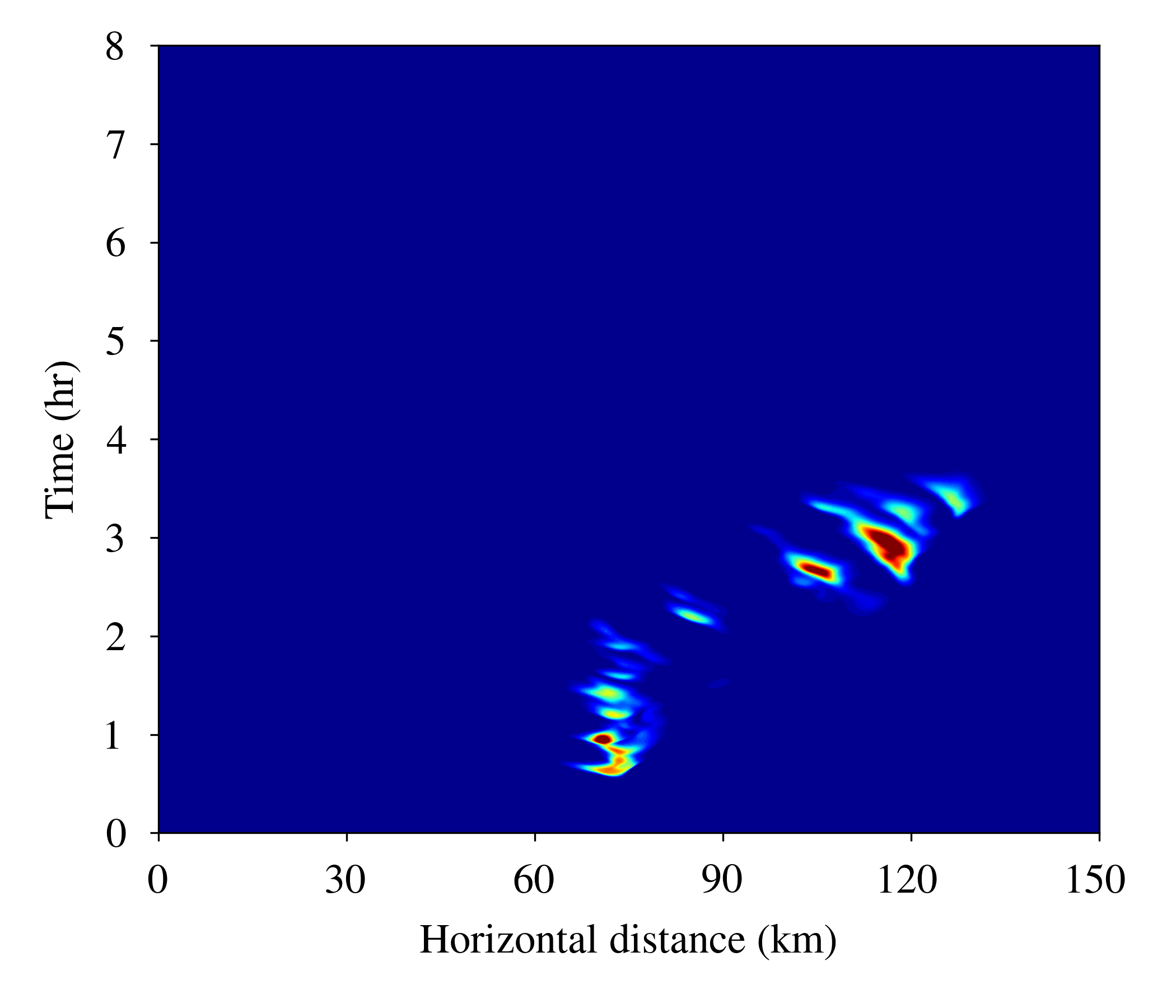}
    \caption{Standard fine}
  \end{subfigure}
  \begin{subfigure}[b]{\figureSLRainWidth\textwidth}
    \includegraphics[width=\textwidth]{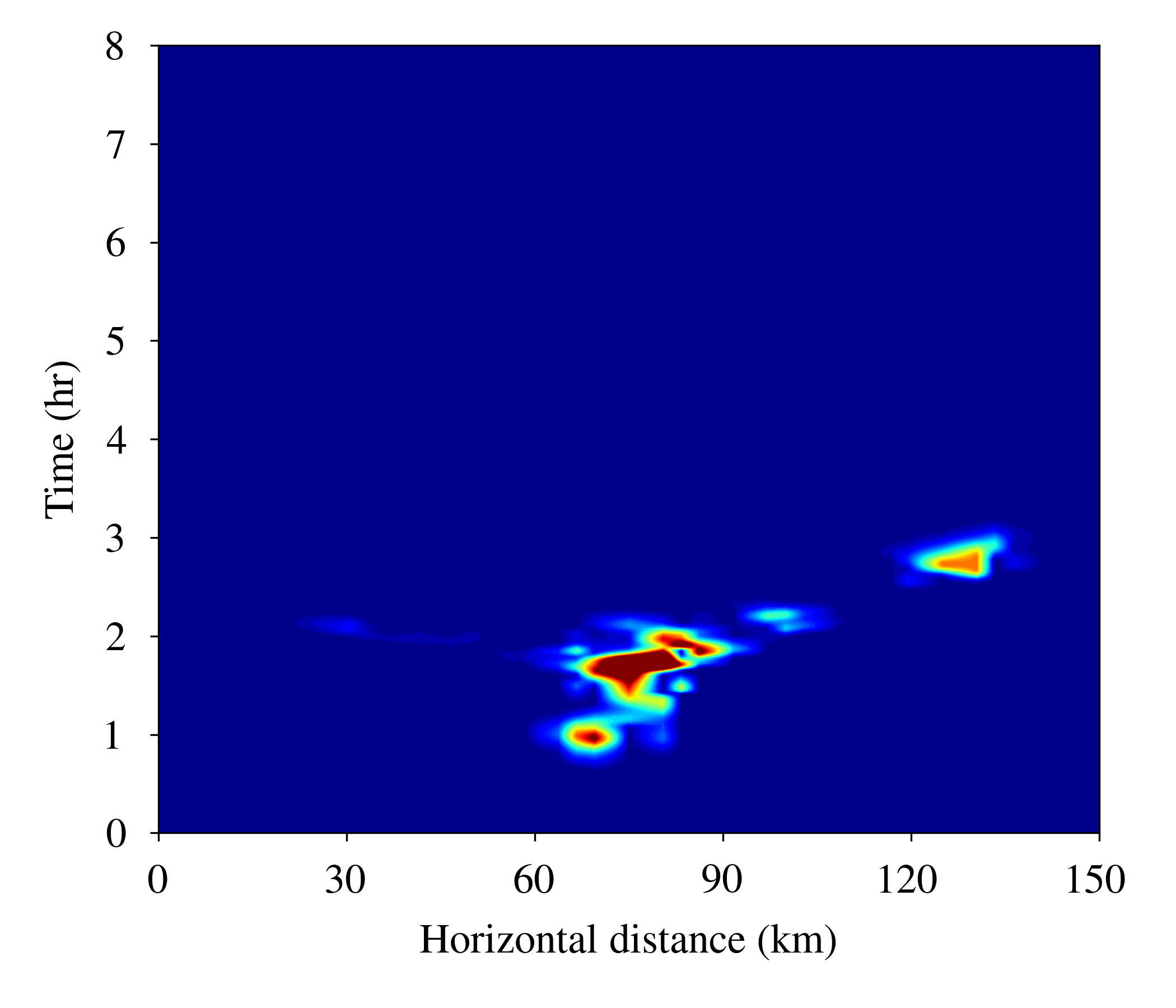}
    \caption{Standard coarse}	
  \end{subfigure}
  \begin{subfigure}[b]{\figureSLRainWidth\textwidth}
    \includegraphics[width=\textwidth]{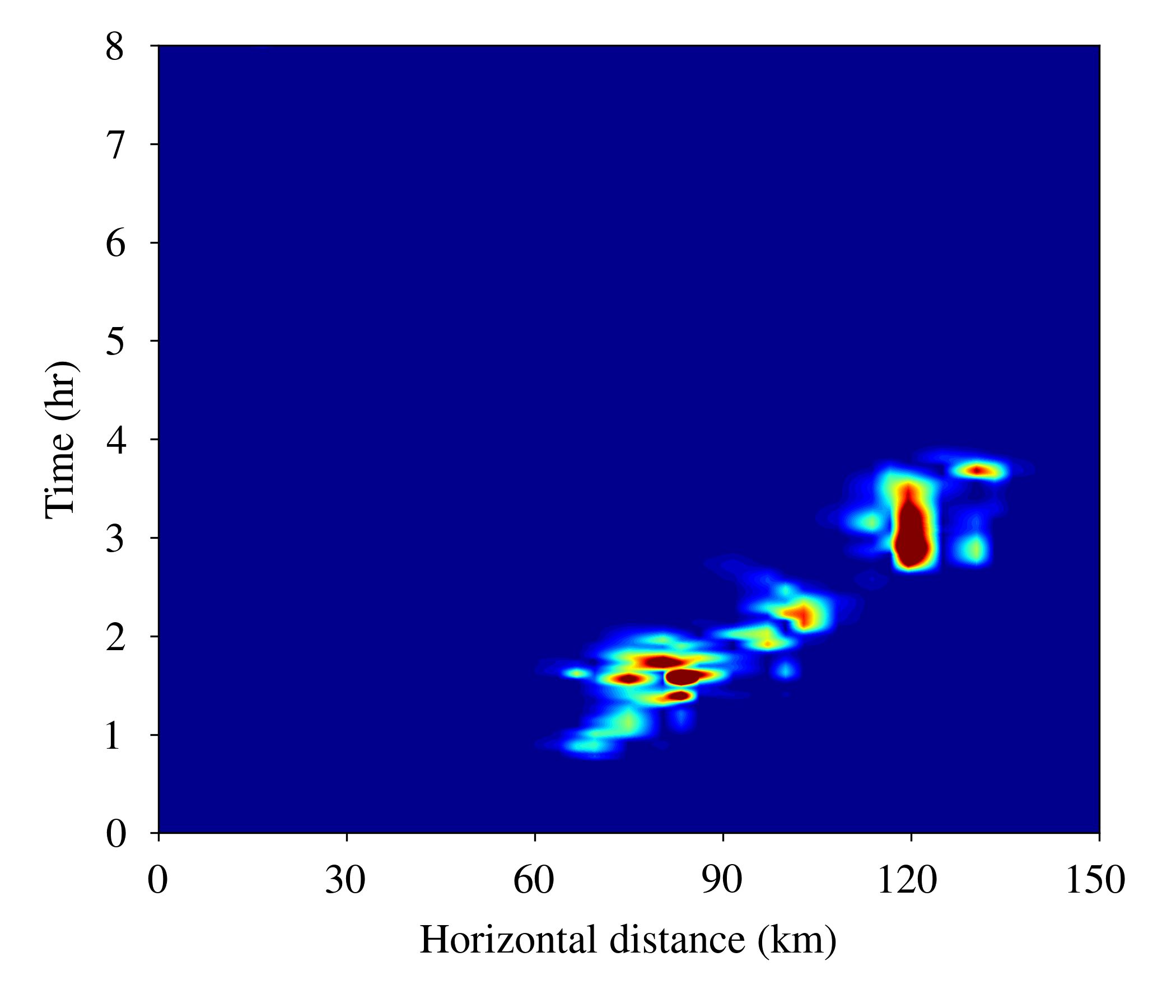}
    \caption{MMF}	
  \end{subfigure}
  \begin{subfigure}[b]{0.32\textwidth}
    \includegraphics[width=\textwidth]{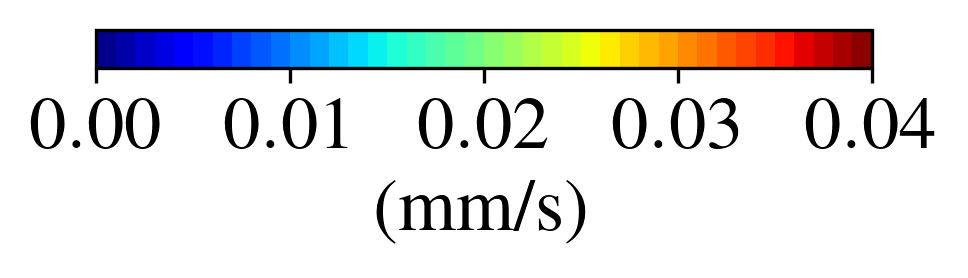}
  \end{subfigure}
  \caption{Contours of the surface precipitation in the squall line.}\label{fig:squall:rain}
\end{figure}

% Rain accumulation

\begin{figure}
  \centering
    \includegraphics[width=0.5\textwidth]{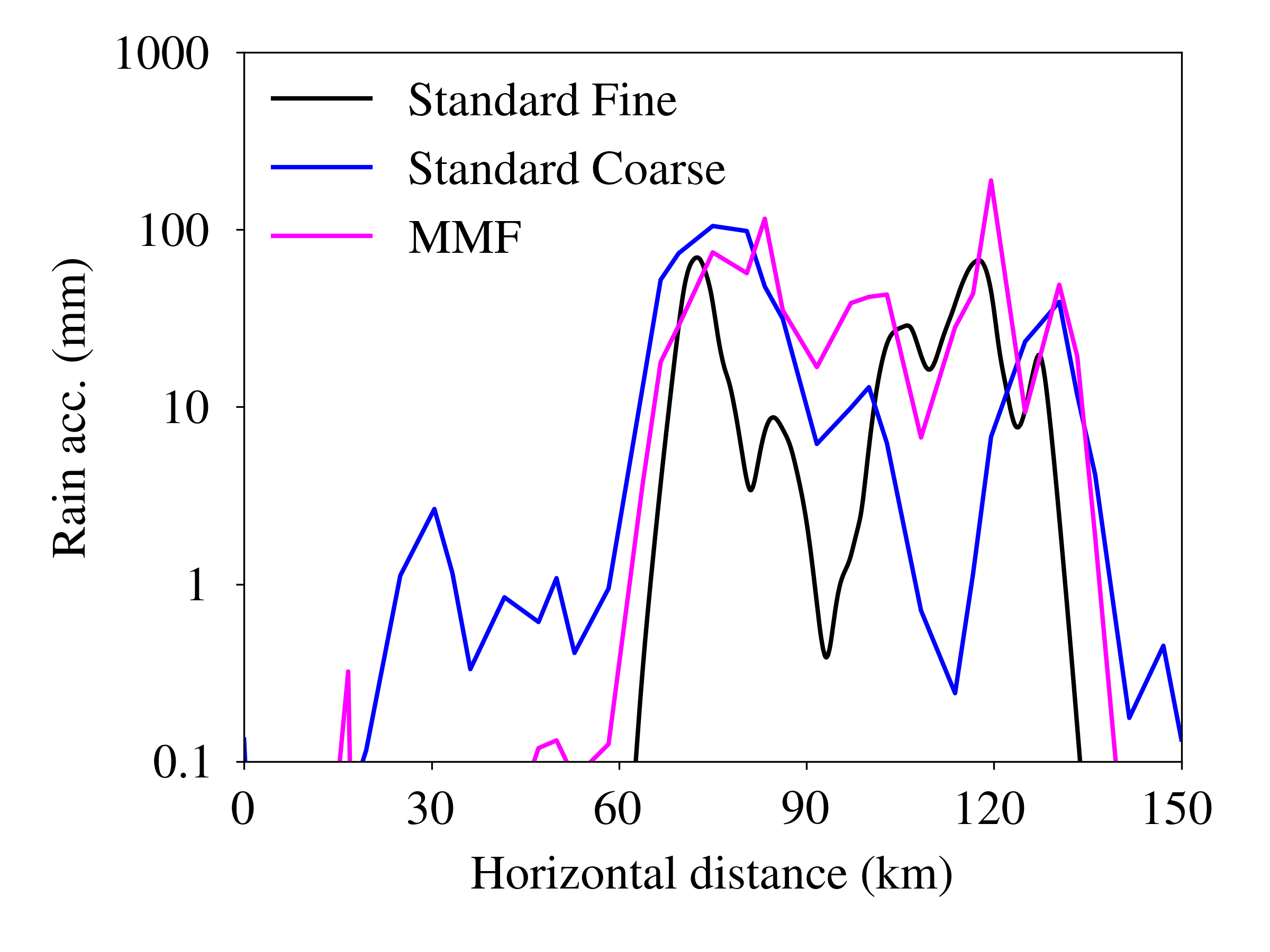}
  \caption{Distribution of the rain accumulation at $t=8$ hours for the squall line.}\label{fig:squall:rainacc}
\end{figure}

Figure \ref{fig:squall:rainacc} compares the distributions of the rain accumulation computed from the three simulations. The standard fine grid simulation yields a narrow and centered profile, while the standard coarse grid simulation produces a wide-spread profile. The profile of the MMF case is closer to the standard fine grid case in terms of shape. This result suggests the possibility that the MMF enhances accuracy in predicting the surface precipitation.

% Surface Theta (Cold pool)

\begin{figure}
  \centering
  \begin{subfigure}[b]{\figureSLRainWidth\textwidth}
    \includegraphics[width=\textwidth]{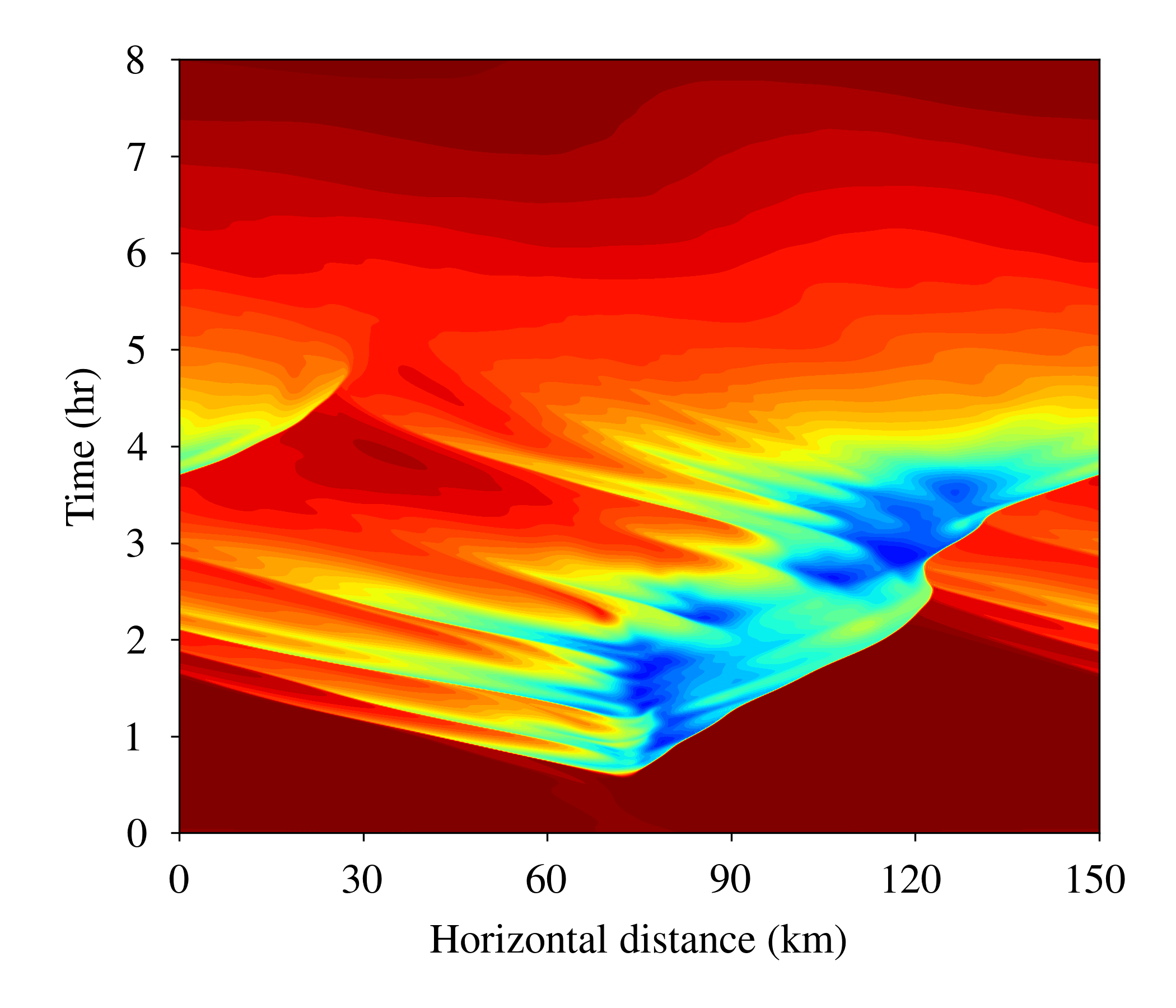}
    \caption{Standard fine}
  \end{subfigure}
  \begin{subfigure}[b]{\figureSLRainWidth\textwidth}
    \includegraphics[width=\textwidth]{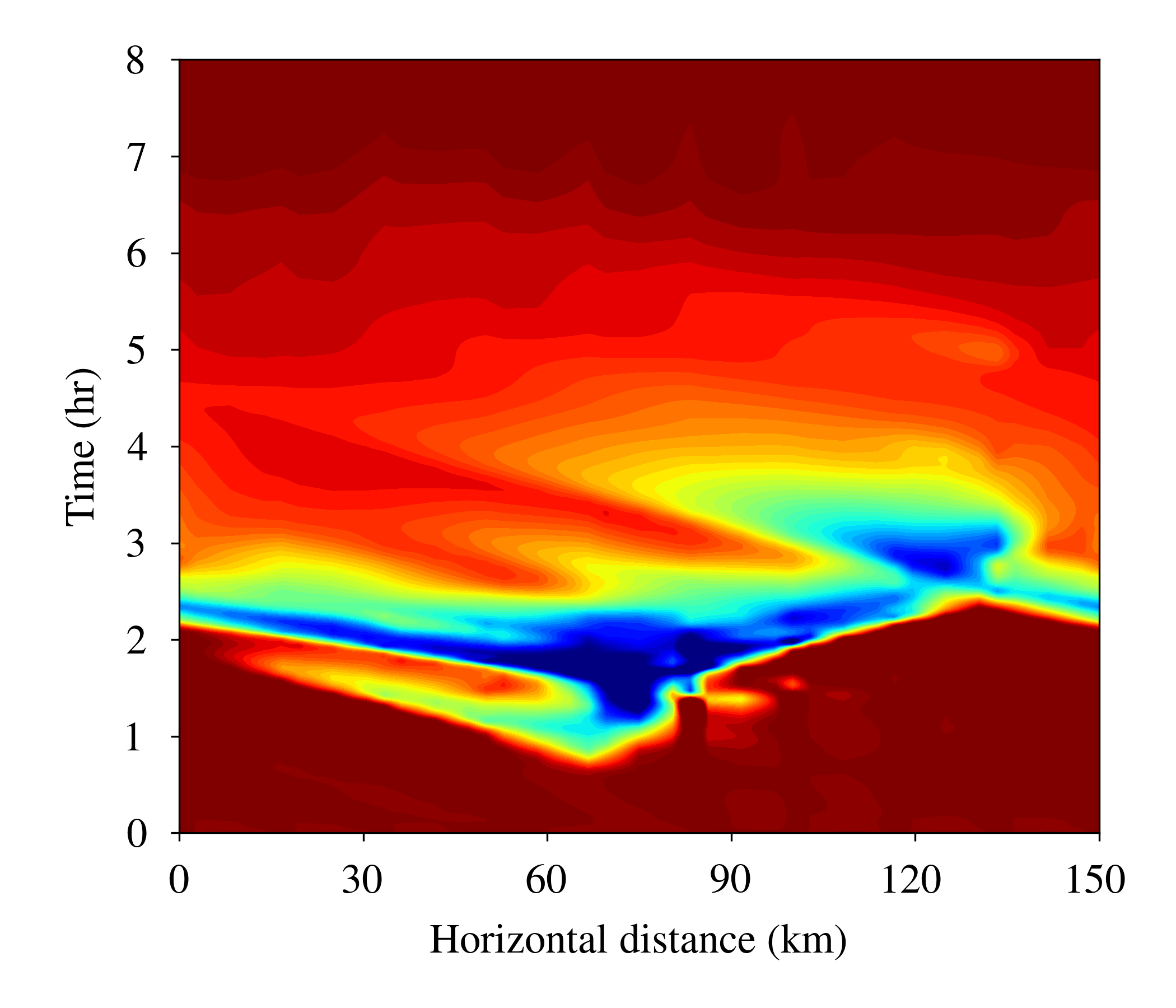}
    \caption{Standard coarse}	
  \end{subfigure}
  \begin{subfigure}[b]{\figureSLRainWidth\textwidth}
    \includegraphics[width=\textwidth]{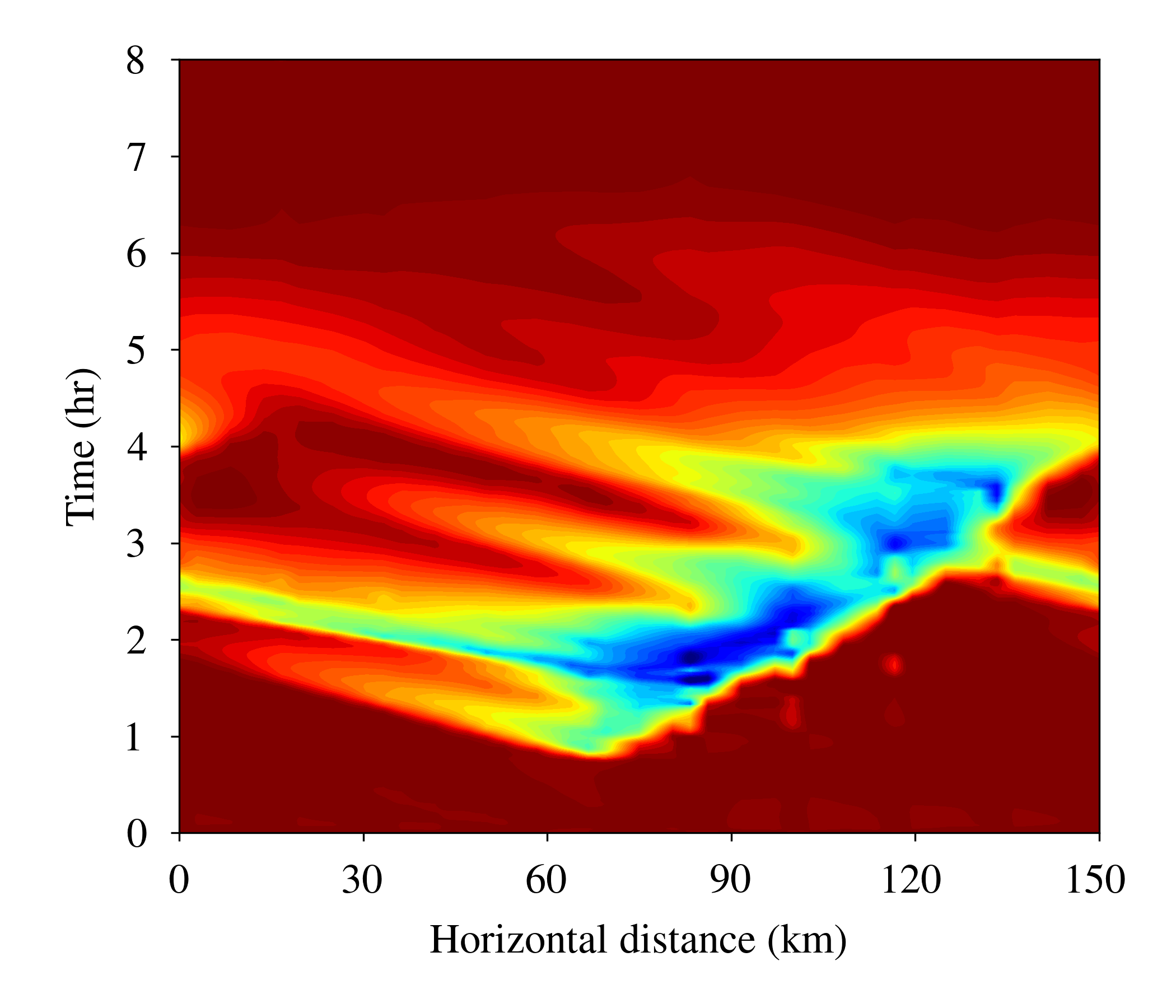}
    \caption{MMF}	
  \end{subfigure}
  \begin{subfigure}[b]{0.32\textwidth}
    \includegraphics[width=\textwidth]{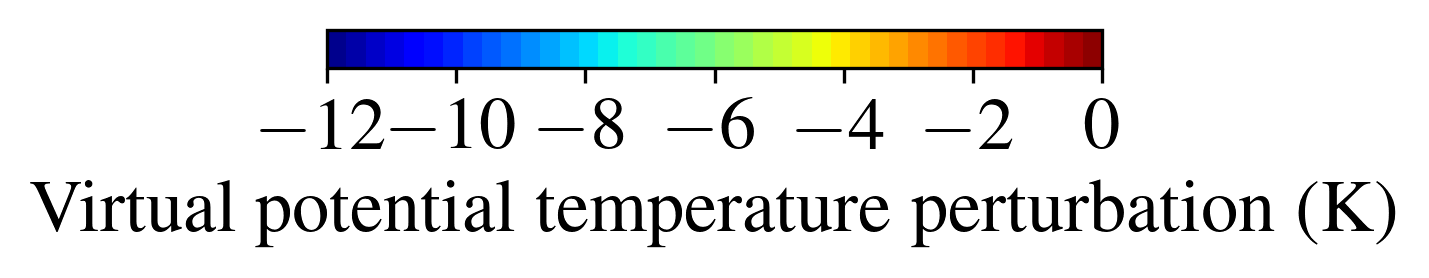}
  \end{subfigure}
  \caption{Contours of the virtual potential temperature perturbations at the surface for the squall line.}\label{fig:squall:theta}
\end{figure}

One of the key parameters in cloud processes is the cold pool, which is a cold pocket of air formed underneath a storm cloud. To identify the cold pool pattern in the simulations, the contours of the potential temperature at the surface are plotted in Figure \ref{fig:squall:theta}. We observe that the MMF captures patterns of cold pools similarly to the fine simulation despite a much larger grid spacing, while the standard coarse simulation produces a wider distribution of cold pools.

% Kinetic Energy

\begin{figure}
  \centering
    \includegraphics[width=0.5\textwidth]{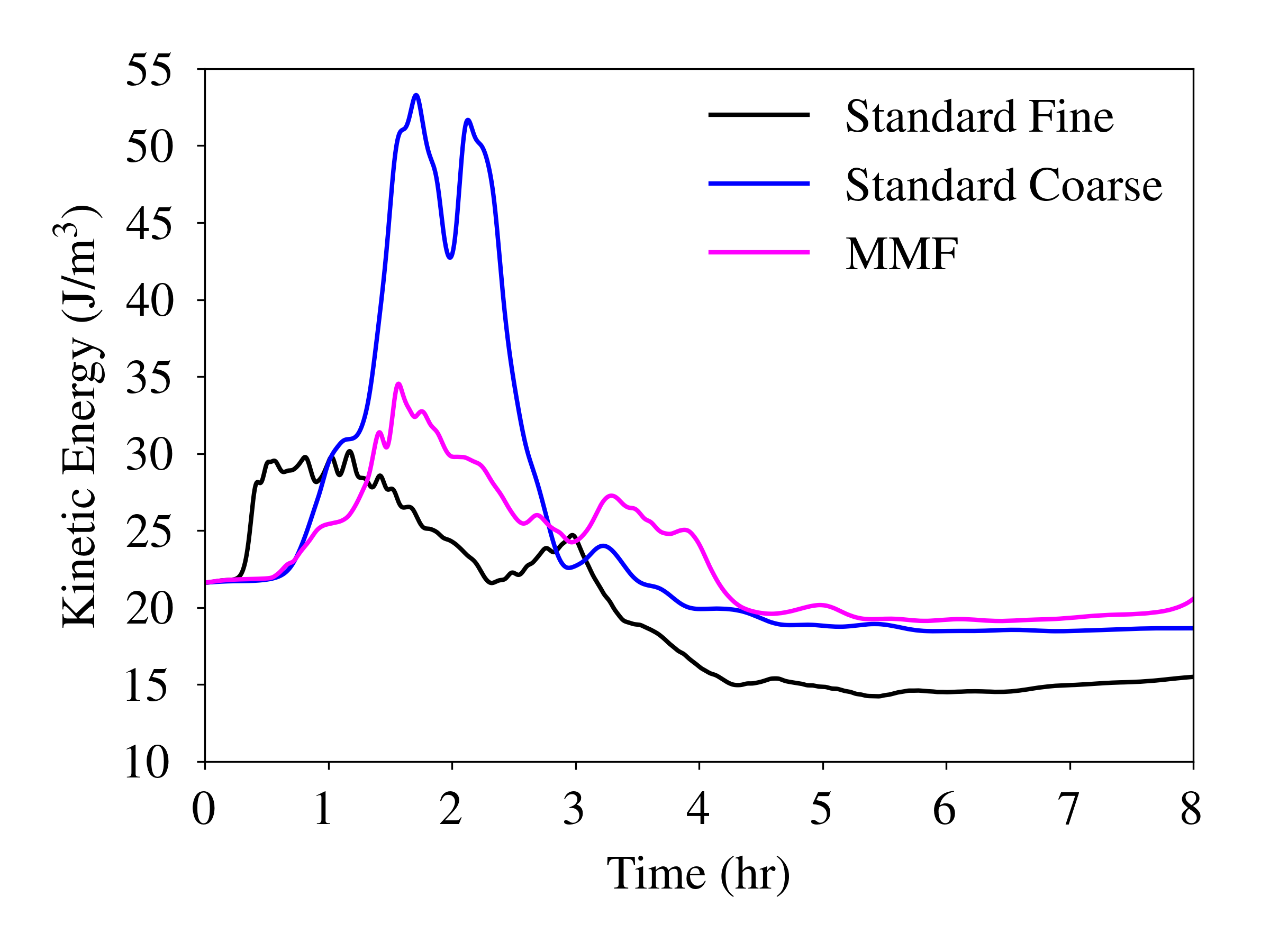}
  \caption{Averaged kinetic energy per volume vs. time for the squall line.}\label{fig:squall:kinetic}
\end{figure}

Figure \ref{fig:squall:kinetic} compares the averaged kinetic energy computed from the three simulations. The averaged kinetic energy per volume is defined as 
\begin{equation}
  \text{KE}=\dfrac{\intO{\tfrac{1}{2}\rho\abs{\ub}^2}} {\intO{}},
\end{equation}
and the integral is computed as in Eq.\ \eqref{eq:integral}. In comparison to the standard fine case, the amplitude of the primary peak in the standard coarse case is much greater. In contrast, the MMF closely matches the kinetic energy profile of the fine simulation in terms of amplitude and pattern.

\begin{figure}
  \centering
  \begin{subfigure}[t]{0.36\textwidth}
    \includegraphics[width=\textwidth]{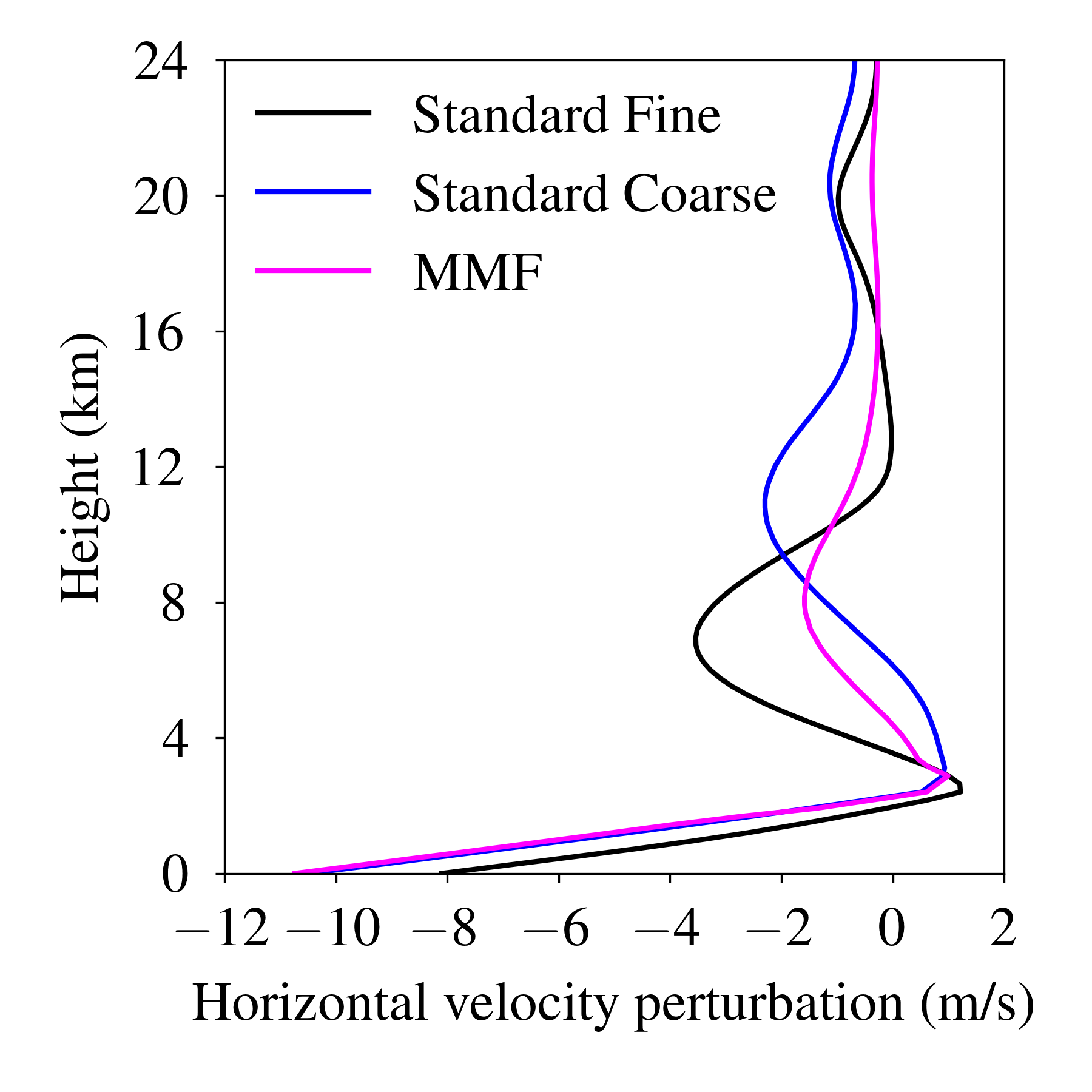}
    \caption{Horizontal velocity}
  \end{subfigure}
  \begin{subfigure}[t]{0.36\textwidth}
    \includegraphics[width=\textwidth]{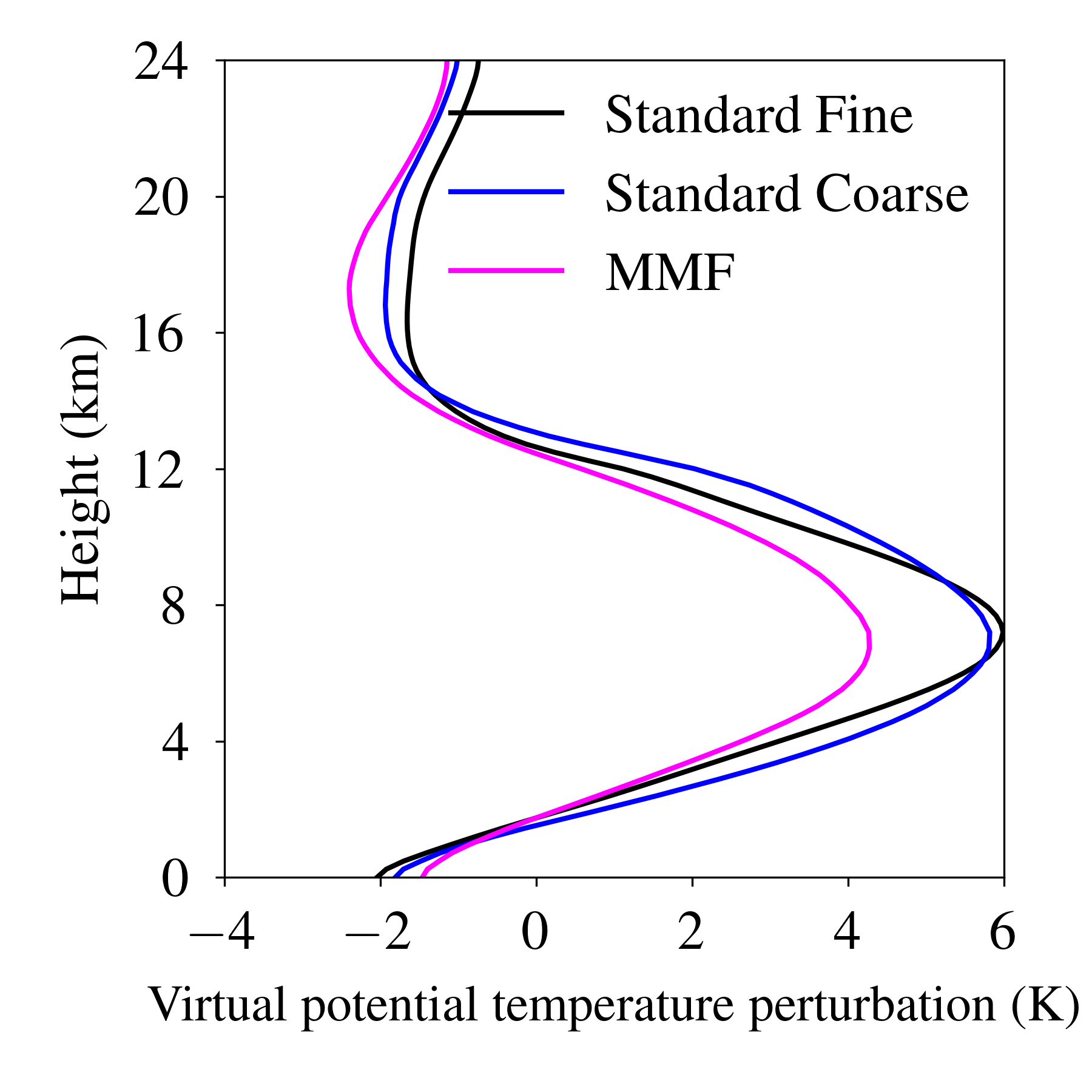}
    \caption{Virtual potential temperature perturbation}
  \end{subfigure}
  \caption{Averaged profiles of horizontal velocity and virtual potential temperature perturbation with respect to height for the squall line.}\label{fig:squall:avg}
\end{figure}

Figure \ref{fig:squall:avg} compares the averaged profiles of the horizontal velocity and virtual potential temperature perturbation along the vertical direction. These profiles are calculated by averaging in time over the whole duration and averaging along the horizontal at each height. We observe that the MMF slightly improves the results for the horizontal velocity perturbation, while underestimating the virtual potential temperature in the troposphere. The reasons for this behavior would make for interesting future work.

% LSP VS SSP

To verify whether the LSP and SSP simulations are properly coupled, Figure \ref{fig:squall:couple} compares the horizontal velocity and potential temperature calculated along the vertical line at the center of the LSP grid and the corresponding SSP. The two curves of the LSP and SSP overlap each other, which demonstrates the accurate imposition of the coupling condition in Eq.\ \eqref{eq:coupling_condition} despite the difference in the vertical grid spacing in the LSP and SSP grids. The same procedure defined by Eqs.\ \eqref{eq:projection:S2L}--\eqref{eq:projection:L2S} is carried out for all the columns, although we only show the coupling at one specific column.

\newcommand\figurewidthtwo{0.3}

\begin{figure}
  \centering
  \begin{subfigure}[b]{\figurewidthtwo\textwidth}
    \includegraphics[width=\textwidth]{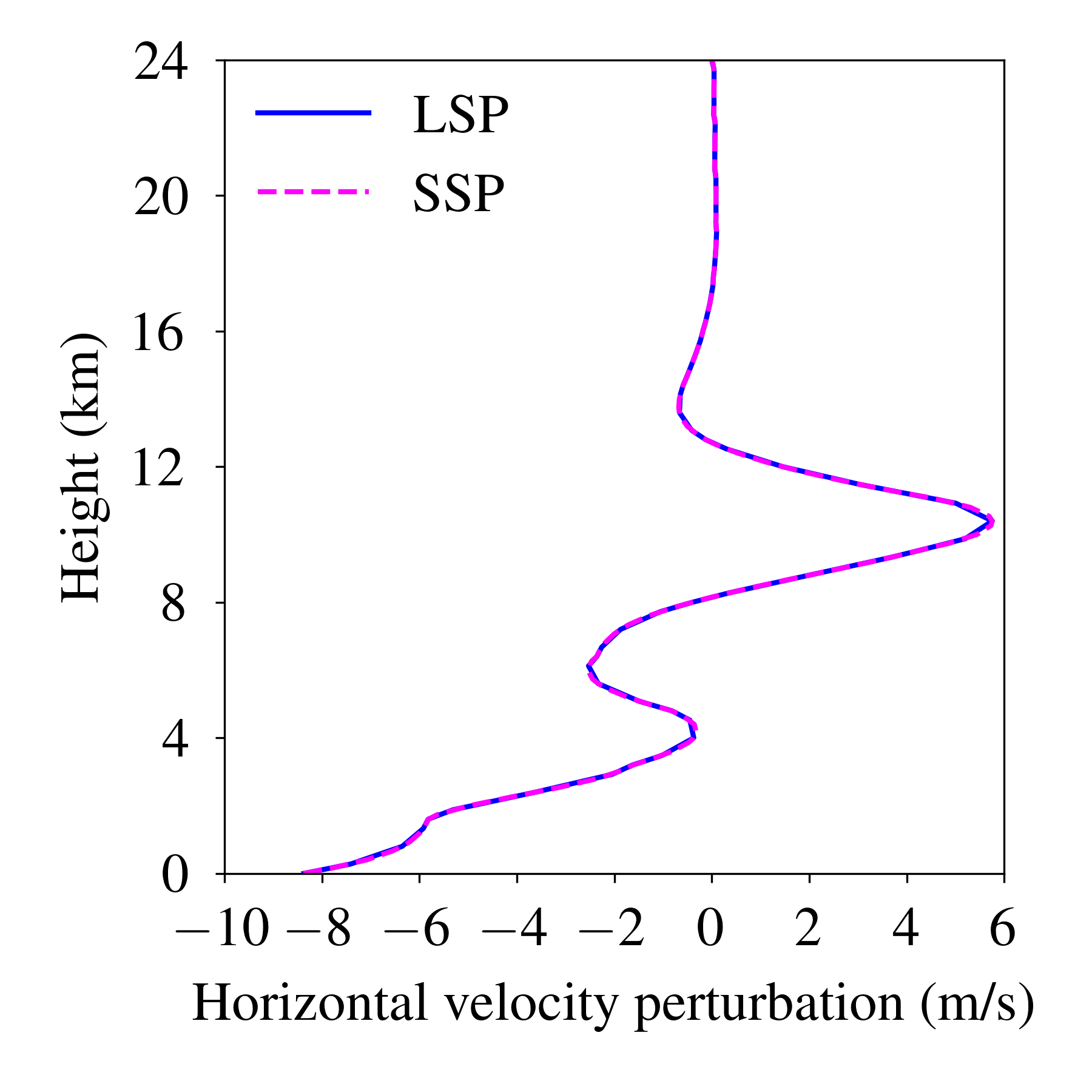}
    \caption{$u$ ($t=3000$)}
  \end{subfigure}
  \begin{subfigure}[b]{\figurewidthtwo\textwidth}
    \includegraphics[width=\textwidth]{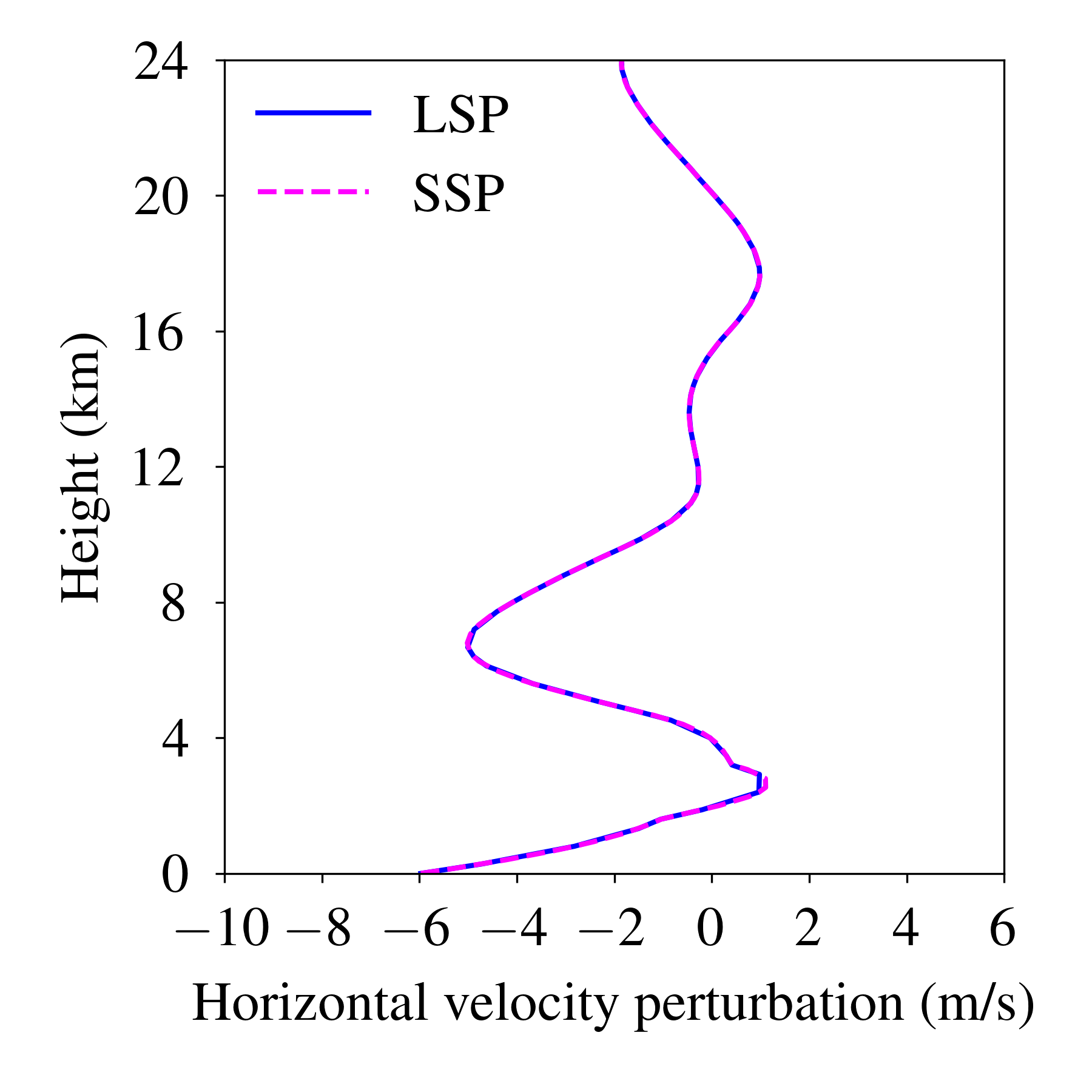}
    \caption{$u$ ($t=9000$)}
  \end{subfigure}
  \\  
  \begin{subfigure}[b]{\figurewidthtwo\textwidth}
    \includegraphics[width=\textwidth]{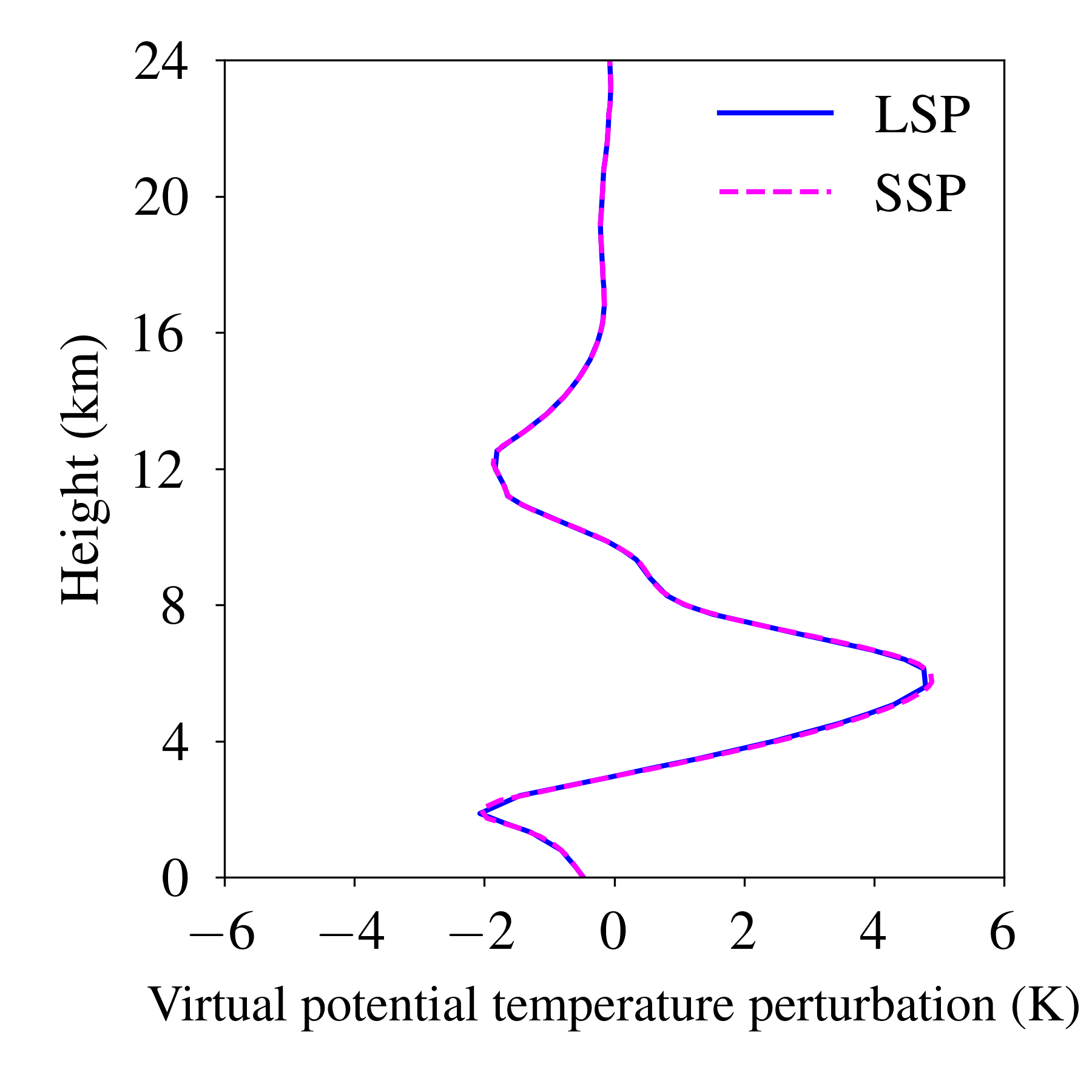}
    \caption{$\theta_v'$ ($t=3000$)}
  \end{subfigure}
  \begin{subfigure}[b]{\figurewidthtwo\textwidth}
    \includegraphics[width=\textwidth]{figures/SquallLine_LSP_vs_SSP/Figure_SquallLine_avg_theta_center_T3000.png}
    \caption{$\theta_v'$ ($t=9000$)}
  \end{subfigure}  
  \caption{Comparisons between the LSP and SSP variables for horizontal velocity and virtual potential temperature perturbation along the vertical line at the domain center at $t=3000$ and 9000 seconds for the squall line. (The SSP variables are horizontally averaged.)}\label{fig:squall:couple}
\end{figure}

%% file: sections/Supercell.tex
\subsection{3D Supercell}

This test case is a benchmark problem that investigates the three-dimensional evolution of a supercell storm \cite{tissaoui2023non}. The dimensions of the domain are 150 km and 100 km in the horizontal directions, and 24 km in height. The $x$ axis is in the streamwise direction, the $y$ axis is in the crosswind direction, and the $z$ axis is in the vertical direction. The impermeable boundary condition is applied at the bottom surface, and periodic boundary conditions are applied at the lateral boundaries in both the $x$ and $y$ directions. At the model top, we apply the implicit sponge layer with the same thickness and coefficient values used in Sec.\ \ref{section:squall}.

For initializing the problem, a thermal perturbation defined by Eqs.\ \eqref{eq:squal:theta_pert} and \eqref{eq:supercell:r} is added on top of the reference field constructed based on the atmosphere sounding from \cite{tissaoui2023non}, similarly to the squall line test case. The radial distance in the elliptical thermal bubble is defined as
\begin{equation}
  \label{eq:supercell:r}
  r = \sqrt{ \qty(\frac{x-x_c}{r_x})^2 + \qty(\frac{y-y_c}{r_y})^2 + \qty(\frac{z-z_c}{r_z})^2 },
\end{equation}
where the coordinates of the center are $(x_c,y_c,z_c)=(75,50,2)$ km, and the semi-major axes are $r_x=r_y=10$ km and $r_z=2$ km, respectively. The SSP models are initialized with the initial state values along the LSP columns. Likewise in the squall line case, we add a random perturbation to the initial virtual potential temperature with an amplitude of 0.3K, as defined in Eq.\ \eqref{eq:rand_pert}.

We compare the standard coarse and MMF simulations with a uniform resolution of 2.5 km by 500 m along the horizontal and vertical directions, respectively. The standard fine grid simulation, with a resolution of 500 m in all directions, serves as the benchmark for this comparison. The width and height of the SSP domains in the MMF are 8 km and 24 km, respectively. The 2D SSP domains are aligned with the background wind direction. The order of the basis functions is 4 in all directions. Numerical diffusion with a magnitude of $\nu=200$ m$^2$/s is added for the standard models and the LSP and SSP models in the MMF. The Boyd-Vandeven filter \cite{boyd1996} with strength of 4\% is applied to the standard models and the LSP model in the MMF. The semi-implicit ARK2 method \cite{giraldo2013implicit} is used to integrate the equations in time until $t=9600$ seconds. The time-step size is $\Delta t=0.5$ seconds for the standard fine grid simulation, and $\Delta t=2.0$ seconds for the standard coarse grid simulation. In the LSP model of the MMF simulation, the time-step size is $\Delta t=2.0$ second, and each time step involves 4 sub-steps for updating the SSP variables.

\begin{figure}
  \centering
  \begin{subfigure}[b]{\figurewidth\textwidth}
    \includegraphics[width=\textwidth]{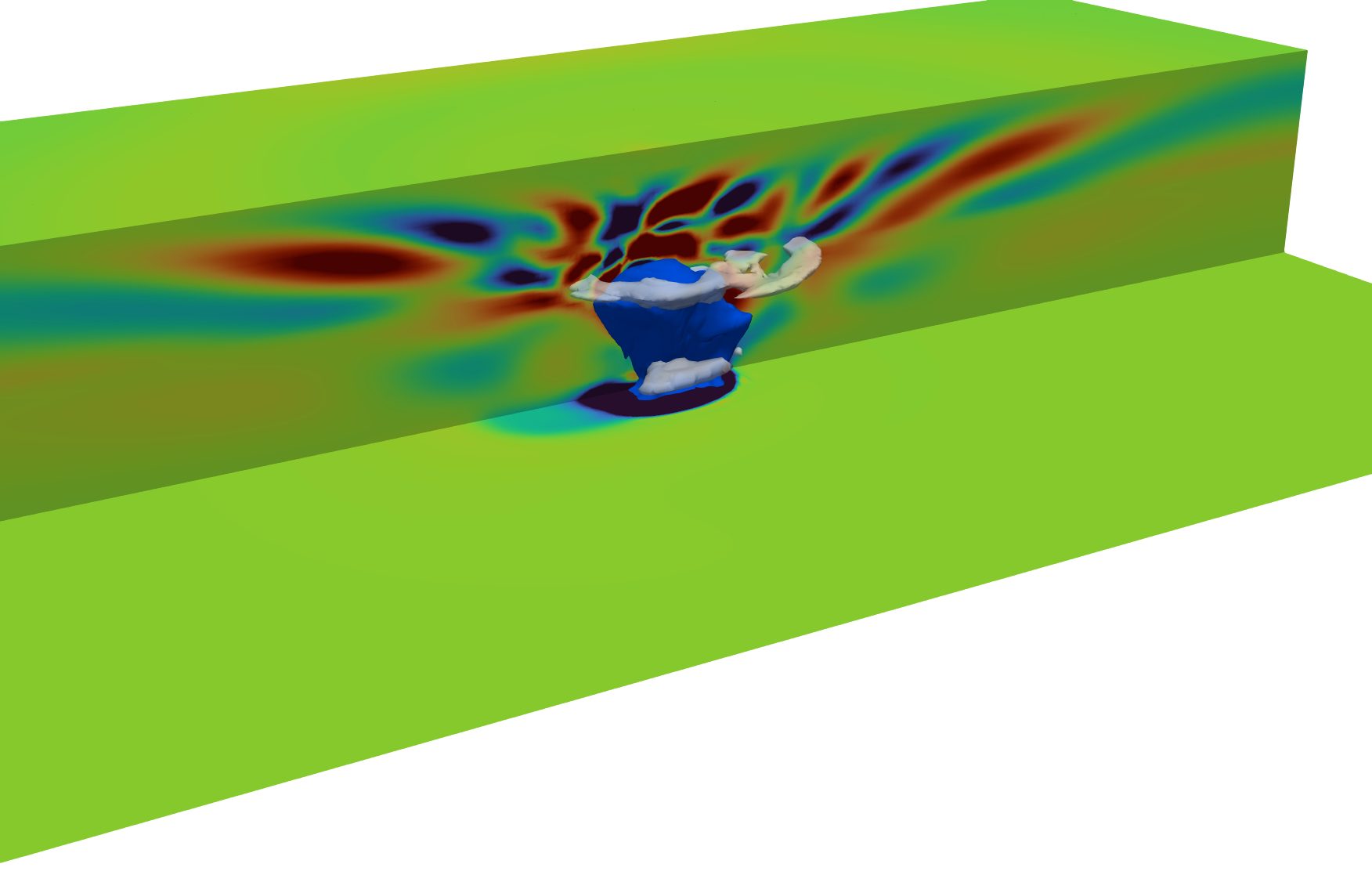}
    \caption{SF ($t=2400$)}
  \end{subfigure}
  \begin{subfigure}[b]{\figurewidth\textwidth}
    \includegraphics[width=\textwidth]{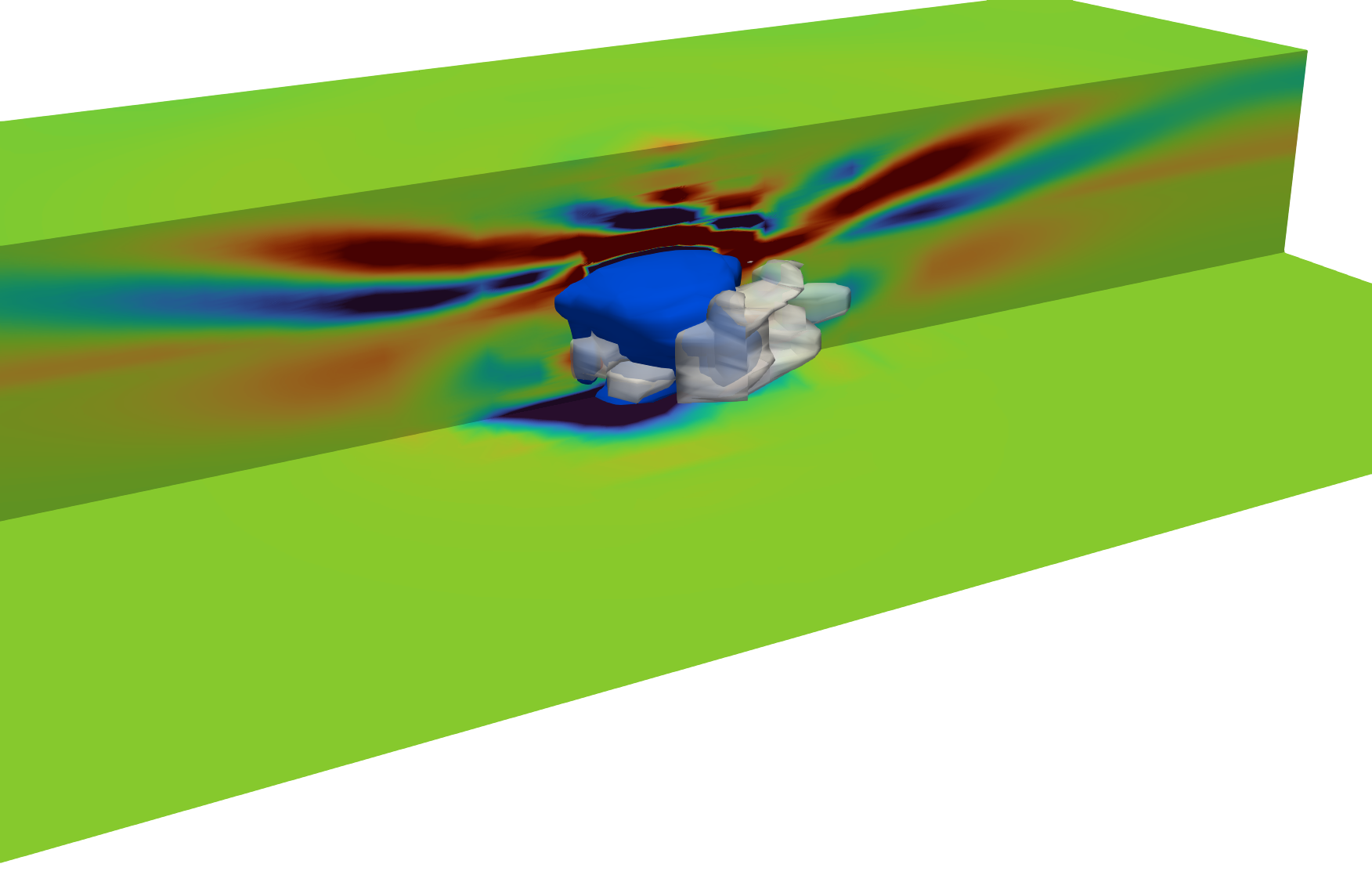}
    \caption{SC ($t=2400$)}
  \end{subfigure}
  \begin{subfigure}[b]{\figurewidth\textwidth}
    \includegraphics[width=\textwidth]{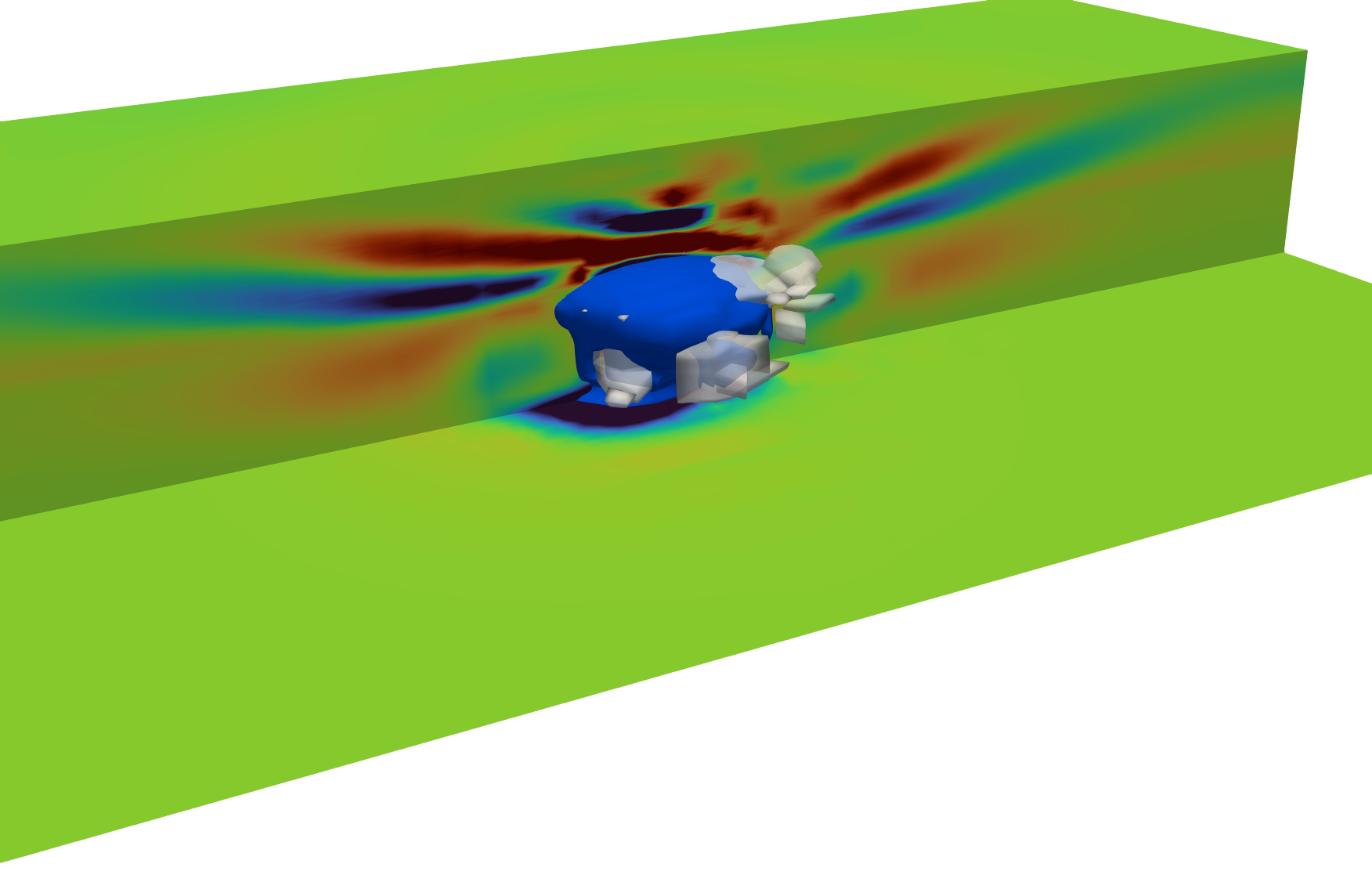}
    \caption{MMF ($t=2400$)}
  \end{subfigure}
  \\
  \begin{subfigure}[b]{\figurewidth\textwidth}
    \includegraphics[width=\textwidth]{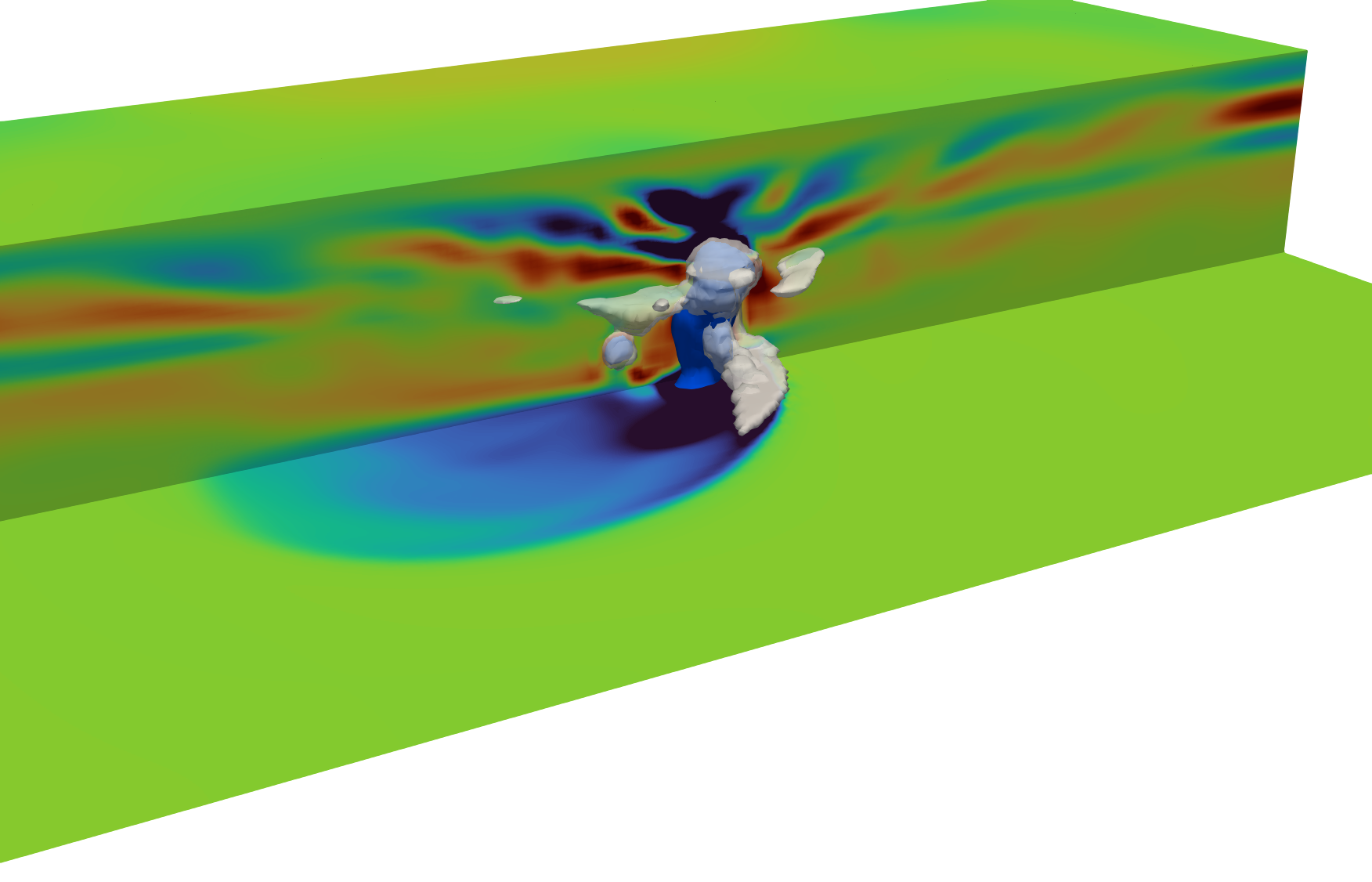}
    \caption{SF ($t=4800$)}
  \end{subfigure}
  \begin{subfigure}[b]{\figurewidth\textwidth}
    \includegraphics[width=\textwidth]{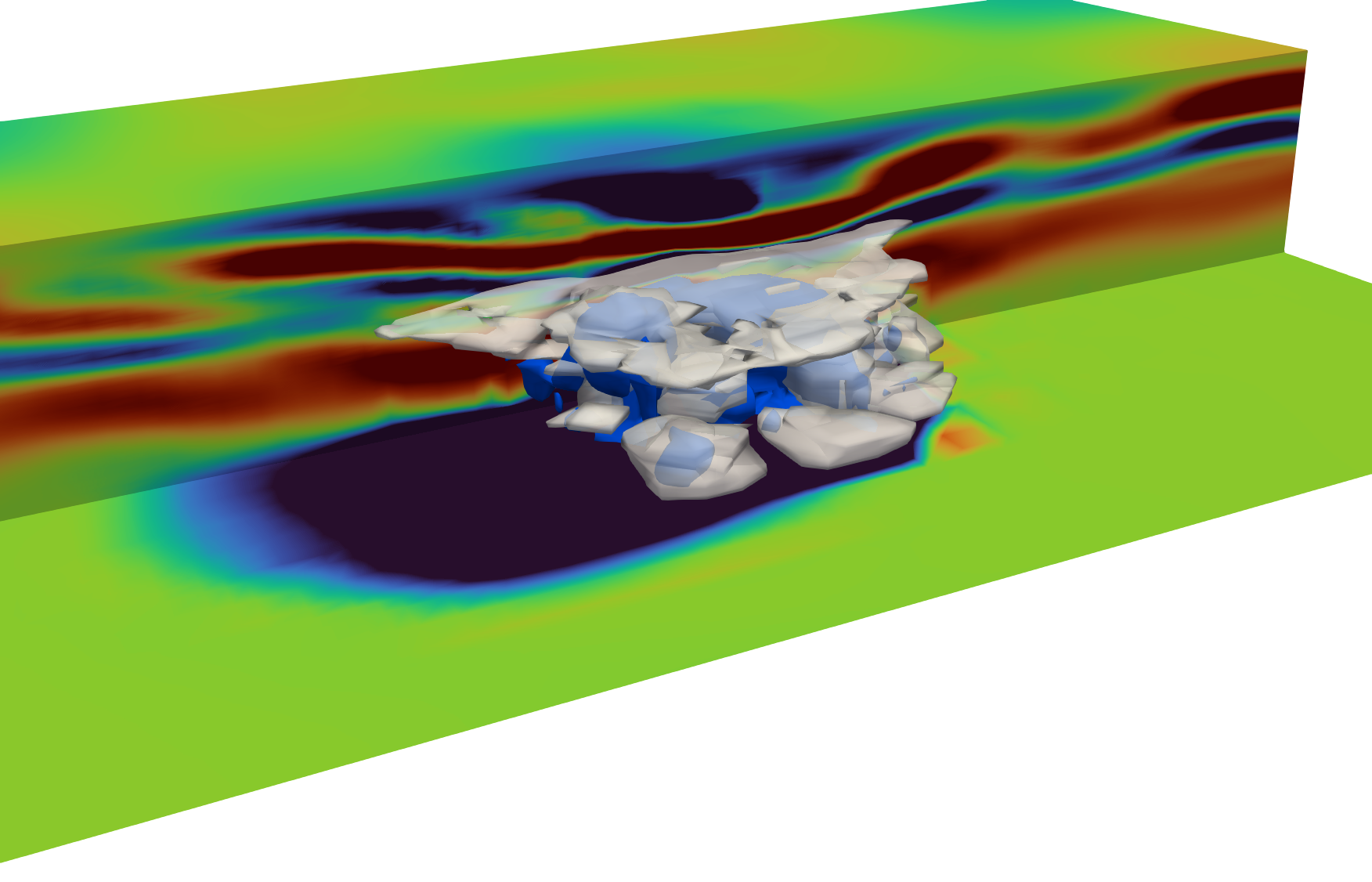}
    \caption{SC ($t=4800$)}
  \end{subfigure}
  \begin{subfigure}[b]{\figurewidth\textwidth}
    \includegraphics[width=\textwidth]{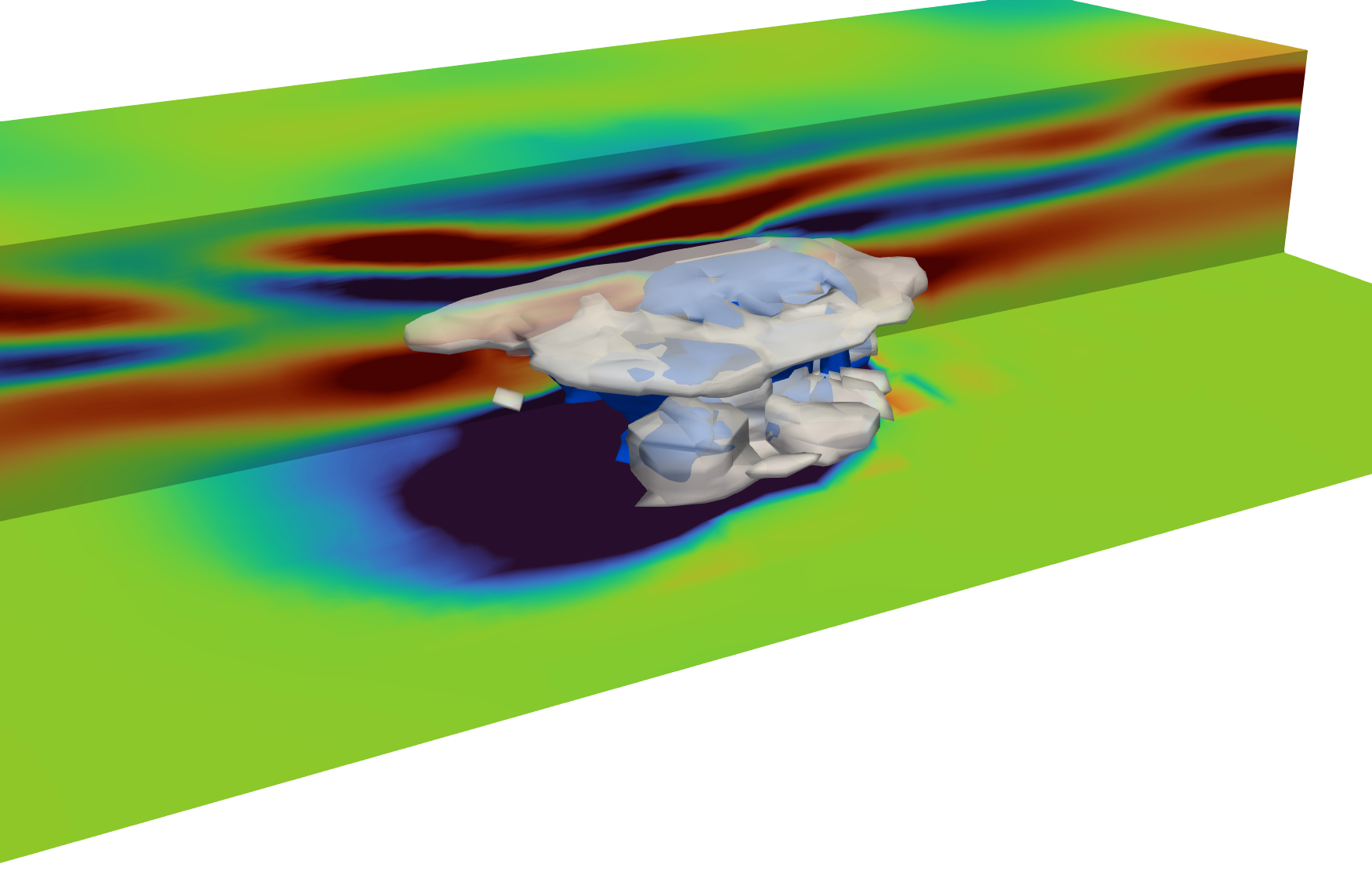}
    \caption{MMF ($t=4800$)}
  \end{subfigure}  
  \\
  \begin{subfigure}[b]{\figurewidth\textwidth}
    \includegraphics[width=\textwidth]{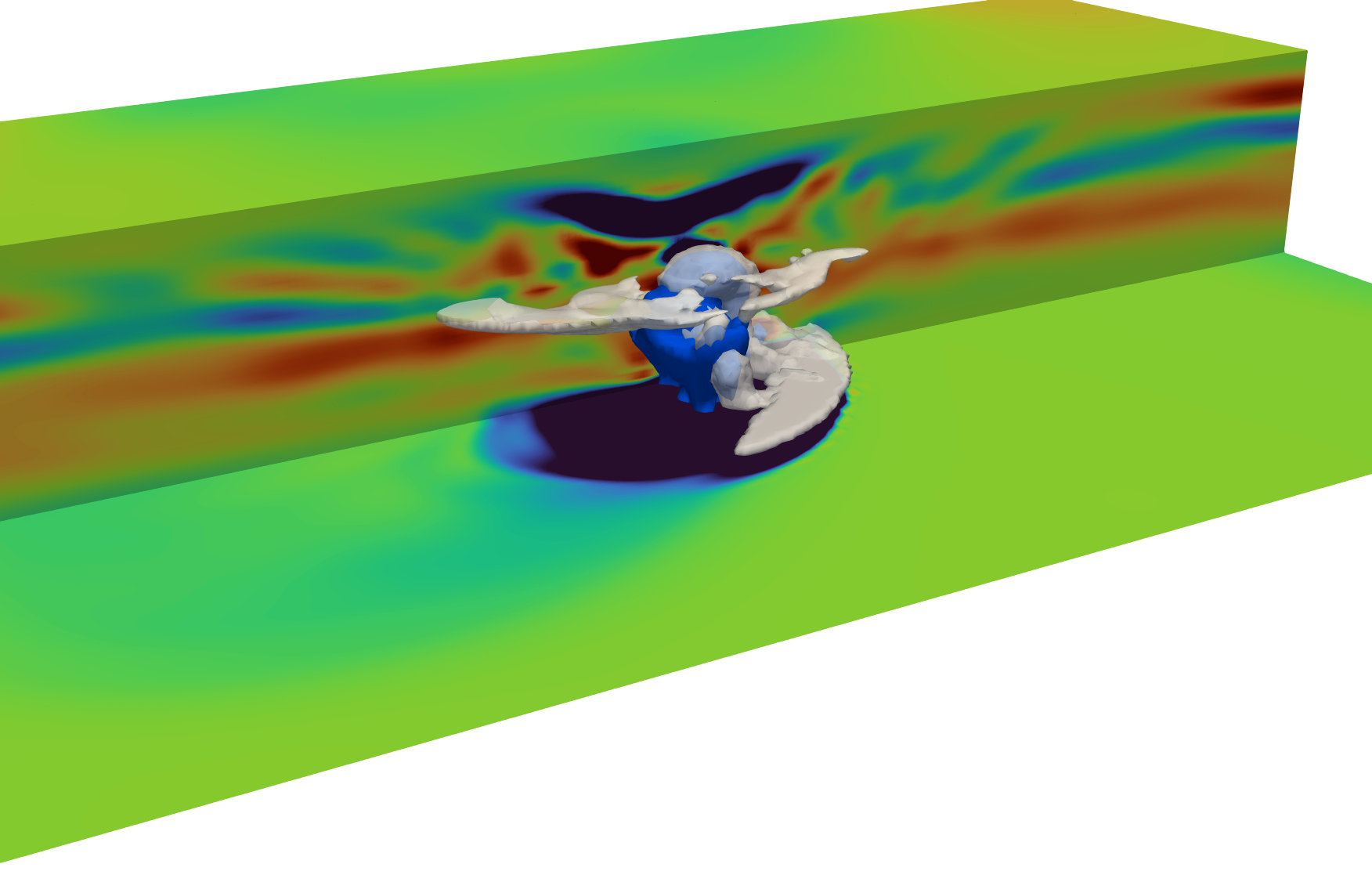}
    \caption{SF ($t=7200$)}
  \end{subfigure}
  \begin{subfigure}[b]{\figurewidth\textwidth}
    \includegraphics[width=\textwidth]{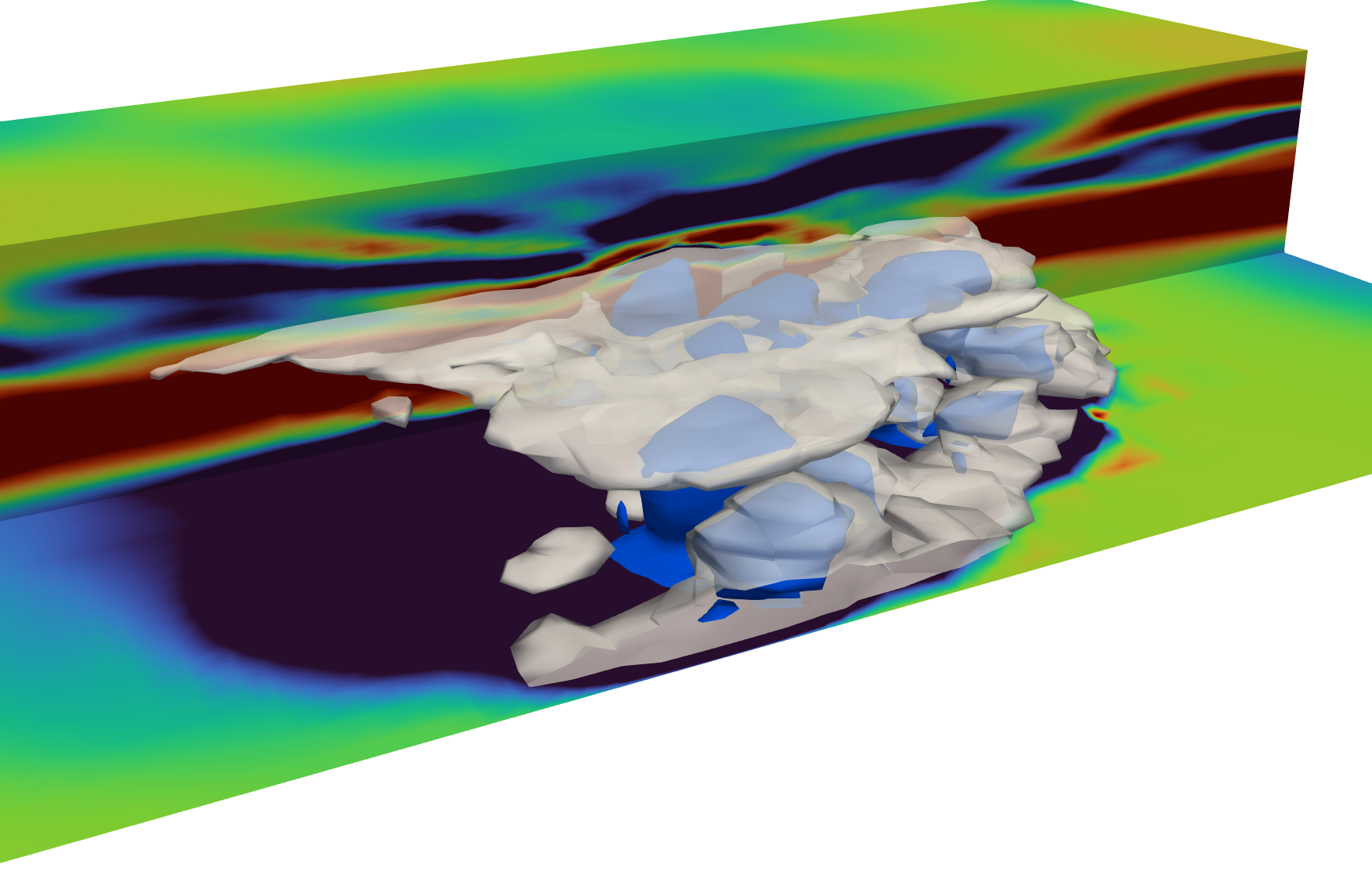}
    \caption{SC ($t=7200$)}
  \end{subfigure}
  \begin{subfigure}[b]{\figurewidth\textwidth}
    \includegraphics[width=\textwidth]{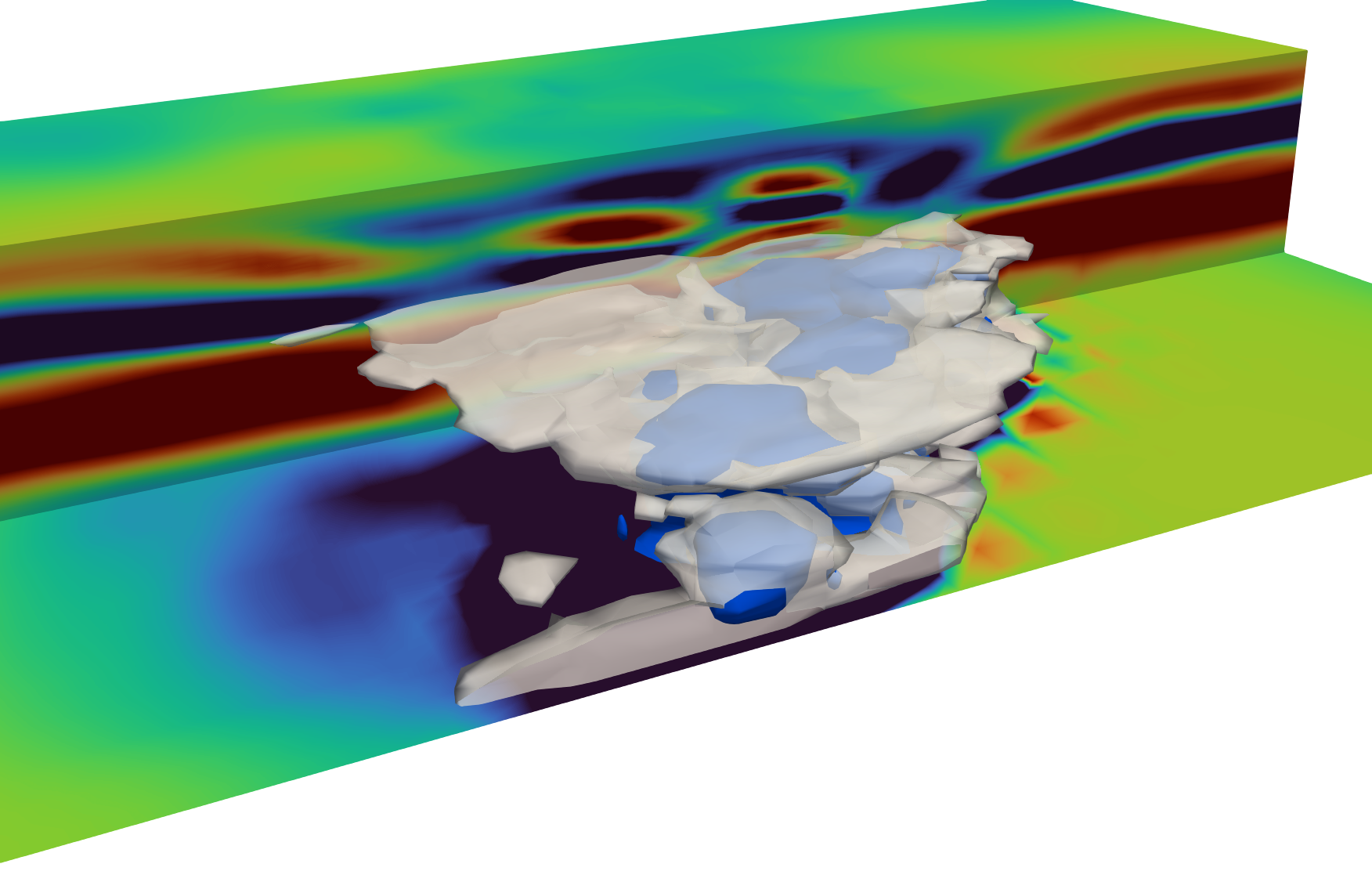}
    \caption{MMF ($t=7200$)}
  \end{subfigure} 
  \\
  \begin{subfigure}[b]{\figurewidth\textwidth}
    \includegraphics[width=\textwidth]{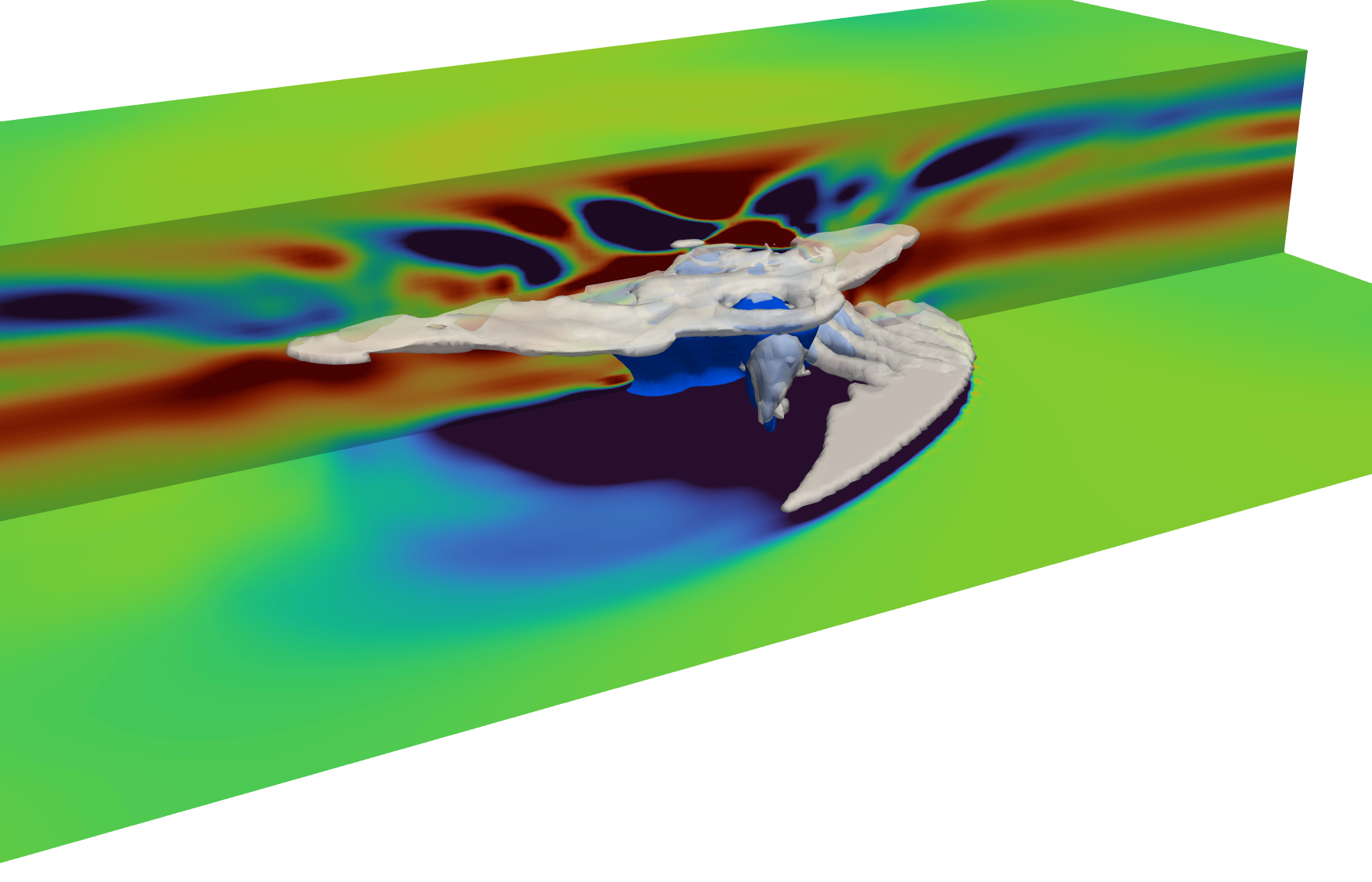}
    \caption{SF ($t=9600$)}
  \end{subfigure}
  \begin{subfigure}[b]{\figurewidth\textwidth}
    \includegraphics[width=\textwidth]{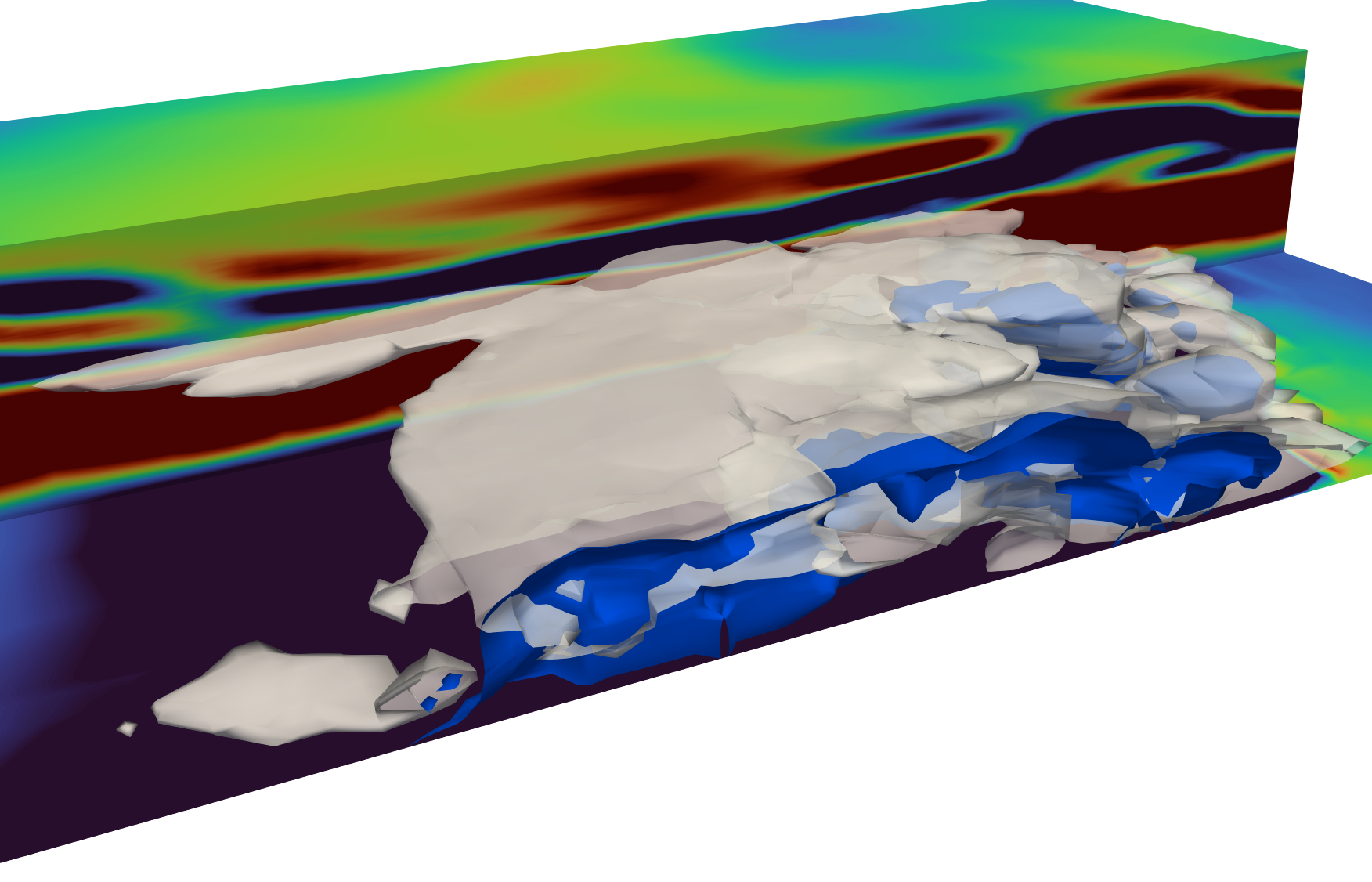}
    \caption{SC ($t=9600$)}
  \end{subfigure}
  \begin{subfigure}[b]{\figurewidth\textwidth}
    \includegraphics[width=\textwidth]{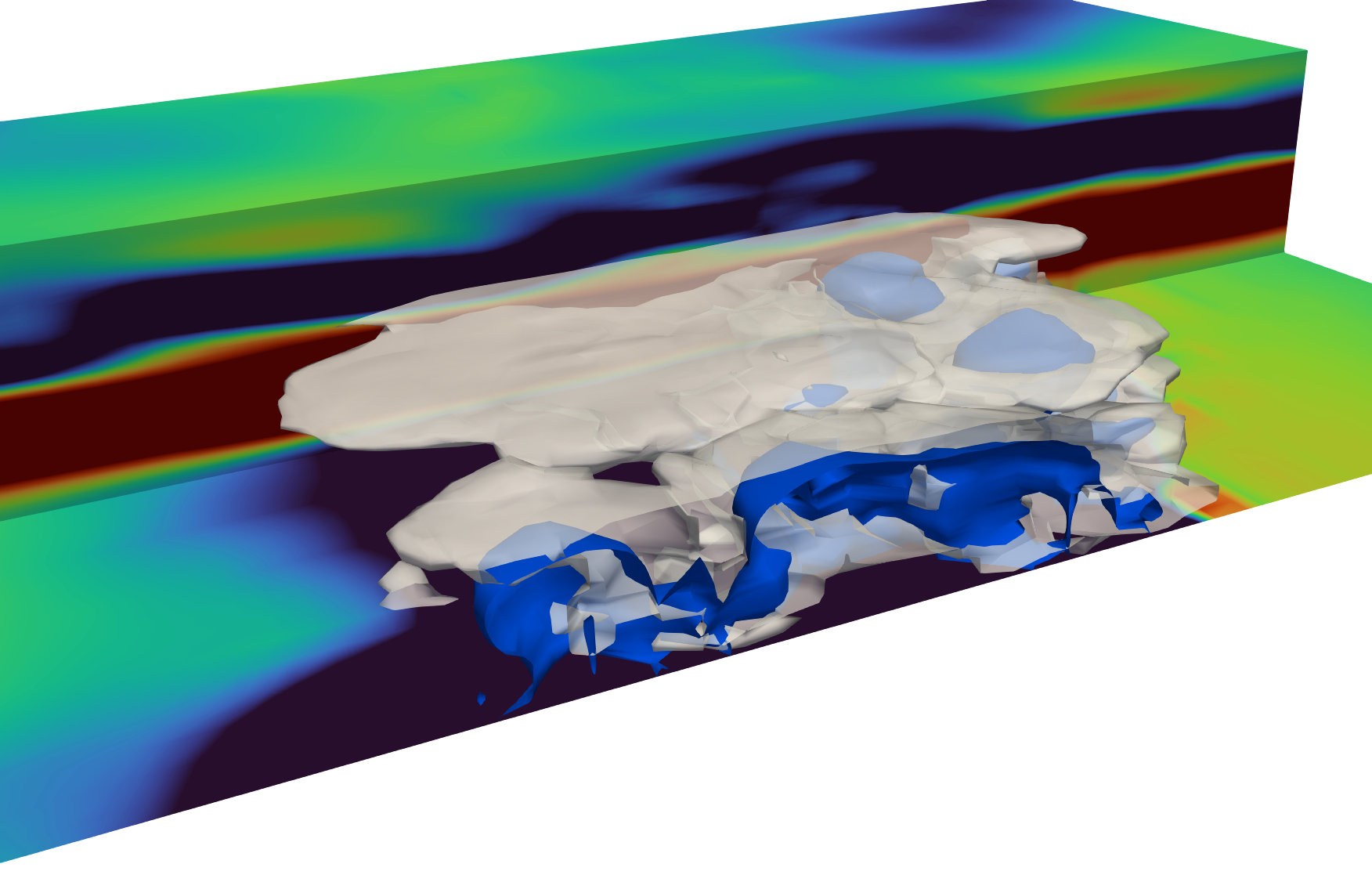}
    \caption{MMF ($t=9600$)}
  \end{subfigure} 
  \\ 
  \begin{subfigure}[b]{0.32\textwidth}
    \includegraphics[width=\textwidth]{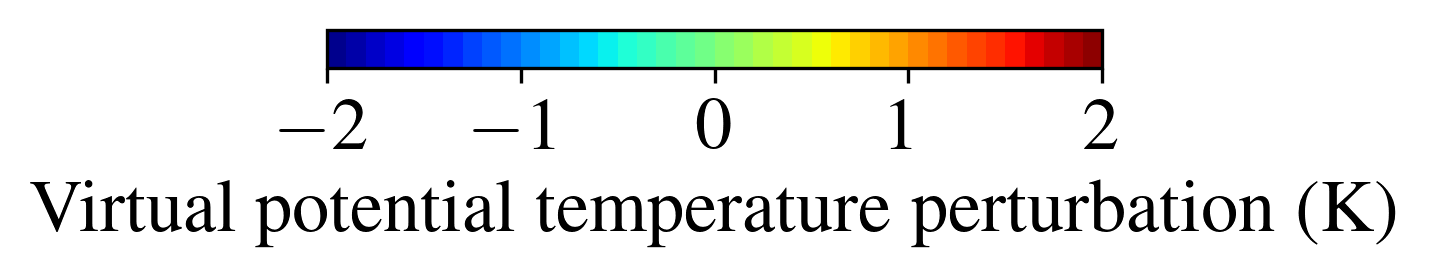}
  \end{subfigure}
  \caption{Instantaneous virtual potential temperature perturbation field and the contours of the cloud and rain concentrations at $t$=2400, 4800, 7200, and 9600 seconds computed from the standard fine (SF), standard coarse (SC), and MMF simulations for the supercell. The gray surface represents the contour of the cloud mixing ratio at $q_c=10^{-5}$, and the light-blue surface represents the contours of the rain mixing ratio at $q_r=10^{-3}$.}\label{fig:supercell:snapshot}
\end{figure}

The supercell cloud is initiated by a thermal perturbation and rises developing into cumulonimbus clouds, eventually merging into the wider cloud around the inversion layer. The position of the rain concentration follows the location of the convective towers, with the rain falling below the cloud towers. Figure \ref{fig:supercell:snapshot} displays the instantaneous virtual potential temperature field at various time steps. Since the evolution of the supercell is symmetric with respect to the mid $x$-$z$ plane, we only show the storm patterns in one half of the domain. Like the 2D squall line test, the standard coarse grid produces more clouds than the benchmark fine grid simulation. It is conjectured that these spurious clouds in the standard coarse simulation are attributed to the large grid size that degrades the accuracy in the parametrization of the moist microphysics. The MMF simulation produces improved results in terms of the size of the clouds and the patterns of the temperature field at the ground. Even though the LSP grid of the MMF case has the same 2 km resolution, the SSP models resolve the cloud motions locally and parameterizes them according to the appropriate grid size. As a result, the quality of the solution within the LSP domain is enhanced in the MMF case.

\begin{figure}
  \centering
  \begin{subfigure}[b]{\figurewidth\textwidth}
    \includegraphics[width=\textwidth]{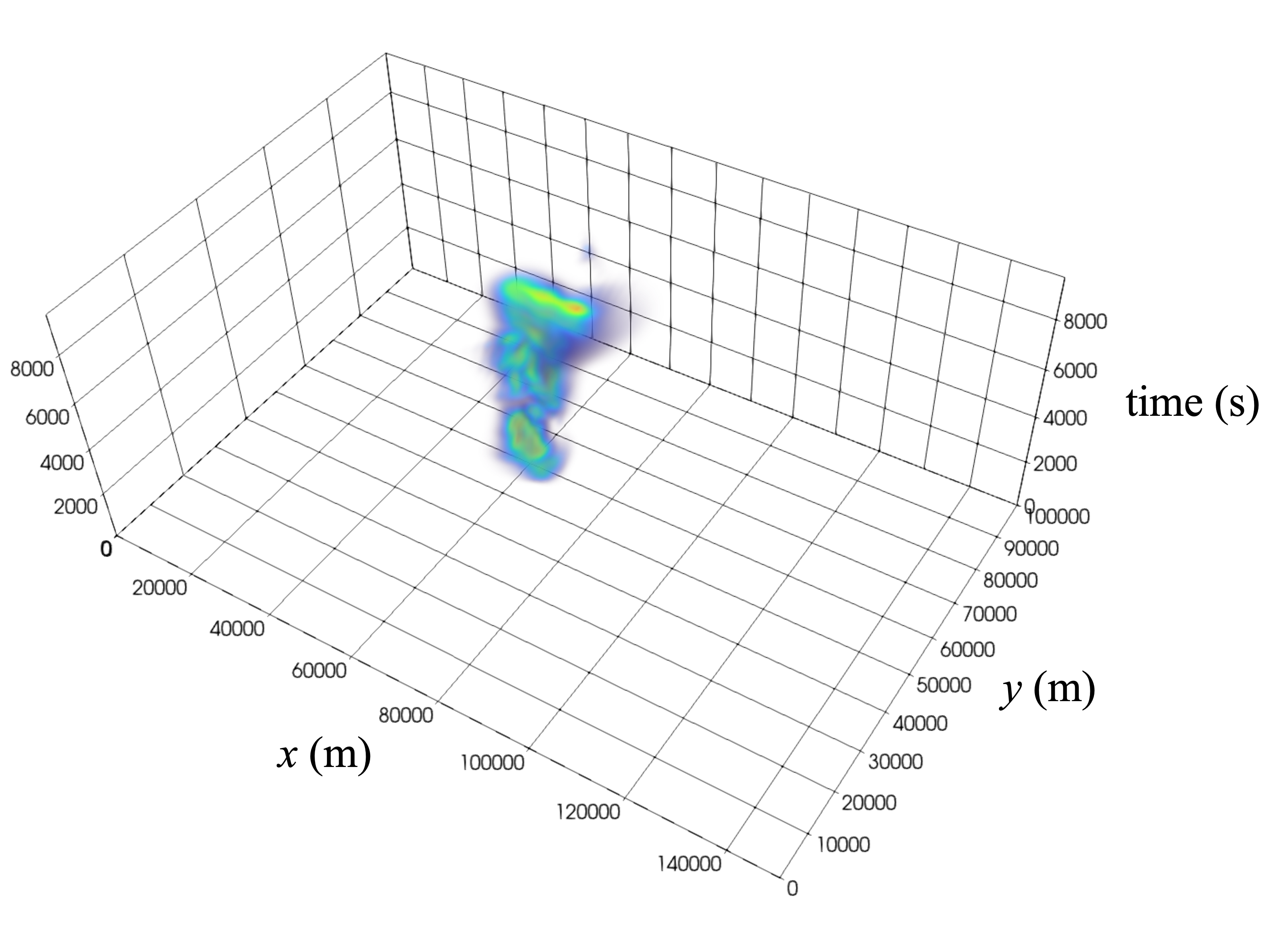}
    \caption{Standard fine}
  \end{subfigure}
  \begin{subfigure}[b]{\figurewidth\textwidth}
    \includegraphics[width=\textwidth]{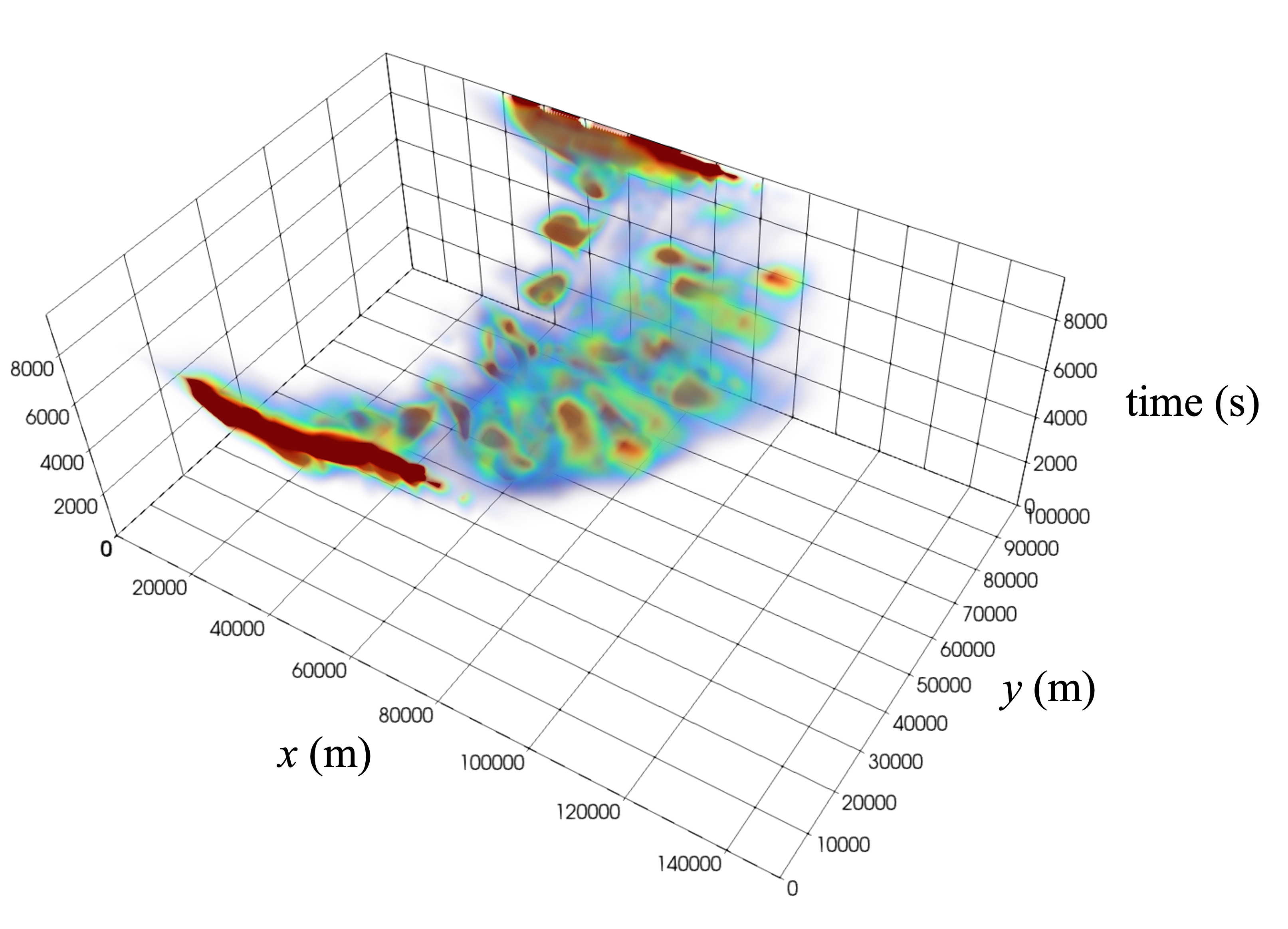}
    \caption{Standard coarse}	
  \end{subfigure}
  \begin{subfigure}[b]{\figurewidth\textwidth}
    \includegraphics[width=\textwidth]{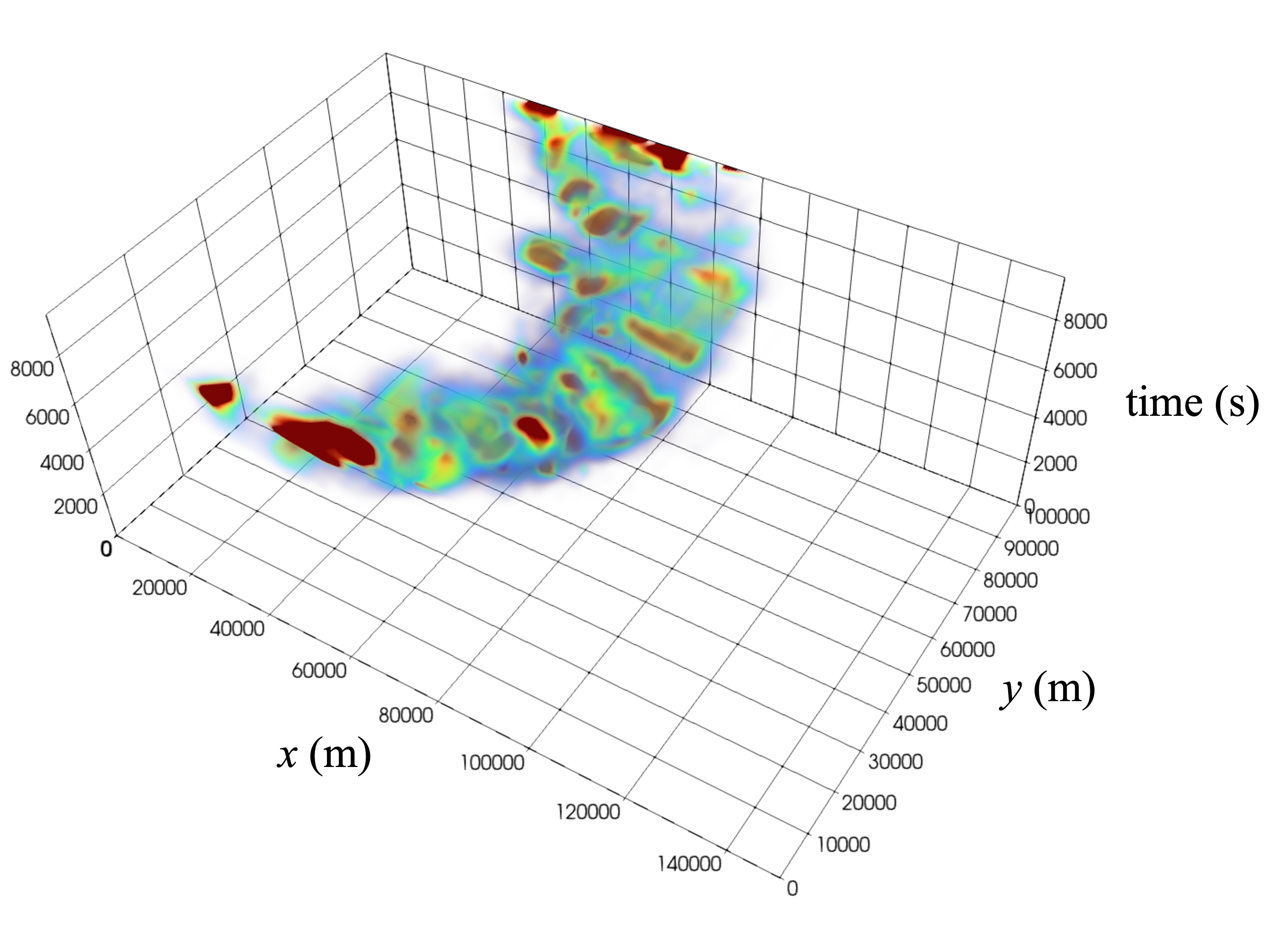}
    \caption{MMF}	
  \end{subfigure}
  \begin{subfigure}[b]{0.25\textwidth}
    \includegraphics[width=\textwidth]{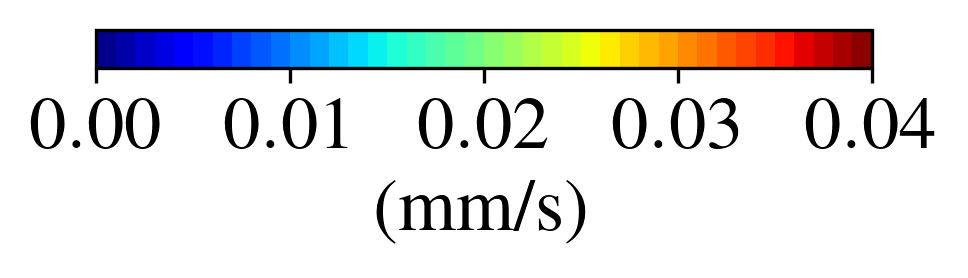}
  \end{subfigure}
  \caption{Volume rendering of surface precipitation for the supercell. The horizontal axes are the distances in the $x$ and $y$ directions (unit:meters), and the vertical axis is time (unit: seconds).}\label{fig:supercell:rain}
\end{figure}

Figure \ref{fig:supercell:rain} displays contours of the surface precipitation patterns in $x$-$y$-$t$ space. While the standard fine case shows a narrow spread of surface precipitation, the standard coarse case shows a wider spread in both the streamwise and crosswind directions and the rain reaches the lateral boundaries in the $y$ direction at around $t=9600$ seconds. The surface rain in the MMF simulation also reaches the lateral boundaries, but it exhibits a narrower distribution, particularly along the SSP grid direction, compared to the standard coarse case. These results verify that the MMF improves the representation of the moist microphysics and enhances accuracy in predicting surface precipitation.

\begin{figure}
  \centering
  \begin{subfigure}[b]{\figurewidth\textwidth}
    \includegraphics[width=\textwidth]{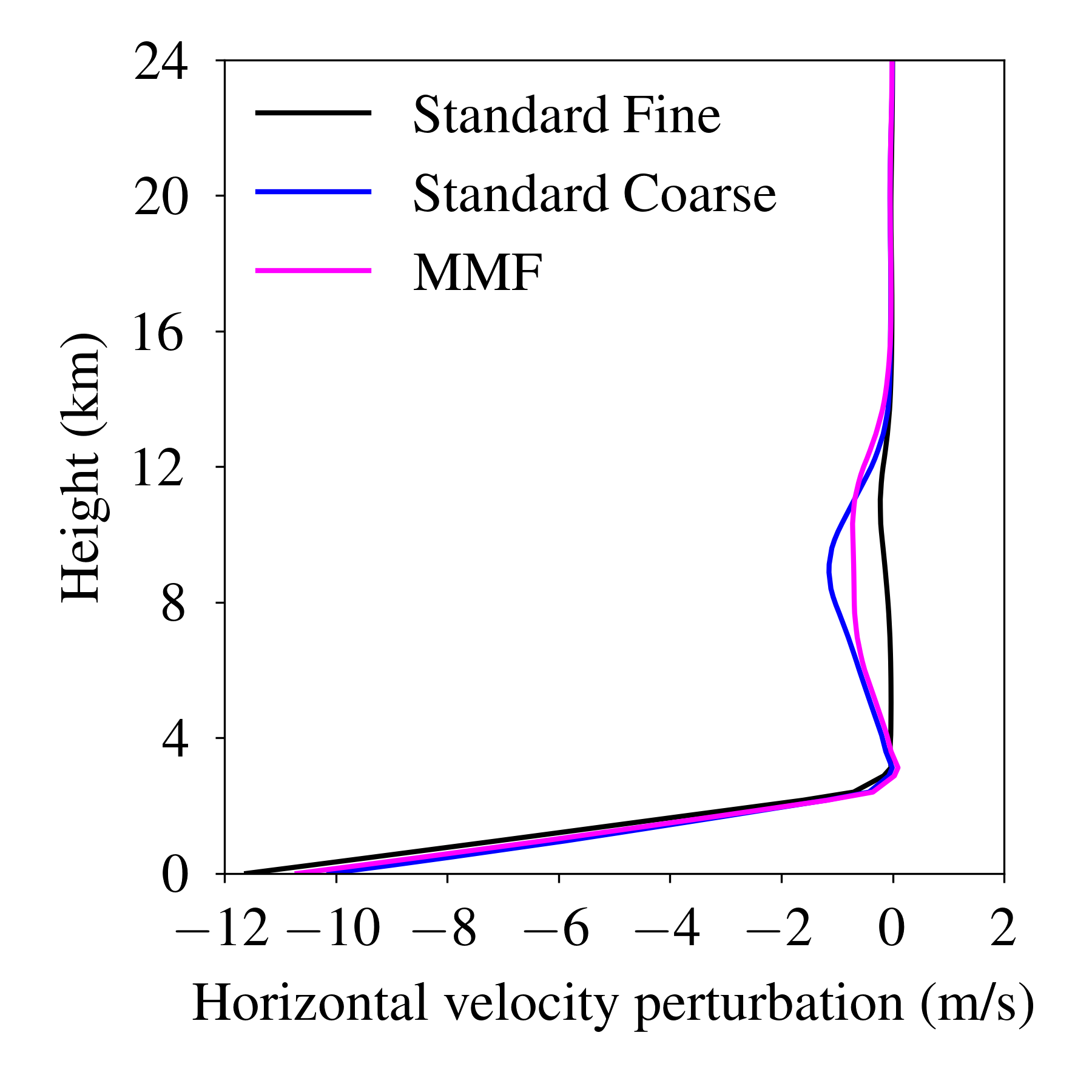}
    \caption{Horizontal velocity perturbation}
  \end{subfigure}
  \begin{subfigure}[b]{\figurewidth\textwidth}
    \includegraphics[width=\textwidth]{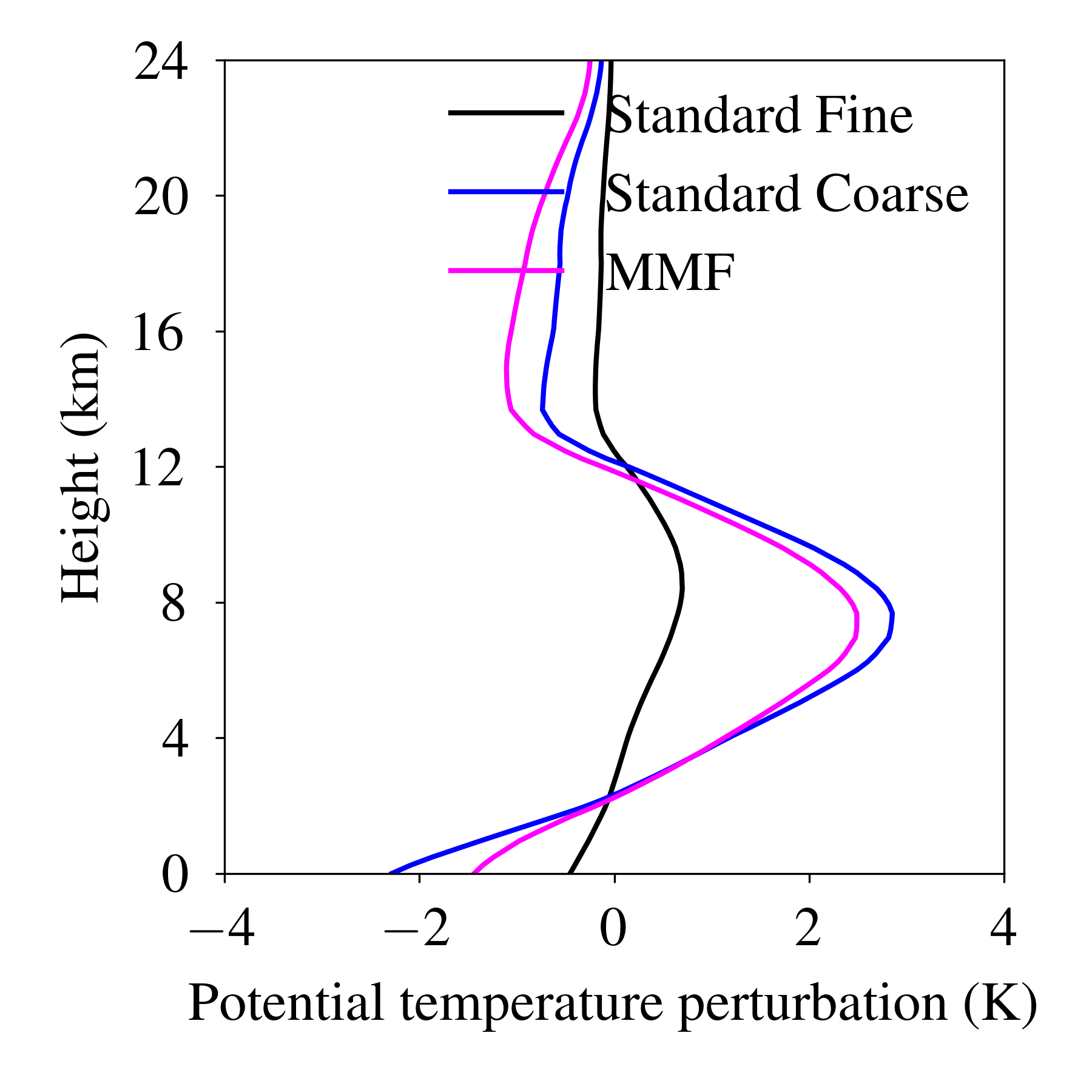}
    \caption{Virtual potential temperature perturbation}
  \end{subfigure}
  \caption{Averaged profiles of horizontal velocity and virtual potential temperature perturbation with respect to height for the supercell. The results are averaged over time and the horizontal plane.}\label{fig:supercell:avg}
\end{figure}

Figure \ref{fig:supercell:avg} compares the profiles of potential temperature perturbation and horizontal velocity obtained from the standard fine and coarse, and the MMF simulations. The numerical results are averaged over time and horizontal plane to condense them into 1D profiles. These plots show that the MMF improves (if only slightly) the results of horizontal velocity and virtual potential temperature that are parameterized by the SSPs.

%% file: sections/Conclusions.tex
\section{Conclusions}\label{sec:conclusion}

We presented a multiscale modeling framework to resolve distinct scales in large-scale and small-scale processes in the moist atmosphere. This approach aims to inject higher-dimensional feedback into a lower-dimension (coarse) model. In the MMF, we utilized the compressible Navier-Stokes equations and an element-based Galerkin method for both large-scale and small-scale models, which, to our knowledge, is the first time that a consistent approach and form of the governing equations have been used in MMF modeling. Moreover, our models are built within the same code-base. Each LSP grid column is coupled with a two dimensional SSP model, where non-conforming vertical discretization is allowed. This enables the adjustment of the grid spacing or the order of basis functions in the SSP model. Through a complexity analysis, we demonstrated that the MMF algorithm exhibits higher arithmetic intensity compared to the standard approach which hints at a likely better performance on exascale computing hardware. In the squall line and supercell test cases, the MMF results showed enhanced representation of moist dynamics for cloud processes compared to the coarse model.

This research is our initial endeavor on MMF for the moist atmosphere and admits many potential avenues for extension. One such avenue is to examine ways to reduce the cost of MMF simulations. The MMF algorithm can be optimized by simplifying the SSP computations using lower-order basis functions or a fully explicit time integrator. Developing reduced-order approximations, such as via proper orthogonal decomposition or neural networks \cite{hannah2020initial}, is a viable option to substitute the SSP models. Another avenue is to explore new configurations of MMF, in which 3D SSP models are coupled with each element column of the LSP model instead of the grid column for dimensional consistency. Finally, we are deploying the SSP simulations onto GPUs in order to accelerate time to solution for our MMF simulations. 

% \rc{(FXG: how about the next steps? I propose saying that this is the first attempt at using the new infrastructure, xNUMA, to develop a standard super-parameterization approach. Now that we have done this, we think we understand how to improve the results. E.g., align the x-z planes with the maximum horizontal flow velocity or the element-based MMF approach that you are developing.  Perhaps not here but somewhere else - Intro maybe - we need to state why MMF is useful. E.g., we would like to avoid parameterizations as much as possible because it greatly simplifies simulating subgrid-scale processes without having to tune.  Jim should chime in here to help us frame the challenges faced by operational weather models in the tuning and development of parameterizations. We want to solve the governing equations of motion. Might also say something about how the MMF model can be used to train neural networks. Anthony, does that sound like something we should say?)}